\documentclass{amsart}
\usepackage{graphicx}
\newtheorem{thm}{Theorem}[section]
\newtheorem*{thmnonnum}{Theorem}
\newtheorem*{propnonnum}{Proposition}
\newtheorem*{lemnonnum}{Lemma}
\newtheorem{cor}[thm]{Corollary}
\newtheorem{lem}[thm]{Lemma}
\newtheorem{prop}[thm]{Proposition}
\newtheorem{defi}[thm]{Definition}
\theoremstyle{remark}
\newtheorem{domanda}[thm]{Question}
\newtheorem{rem}[thm]{Remark}
\newenvironment{demo}{\noindent{\bf Proof: }}{\hfill$\Box$\medskip}
\numberwithin{equation}{section}
\newcommand\NN {{\mathbb N}}
\newcommand\ZZ {{\mathbb Z}}

\newcommand\RR {{\mathbb R}}
\newcommand\CC {{\mathbb C}}

\newcommand\PP {{\mathbb P}}

\newcommand\sldz{{\rm SL(d,\ZZ)}}

\newcommand\distanza{{\rm dist }}

\newcommand{\sys}{\texttt{Sys}}

\newcommand\leb{{\rm Leb }}

\newcommand\sign{{\rm sign }}
\newcommand\proiettivo{{\rm Proj }}

\newcommand\spanne{{\rm Span }}
\newcommand\hol {{\rm Hol}}
\newcommand\modphi{{\rm mod \varphi}}

\newcommand\al{\alpha}
\newcommand\be{\beta}
\newcommand\ga{\gamma}

\newcommand\la{\lambda}

\newcommand\te{\theta}

\newcommand\Ga{\Gamma}
\newcommand\De{\Delta}

\newcommand\cA{{\mathcal{A}  }}
\newcommand\cB{{\mathcal{B}  }}
\newcommand\cC{{\mathcal{C}  }}
\newcommand\cD{{\mathcal{D}  }}
\newcommand\cE{{\mathcal{E}  }}
\newcommand\cF{{\mathcal{F}  }}
\newcommand\cG{{\mathcal{G}  }}
\newcommand\cH{{\mathcal{H}  }}
\newcommand\cI{{\mathcal{I}  }}

\newcommand\cN{{\mathcal{N}  }}

\newcommand\cP{{\mathcal{P}  }}

\newcommand\cR{{\mathcal{R}  }}
\newcommand\cS{{\mathcal{S}  }}

\newcommand\cZ{{\mathcal{Z}  }}

\begin{document}

\title[Khinchin Theorem for i.e.t.s]{Khinchin Theorem for interval exchange transformations.}
\author{Luca Marchese}
\email{l.marchese@sns.it}
\address{Section de math\'ematiques, case postale 64, 2-4 Rue du Li\`evre, 1211 Geneve, Suisse.}


\date{April 28, 2010}

\begin{abstract}
We define a diophantine condition for interval exchange transformations (i.e.t.s). When the number of intervals is two, that is for rotations on the circle, our condition coincides with classical Khinchin condition. We prove for i.e.t.s the same dichotomy as in Khinchin Theorem. We also develop several results relating the Rauzy-Veech algorithm with homogeneous approximations for i.e.t.s.
\end{abstract}



\maketitle

\tableofcontents

\section{Introduction}

A rotation on the circle is completely described by its \emph{rotation number} and we can establish a very precise dictionary between the dynamical properties of the transformations in this class and the diophantine properties of the rotation number. For example a rotation is periodic if and only if its rotation number is rational and conversely all irrational rotations are minimal. In general a diophantine condition on $\te\in [0,1)$ concerns the set of solutions $n\in\NN$ of
$$
\{n\te\}<\varphi(n),
$$
where $\{\cdot\}$ denotes the fractionary part and $\{\varphi(n)\}_{n\in\NN}$ is a positive sequence. Khinchin proved the following classical result (see \cite{kin}).

\begin{thmnonnum}[Khinchin]
Let $\{\varphi(n)\}_{n\in\NN}$ be a positive sequence.
\begin{itemize}
\item
If $\sum_{n=1}^{+\infty}\varphi(n)<+\infty $ then for almost any $\te$ there exists just finitely many $n\in\NN$ such that $\{n\te\}<\varphi(n)$.
\item
If $n\varphi(n)$ is decreasing monotone and $\sum_{n=1}^{+\infty} \varphi(n)=+\infty $ then for almost any $\te$ there exist infinitely many $n\in\NN$ such that $\{n\te\}<\varphi(n)$.
\end{itemize}
\end{thmnonnum}

\subsection*{Interval exchange transformations}

An alphabet is a finite set $\cA$ with $d\geq2$ elements. An \emph{interval exchange transformation} (also called i.e.t.) is a map $T$ from an interval $I$ to itself such that $I$ admits two partitions $\cP_{t}:=\{I^{t}_{\xi}\}_{\xi\in\cA}$ and $\cP_{b}:=\{I^{b}_{\xi}\}_{\xi\in \cA}$ into $d$ open intervals and for any $\xi\in\cA$ the restriction of $T$ to the interval $I^{t}_{\xi}$ is a translation with image the interval $I^{b}_{\xi}$. The map $T:I\to I$ is completely defined by the following data:
\begin{enumerate}
\item
The lengths of the intervals, called \emph{length data}. They are given by a vector $\la$ in $\RR_{+}^{\cA}$, where for any $\xi\in\cA$ the coordinate $\la_{\xi}$ equals the length of $I^{t}_{\xi}$, which is also equal to the length of $I^{b}_{\xi}$.
\item
The order before and after rearranging, called \emph{combinatorial data}. They are given by a pair of bijections $\pi=(\pi^{t},\pi^{b})$ from $\cA$ to $\{1,\dots,d\}$. For any $\xi\in\cA$, if we count starting from the left, the interval $I^{t}_{\xi}$ is in $\pi^{t}(\xi)$-th position in $\cP_{t}$ and the interval $I^{b}_{\xi}$ is in $\pi^{b}(\xi)$-th position in $\cP_{b}$.
\end{enumerate}

For a combinatorial datum $\pi$ set $\De_{\pi}:=\{\pi\}\times \RR^{\cA}_{+}$, which is the set of all i.e.t.s with combinatorial datum $\pi$. Consider any $T$ in $\De_{\pi}$ and write $T=(\pi,\la)$, where $\la$ is the corresponding length datum. For any $\xi\in\cA$ with $\pi^{t}(\xi)>1$ we call $u^{t}_{\al}$ the left endpoint of $I^{t}_{\al}$. In general $T$ is not continuous at $u^{t}_{\al}$. Similarly for any $\be\in\cA$ with $\pi^{b}(\be)>1$ we call $u^{b}_{\be}$ the left endpoint of $I^{b}_{\be}$. In general the inverse $T^{-1}$ of $T$ is not continuous at $u^{b}_{\be}$. If we identify the interval $I$ with $(0,\sum_{\al\in\cA}\la_{\al})$ then the position of the singularities $u^{t}_{\al}$ and $u^{b}_{\be}$ is given by
$$
u^{t}_{\al}:=
\sum_{\pi^{t}(\xi)<\pi^{t}(\al)}\la_{\xi}
\textrm{ and }
u^{b}_{\be}:=
\sum_{\pi^{b}(\xi)<\pi^{b}(\be)}\la_{\xi}.
$$
We say that the combinatorial datum $\pi$ is \emph{admissible} if there is no proper subset $\cA'\subset \cA$ with $k<d$ elements such that $\pi^{t}(\cA')= \pi^{b}(\cA')=\{1,\dots,k\}$. A \emph{connection} for the i.e.t. $T$ is a triple $(\be,\al,n)$ with $\pi^{b}(\be)>1$, $\pi^{t}(\al)>1$ and $n\in \NN$ such that $T^{n}u^{b}_{\be}=u^{t}_{\al}$. In particular, if $T$ has no connections then $\pi$ is admissible.

Rauzy \cite{rauzy}, Veech \cite{veech} and then Zorich \cite{zorich} introduced a continued fraction algorithm for i.e.t.s, called \emph{Rauzy-Veech(-Zorich) algorithm}. The algorithm is performed by a map $Q$ acting on the parameter space of all i.e.t.s, which generalizes the \emph{Gauss map}. Starting from $T$ it produces a sequence of i.e.t.s $T^{(1)}=Q(T),\dots,T^{(r)}=Q^{r}(T),\dots$ with the same number of intervals. I.e.t.s admitting infinitely many steps of the algorithm are those $T$ without connections, exactly as irrational real numbers are the points where the Gauss map can be iterated infinitely many times. For this reason we say that i.e.t.s with connections play the role of rational elements and we state a diophantine condition on i.e.t.s in terms of \emph{homogeneous approximation}, that is approximation of connections. Such a point of view is also motivated be Keane's Theorem (see \cite{keane}), which says that if $T$ has no connections, then it is minimal (we remark anyway that there is not a dichotomy as for rotations, since there exist minimal i.e.t.s with connections).

\subsubsection*{Khincin type condition for i.e.t.s}

Consider a positive sequence $\{\varphi(n)\}_{n\in\NN}$ and an admissible combinatorial datum $\pi$. Let $T$ be an i.e.t. in $\De_{\pi}$, acting on an interval $I$. We consider triples $(\be,\al,n)$ with $\pi^{b}(\be)>1$, $\pi^{t}(\al)>1$ and $n\in\NN$ such that
\begin{equation}\label{eqKhinTSI}
|T^{n}(u_{\be}^{b})-u_{\al}^{t}|<\varphi(n).
\end{equation}
For a triple $(\be,\al,n)$ which is not a connection for $T$ we denote $I(\be,\al,n)$ the open subinterval of $I$ whose endpoints are $T^{n}(u^{b}_{\be})$ and $u^{t}_{\al}$.

\begin{defi}\label{defreducedtriple}
Let $\pi$ be an admissible combinatorial datum and $(\be,\al,n)$ be a triple with $n\in\NN$, $\pi^{b}(\be)>1$ and $\pi^{t}(\al)>1$. We say that $(\be,\al,n)$ is a \emph{reduced triple} for $T$ in $\De_{\pi}$ if it is not a connection for $T$ and moreover for any $k\in \{0,\dots,n\}$ the pre-image $T^{-k}(I(\be,\al,n))$ does not contain any singularity $u^{t}_{\xi}$ of $T$ or any singularity $u^{b}_{\xi}$ of $T^{-1}$.
\end{defi}

Note that if $T$ has no connections, then the singularities $u^{t}_{\xi}$ and $u^{b}_{\xi}$ cannot be in the boundary of $T^{-k}(I(\be,\al,n))$.

\begin{defi}\label{defdiofanteoliouvilletsi}
Let $\{\varphi(n)\}_{n\in\NN}$ be a positive sequence and $\pi$ be admissible.
\begin{itemize}
\item
An i.e.t. $T$ in $\De_{\pi}$ is said $\modphi$\emph{-Diophantine} if condition (\ref{eqKhinTSI}) is satisfied just by finitely many triples.
\item
Conversely $T$ is said $\modphi$\emph{-Liouville} if for any pair of letters $(\be,\al)$ with $\pi^{b}(\be)>1$ and $\pi^{t}(\al)>1$ there exists infinitely many $n\in\NN$ such that the triple $(\be,\al,n)$ is reduced for $T$ and satisfies (\ref{eqKhinTSI}).
\end{itemize}
\end{defi}

The main result in this paper is the following generalization of Khinchin Theorem.

\begin{thm}\label{teoremaa}
Let $\pi_{0}$ be an admissible combinatorial datum and $\{\varphi(n)\}_{n\in\NN}$ be a positive sequence.
\begin{itemize}
\item
If $\varphi(n)$ is decreasing monotone with $\sum_{n=1}^{+\infty}\varphi(n)<+\infty$ then almost any $T$ in $\De_{\pi_{0}}$ is $\modphi$-Diophantine.
\item
If $n\varphi(n)$ is decreasing monotone and $\sum_{n=1}^{+\infty}\varphi(n)=+\infty$ then almost any $T$ in $\De_{\pi_{0}}$ is $\modphi$-Liouville.
\end{itemize}
\end{thm}

In general, even without any assumption on monotonicity, $\modphi$-Liouville i.e.t.s exist by the following proposition.

\begin{prop}\label{proposizioneg}
For any positive sequence $\{\varphi(n)\}_{n\in\NN}$ there exists a residual set in $\De_{\pi_{0}}$ of $\varphi$-liouville i.e.t.s.
\end{prop}

\subsubsection*{Borel-Cantelli setting and normalization of lengths}

The Rauzy-Veech algorithm has interesting recurrence properties just on the space of rays in $\De_{\pi}$ rather than on $\De_{\pi}$ itself, therefore it is useful to consider i.e.t.s $T$ acting on an interval $I$ with length one. For $\la\in\RR^{\cA}_{+}$ set $\|\la\|:=\sum_{\xi\in\cA}\la_{\xi}$ and denote $\De^{(1)}_{\pi}$ the $(d-1)$-simplex of those $T=(\pi,\la)$ in $\De_{\pi}$ with $\|\la\|=1$. The left-hand side of equation (\ref{eqKhinTSI}) is homogeneous in $\la$, indeed if for some positive $\rho$ we change $T=(\pi,\la)$ with $T'=(\pi,\la')$ where $\la'=\rho\la$, then the quantity $|T^{n}u^{b}_{\be}-u^{t}_{\al}|$ changes by the same factor. It is also obvious that $(\be,\al,n)$ is a reduced triple for $T$ if and only if it is reduced for $T'$. On the other hand, if $\{\varphi(n)\}_{n\in\NN}$ is a positive sequence satisfying the assumption in Theorem \ref{teoremaa} (either for the convergent or for the divergent case), dividing $\varphi(n)$ by $\rho$ the same assumption are still satisfied. For this reason Theorem \ref{teoremaa} admits an equivalent statement on the simplex $\De^{(1)}_{\pi_{0}}$ with respect to its Lebesgue measure and we can decompose the proof of the theorem into the proof of the following two propositions (see discussion after Proposition \ref{propcasoconvergente} and Proposition \ref{propcasodivergente}).

\begin{prop}\label{propcasoconvergente}
Let $\pi_{0}$ and $(\be,\al)$ be respectively a combinatorial datum and a pair of letters as in Theorem \ref{teoremaa}. Let $\varphi(n)$ be a positive and decreasing monotone sequence such that $\sum_{n=1}^{\infty}\varphi(n)<+\infty$. Then for almost any $T\in\De^{(1)}_{\pi_{0}}$ there exist just finitely many triples $(\be,\al,n)$ which satisfy equation (\ref{eqKhinTSI}).
\end{prop}

\begin{prop}\label{propcasodivergente}
Let $\pi_{0}$ and $(\be,\al)$ be respectively a combinatorial datum and a pair of letters as in Theorem \ref{teoremaa}. Let $\varphi(n)$ be a positive sequence such that $n\varphi(n)$ is decreasing monotone and $\sum_{n=1}^{+\infty}\varphi(n)=+\infty$, then for almost any $T\in\De^{(1)}_{\pi_{0}}$ there exist infinitely many triples $(\be,\al,n)$ which are reduced for $T$ and satisfy (\ref{eqKhinTSI}).
\end{prop}

For any triple $(\be,\al,n)$ as in Theorem \ref{teoremaa}, the set of those $T$ in $\De^{(1)}_{\pi_{0}}$ such that $(\be,\al,n)$ is reduced for $T$ and satisfies equation (\ref{eqKhinTSI}) defines an event in $\De^{(1)}_{\pi_{0}}$. Roughly speaking our strategy is to prove that such event has probability proportional to $\varphi(n)$ and that, when the triple $(\be,\al,n)$ varies, we have some weak form of independence for the family of the associated events. Proposition \ref{propcasoconvergente} and Proposition \ref{propcasodivergente} then follows as consequences according to the well-known \emph{Borel-Cantelli lemma} (there are many good references, see for example \cite{billi}).

\begin{thmnonnum}
Let $(X,\PP)$ be a probability space and let $(X_{n})_{n\in\NN}$ be a countable family of events in $X$.
\begin{itemize}
\item
If $\sum_{n=1}^{\infty}\PP(X_{n})<+\infty$ then almost any $x\in X$ belongs to finitely many events $X_{n}$.
\item
On the other hand, if $\sum_{n=1}^{\infty}\PP(X_{n})=+\infty$ and the events $X_{n}$ are each other independent, then almost any $x\in X$ belongs to infinitely many $X_{n}$.
\end{itemize}
\end{thmnonnum}

\subsubsection*{Reduced triples}

An i.e.t. $T:I\to I$ with $d=2$ intervals is a rotation. If $I=[0,1)$ and $\te\in[0,1)$ is the rotation number, then the singularities of $T$ and $T^{-1}$ are respectively $1-\te$ and $\te$. Suppose that $\te$ is irrational, which is equivalent to say that $T$ has not connections. For any $n\in\NN$ we have a pair $q,p$ of non-negative integers such that $|T^{n}(\te)-(1-\te)|=|q\te-p|$. In particular a reduced triple corresponds to a pair $q,p$ such that
$$
|q\te-p|<|q'\te-p'|
$$
for any pair of integers $q',p'$ with $q'\leq q$, $p'\leq p$ and such that $p'/q'\not=p/q$ and $\sign(q'\te-p')=\sign(q\te-p)$ (see \cite{gugu}). It is well known that in this case $p/q$ is equal to a convergent $p_{k}/q_{k}$ of the continued fraction expansion of $\te$. We say therefore the that \emph{the continued fraction algorithm detects reduced triples}. On the other hand for any convergent $p_{k}/q_{k}$ of $\te$ the quantity $|q_{k}\te-p_{k}|$ corresponds to a reduced triple for $T$, thus we say that \emph{the continued fraction algorithm produces infinitely many reduced triples}. For i.e.t.s with $d>2$ intervals the notion or reduced triple in Definition \ref{defreducedtriple} satisfies a generalization of the above property.

\begin{thm}\label{teoremaf}
Let $\pi$ be an admissible combinatorial datum and $(\be,\al)$ be a pair with $\pi^{b}(\be)>1$ and $\pi^{t}(\al)>1$.
\begin{enumerate}
\item
Let $T$ in $\De_{\pi}$ without connections. Any triple $(\be,\al,n)$ reduced for $T$ is \emph{detected} at some step $T^{(r)}$ of the Rauzy-Veech algorithm.
\item
Conversely for almost any $T$ in $\De_{\pi}$ there are infinitely many steps $\{T^{(r_{k})}\}_{k\in\NN}$ of the Rauzy-Veech algorithm producing a reduced triple $(\be,\al,n(k))$ for $T$.
\end{enumerate}
\end{thm}

Detection and production of reduced triples is defined rigorously in Definition \ref{defdetectingproducing}. In Proposition we give an estimate on the total measure of reduced triples. Theorem \ref{teoremaf} says that reduced triples are the good notion to study homogeneous approximation via the Rauzy-Veech algorithm. They have an important role also in relation to \emph{translation surfaces} (see next subsection). In particular in \cite{luca3} it is shown that for any triple reduced for $T$ we can obtain a saddle connection for a proper translation surface $X$.

\subsection*{Some related results}

\subsubsection*{Translation surfaces and Teichm\"uller flow}

A \emph{translation surface} (also said \emph{abelian differential}) is a pair $(X,w)$, where $X$ is a compact Riemann surface and $w$ an holomorphic one-form on $X$. Equivalently $X$ has a flat metric with isolated conical singularities (the zeros of $w$) whose angle is an integer multiple of $2\pi$. I.e.t.s are strictly linked to translation surfaces and to the \emph{Teichm\"uller flow} on their \emph{moduli space} (see \cite{veech}, \cite{masuruno} and \cite{zorichflatsurf}). More precisely, any translation surface $X$ has an unitary constant vector field whose \emph{first return map} to a transverse segment $I$ in $X$ is an i.e.t..

A \emph{saddle connection} is a geodesic $\ga$ for the flat metric of $X$ starting and ending in two conical singularities and not containing any other conical singularity in its interior. The set $\hol(X)$ of \emph{periods} of $X$ is the set of complex numbers $v:=\int_{\ga}w$, where $\ga$ is a saddle connection for $X$ and $w$ is the holomorphic one-form. Around any conical singularity $p_{i}$ with angle $2k_{i}\pi$ we can track $k_{i}$ angular sectors with amplitude $2\pi$. If $p_{1},\dots,p_{r}$ are the conical singularities of $X$ then there are in total $d-1=k_{1}+\dots+k_{r}$ such sectors, where $d$ is the number of intervals of some i.e.t.. Therefore $\hol(X)$ splits in subsets $\hol_{(j,i)}(X)$ with $j,i\in\{1,\dots,d-1\}$, where $v\in\hol_{(j,i)}(X)$ if and only if the corresponding saddle connection starts in the sector $j$ and ends in the sector $i$. For a period $v$ consider the condition
\begin{equation}\label{e3}
|\Re(v)|\leq\varphi(|v|),
\end{equation}
where $|v|$ is the modulus of $v$ and $\varphi:(0,+\infty)\to (0,+\infty)$ is a positive function bounded from above. In \cite{luca3} we develop the counterpart of Theorem \ref{teoremaa} for translation surfaces. In particular the following dichotomy holds.

\begin{thmnonnum}
If $\varphi(t)$ is decreasing monotone with $\int_{0}^{+\infty}\varphi(t)dt<+\infty$ then $\hol(X)$ contains just finitely many solutions $v$ of (\ref{e3}) for almost any translation surface $X$.

If $t\varphi(t)$ is decreasing monotone with $\int_{0}^{+\infty}\varphi(t)dt=+\infty$ then $\hol_{(j,i)}(X)$ contains infinitely many solutions $v$ of (\ref{e3}) for almost any translation surface $X$ and any $j,i\in\{1,\dots,d-1\}$.
\end{thmnonnum}

The expression \emph{almost any} $X$ in the theorem above is meant with respect to Lebesgue measure on the \emph{stratum} $\cH(k_{1},\dots,k_{r})$ of the \emph{moduli space} of translation surfaces with $r$ conical singularities whose angles are $2k_{1}\pi,\dots,2k_{r}\pi$, which is a non-compact complex \emph{orbifold} with $\dim_{\CC}=k_{1}+\dots +k_{r}+1$. Compact subset in the stratum are characterized by the \emph{systole} function $X\mapsto\sys(X)$, defined as the length of the shortest saddle connection of the translation surface $X$: a sequence $X_{n}$ diverges in $\cH(k_{1},\dots,k_{r})$, that is it leaves any compact set, if and only if $\sys(X_{n})\to 0$. In \cite{luca3} we relate the diophantine property in the last theorem to the asymptotic behavior of orbits of the \emph{Teichm\"uller flow} $\cF_{t}$ on $\cH(k_{1},\dots,k_{r})$ and we prove the following sharp estimate.

\begin{thmnonnum}
Let $\psi:[0,+\infty)\to(0,+\infty)$ be a monotone decreasing function. If $\int_{0}^{\infty}\psi(t)dt<+\infty$ then for almost any $X$ in $\cH(k_{1},\dots,k_{r})$ we have
$$
\lim_{t\to\infty}
\frac{\sys(\cF_{t}X)}
{\sqrt{\psi(t)}}
=\infty.
$$
On the other hand, if $\int_{0}^{\infty}\psi(t)dt=+\infty$, then for almost any $X$ in $\cH(k_{1},\dots,k_{r})$ we have
$$
\liminf_{t\to\infty}
\frac{\sys(\cF_{t}X)}
{\sqrt{\psi(t)}}
=0.
$$
\end{thmnonnum}
\emph{Note.} In fact in \cite{luca3} it is proved a stronger result. Indeed for any $j,i$ in $\{1,\dots,d-1\}$ we can replace $\sys(X)$ with $\sys_{(j,i)}(X)$, defined as the length of the shortest saddle connection starting in the sector $j$ and ending in the sector $i$.\\

Applying both parts of the last theorem to the family of functions $\psi_{r}(t):=\min\{1,t^{-r}\}$ with $r\geq 1$, we get for almost any $X$:
$$
\limsup_{t\to\infty}
\frac{-\log\circ\sys(\cF_{t}X)}
{\log t}
=1/2.
$$
\emph{Masur's logarithm law} says that the maximal excursion up to time $t$ of $\distanza([X],[\cF_{t}X])$ has amplitude $\log\sqrt{t}$, where $[X]$ and $[\cF_{t}X]$ are the Riemann surfaces underlying respectively to $X$ and $\cF_{t}X$ and where $\distanza$ denotes the \emph{Teichm\"uller distance} in the moduli space of Riemann surfaces of genus $g$ (see \cite{masurtre}). The last estimate shows that $-\log \circ\sys(\cF_{t}X)$ has the same asymptotic behavior as $\distanza([X],[\cF_{t}X])$ and is therefore a natural extension of Masur's result to strata of translation surfaces.

\begin{domanda}
\emph{Linear involutions} are a natural generalization of i.e.t.s introduced in \cite{danno} by Danthony and Nogueira. In \cite{bola} Boissy and Lanneau related linear involutions to \emph{quadratic differentials}, whose moduli space is the co-tangent bundle of the moduli space of complex curves, that is the natural setting for an extension of Masur's logarithmic law. We believe that the techniques introduced in this paper can be extended to linear involutions and we ask if a generalization of Theorem \ref{teoremaa} can be proved for them (more precisely for the subclass of linear involutions which are relevant for quadratic differentials, as it is explained in \cite{bola}). The question is also motivated by the paper of Avila and Resende (\cite{mariajoao}), where the authors generalize some results of \cite{agy} which play an important role in the proof of Theorem \ref{teoremaa}.
\end{domanda}

\subsubsection*{Non-homogeneous results}

Boshernitzan and Chaika studied shrinking target properties related to the diophantine condition in Definition \ref{defdiofanteoliouvilletsi} (see \cite{bocha} and \cite{chaika}). In particular \cite{chaika} it is proven the following non-homogeneous result.

\begin{thmnonnum}[Chaika]
Let $\{\varphi(n)\}_{n\in\NN}$ be a positive sequence such that $n\varphi(n)$ is decreasing monotone and $\sum_{n=1}^{\infty}\varphi(n)=\infty$. Then for almost any i.e.t. $T:I\to I$ with admissible combinatorial datum, for any $x$ in $I$ and for almost any $y\in I$, there are infinitely many $n\in\NN$ such that
$$
|T^{n}(x)-y|<\varphi(n).
$$
\end{thmnonnum}

\subsection*{Contents of this article}

In \emph{Section 2} we recall the basic theory of i.e.t.s. In \S \ref{ss1s2(algorithmR-V)} we introduce the map of Rauzy-Veech and Zorich's acceleration, in particular we define \emph{Rauzy classes}. In \S \ref{Reduction of rauzy classes} we describe a combinatorial operation introduced in \cite{agy} and called \emph{reduction of Rauzy classes}. The normalized Rauzy-Veech map is a piecewise liner-projective map, in paragraph \ref{thedistortionestimate} we describe the connected components of its domains and we give a formula for their volume. The volume has unbounded distortion under iteration of the map and in paragraph \ref{sss1ss3s2(distortionestimate)} we state a result proved in \cite{agy} on the control of the distortion.

\emph{Section 3} contains the proof of Theorem \ref{teoremaf} and Proposition \ref{proposizioneg}. Detection and production of reduced triples are defined in Definition \ref{defdetectingproducing}. Lemma \ref{lem1s3ss1sss0} proves that reduced triples are detected by the Rauzy-Veech algorithm and thus the first part of Theorem \ref{teoremaf}. In \S \ref{s3ss1sss1} we consider the shortest sequence of steps of the algorithm to be applied to $T$ in order to detect a triple $(\be,\al,n)$ reduced for $T$. Such sequence is denoted $\ga=\ga(\be,\al,n,T)$ and identifies a region $\De_{\ga}$ of the parameter space of i.e.t.s. Proposition \ref{prop1s3ss1sss1} gives a geometric description of those $T$ in the region $\De_{\ga}$ having $(\be,\al,n)$ as reduced triple. Then we consider the subspace $\De^{(1)}_{\ga}$ of $\De_{\ga}$ of those i.e.t.s with $\|\la\|=1$. Lemma \ref{lem1s3ss2sss2} gives a local estimate for the measure of those $T$ in $\De^{(1)}_{\ga}$ such that $(\be,\al,n)$ is reduced for $T$ and satisfies $|T^{n}u^{b}_{\be}-u^{t}_{\al}|<\varphi(n)$. Proposition \ref{props3ss2} gives a global estimate (rather implicit) on the measure of all $T$ with $\|\la\|=1$ such that $(\be,\al,n)$ is reduced for $T$. \S \ref{s3ss3} contains the proof of the second part of Theorem \ref{teoremaf}, which is a consequence of Proposition \ref{props3ss3}. For any pair $(\be,\al)$ and for a generic $T$ the proposition gives  a family of sequences $\{\ga_{k}\}_{k\in\NN}$ of steps of the algorithm such that any $\ga_{k}$ produces a triple $(\be,\al,n_{k})$ reduced for $T$. Proposition \ref{props3ss3} requires a combinatorial property of Rauzy classes which is consequence of Theorem \ref{propcombinatoria}. Proposition \ref{proposizioneg} is proved in \S \ref{s3ss3sss4}.

\emph{Section 4} is devoted to the proof of Theorem \ref{teoremaa}. We consider normalized length data and with Lemma \ref{lem1s4ss0sss0} we prove that Theorem \ref{teoremaa} is consequence of Proposition \ref{propcasoconvergente} and Proposition \ref{propcasodivergente}. In \S \ref{s4ss1} we prove Proposition \ref{propcasoconvergente}, which follows directly from the results in \S \ref{s3ss2}. In \S \ref{s4ss2} we prove Proposition \ref{props4ss2}, which establish a sufficient condition implying Proposition \ref{propcasodivergente}. We consider the half-infinite sequence $\{\ga_{k}(T)\}_{k\in\NN}$ of steps of the Rauzy-Veech algorithm generated by a generic $T$. Theorem \ref{teoremaf} implies that for any pair $(\be,\al)$ there exists a \emph{reference sequence} $\eta$ (and thus a region $\De^{(1)}_{\eta}$ in the parameter space of normalized i.e.t.s) with the following property: the algorithm produces a triple $(\be,\al,n_{k})$ reduced for $T$ each time that $\eta$ appears in $\ga_{k}(T)$. This amount to consider the \emph{first return map} $\cF_{\eta}$ to the region $\De^{(1)}_{\eta}$ (see Definition \ref{defreferencepath} and equation (\ref{eqdefmappaeffe})).  Proposition \ref{props4ss2} proves that having infinitely many triples $(\be,\al,n_{k})$ reduced for $T$ and satisfying (\ref{eqKhinTSI}) corresponds to a \emph{shrinking target property} for the map $\cF_{\eta}$. In \S \ref{s4ss3} we treat a technical issue, that is we provide a family of shrinking targets $\{\cE_{k}\}_{k\in\NN}$ which are measurable with respect to the sigma-algebra generated by the connected components of the domain of $\cF_{\eta}$, this is necessary to convert Proposition \ref{props4ss2} in the setting of the \emph{Borel-Cantelli Lemma}. Finally in \S \ref{s4ss4} we prove that the shrinking-target criterion is satisfied (Proposition \ref{props4ss4}).

In \emph{Section 5} we state and prove two general results for i.e.t.s. The first is Theorem \ref{propcombinatoria}, which affirms that for any pair of letters $(\be,\al)$ as in Theorem \ref{teoremaa}, any Rauzy class contains an element $\pi$ where $\al$ and $\be$ are in some required reciprocal position. The second general result is Theorem \ref{teoremastimaprincipale}, which implies that the measure of the targets $\{\cE_{k}\}_{k\in\NN}$ constructed in \ref{s4ss3} does not decay too fast.

\subsection*{Acknowledgements}

The author would like to thank Jean-Christohpe Yoccoz for many discussions, and Stefano Marmi, Giovanni Forni, Pascal Hubert and anonymous referees for many questions and precious remarks.

\section{Background theory}\label{Background Theory}

\subsection{Rauzy-Veech Algorithm}\label{ss1s2(algorithmR-V)}

This subsection is a brief survey of the basic properties of the Rauzy-Veech algorithm. We follow \cite{mmy} and \cite{agy}.\\

Let $\pi=(\pi^{t},\pi^{b})$ and $\lambda$ define an interval exchange transformation $T:I\to I$. Let $\epsilon\in \lbrace t,b \rbrace$, where the letter $t$ stands for \emph{top} and the letter $b$ for \emph{bottom}. If $\epsilon=t$ we put $1-\epsilon:=b$ and if $\epsilon=b$ we put $1-\epsilon:=t$. Let us call $\al_{t}$ and $\al_{b}$ the two letters in $\cA$ such that respectively $\pi^{t}(\al_{t})=d$ and $\pi^{b}(\al_{b})=d$. The rightmost singularity of $T$ is therefore $u^{t}_{\al_{t}}$ and the rightmost singularity of $T^{-1}$ is $u^{b}_{\al_{b}}$. We suppose that
\begin{equation}\label{eq1ss1s2(algorithmR-V}
u^{t}_{\al_{t}}\not= u^{b}_{\al_{b}}
\end{equation}
and we consider the value of $\epsilon\in \lbrace t,b \rbrace$ such that
$$
u^{\epsilon}_{\al_{\epsilon}}
<
u^{1-\epsilon}_{\al_{1-\epsilon}}.
$$
With this Definition of $\epsilon$  we say that $T$ is of \emph{type} $\epsilon$. We also say that the letter $\al_{\epsilon}$ is the \emph{winner} of $T$ and $\al_{1-\epsilon}$ is the \emph{loser}. We consider the subinterval of $I$
$$
\widetilde{I}:=
I\cap(0,u^{1-\epsilon}_{\al_{1-\epsilon}})
$$
and we define $\widetilde{T}:\widetilde{I}\to\widetilde{I}$ as the first return map of $T$ to $\widetilde{I}$. It is easy to check that $\widetilde{T}$ is an i.e.t.. The combinatorial datum $\widetilde{\pi}=(\widetilde{\pi}^{t},\widetilde{\pi}^{b})$ of $\widetilde{T}$ is given by:
\begin{equation}\label{eq2ss1s2(algorithmR-V)}
\begin{array}{lllll}
\widetilde{\pi}^{\epsilon}(\al)=
\pi^{\epsilon}(\al)\,
\forall\al\in\cA\\
\\
\widetilde{\pi}^{1-\epsilon}(\al)=
\pi^{1-\epsilon}(\al)\,
\textrm{ if }
\pi^{1-\epsilon}(\al)\leq
\pi^{1-\epsilon}(\al_{\epsilon})\\
\widetilde{\pi}^{1-\epsilon}(\al_{1-\epsilon})=
\pi^{1-\epsilon}(\al_{\epsilon})+1\\
\widetilde{\pi}^{1-\epsilon}(\al)=
\pi^{1-\epsilon}(\al)+1\,
\textrm{ if }
\pi^{1-\epsilon}(\al_{\epsilon})<\pi^{1-\epsilon}(\al)<d.
\end{array}
\end{equation}
The length datum $\widetilde{\lambda}$ of $\widetilde{T}$ is given by:
\begin{equation}\label{eq3ss1s2(algorithmR-V)}
\begin{array}{ll}
\widetilde{\lambda}_{\al}=
\lambda_{\al}\,
\textrm{ if }
\al\neq\al_{\epsilon}\\
\widetilde{\lambda}_{\al_{\epsilon}}=
\lambda_{\al_{\epsilon}}-\lambda_{\al_{1-\epsilon}}.
\end{array}
\end{equation}
When $T=(\pi,\la)$ satisfies condition (\ref{eq1ss1s2(algorithmR-V}), equations (\ref{eq2ss1s2(algorithmR-V)}) and (\ref{eq3ss1s2(algorithmR-V)}) define a map $T\mapsto Q(T):=\widetilde{T}$ which is known as \emph{Rauzy-Veech map}. We introduce two operations $R^{t}$ and $R^{b}$ on the set of admissible combinatorial data as follows. If $\epsilon$ is the type of $T$ and $\pi$ is its combinatorial datum, then we set $R^{\epsilon}(\pi):=\widetilde{\pi}$, where $\widetilde{\pi}$ is the combinatorial datum of $\widetilde{T}$. It is easy to check that if $\pi$ is an admissible combinatorial datum then both $R^{t}(\pi)$ and $R^{b}(\pi)$ are admissible.

\begin{defi}\label{def1ss1s2(algorithmR-V)}
Let us call $\mathfrak{S}$ the set of all the admissible combinatorial data $\pi$ over some alphabet $\cA$. The maps $R^{t}$ and $R^{b}$ from $\mathfrak{S}$ to itself are called the \emph{Rauzy elementary operations}.
\begin{itemize}
\item
A \emph{Rauzy class} is a minimal non-empty subset $\cR$ of $\mathfrak{S}$ which is invariant under $R^{t}$ and $R^{b}$.
\item
A \emph{Rauzy diagram} is a connected oriented graph $\cD$ whose vertexes are the elements of $\cR$ and whose oriented arcs, or \emph{arrows}, correspond to Rauzy elementary operations $\pi\mapsto R^{\varepsilon}(\pi)$ between elements of $\cR$. \item
An arrow corresponding to $R^{t}$ is called a \emph{top} arrow and we say that $\al_{t}$ is its winner and $\al_{b}$ is its loser. Conversely an arrow corresponding to $R^{b}$ is called a \emph{bottom} arrow and we say that $\al_{t}$ is its loser and $\al_{b}$ is its winner.
\item
A concatenation of compatible arrows in a Rauzy diagram is called a \emph{Rauzy path}. The set of all Rauzy paths connecting elements of $\cR$ is denoted $\Pi (\cR)$. If a path $\ga$ is concatenation of $r$ simple arrows, we say that $\ga$ has length $r$. Length one paths are arrows, we also identify elements of $\cR$ with trivial paths, that is length zero paths.
\item
A partial ordering $\prec$ is defined on $\Pi(\cR)$ saying that $\nu\prec\ga$ iff $\ga$ begins with $\nu$. A subfamily $\Ga$ of $\Pi(\cR)$ is called \emph{disjoint} iff for any two elements $\eta$ and $\nu$ of $\Ga$ we have neither $\eta\prec\ga$ nor $\ga\prec\eta$.
\end{itemize}
\end{defi}

With the notation above, recalling that for any combinatorial datum $\pi$ we defined $\De_{\pi}=\{\pi\} \times \RR_{+}^{\cA}$, we denote $\Delta(\cR):=\bigsqcup_{\pi\in \cR}\Delta_{\pi}$ the set of all the intervals exchange transformations with combinatorial datum in the Rauzy class $\cR$.

\subsubsection{Linear action}\label{sss1ss1s2}

For any Rauzy class $\cR$ and any path $\ga\in \Pi(\cR)$ we define a linear map $B_{\ga}\in \sldz$ as follows. If $\ga$ is trivial then $B_{\ga}=id$. If $\ga$ is an arrow with winner $\al$ and loser $\be$ then we set $B_{\ga}e_{\al}=e_{\al}+e_{\be}$ and $B_{\ga}e_{\xi}=e_{\xi}$ for all $\xi\in \cA\setminus\{\al\}$, where $\{e_{\xi}\}_{\xi\in \cA}$ is the canonical basis of $\RR^{\cA}$. We extend the definition to paths so that $B_{\ga_{1}\ga_{2}}=B_{\ga_{2}}B_{\ga_{1}}$.

Let us fix some element $\pi$ in the Rauzy class $\cR$. For any $\ga\in\Pi(\cR)$ starting at $\pi$ we define the simplicial sub-cone $\De_{\ga}\subset\De_{\pi}$ by
$$
\Delta_{\ga}=\{\pi\}\times^{t}B_{\ga}(\RR_{+}^{\cA}),
$$
where $^{t}B_{\ga}$ denotes the trasposed of the matrix $B_{\ga}$ defined above. For the same $\ga$ we also define the vector $q^{\ga}\in\NN^{\cA}$ by
$$
q^{\ga}:=B_{\ga}\vec{1},
$$
where $\vec{1}$ denotes the column vector of $\NN^{\cA}$ whose entries are all equal to $1$. If $\ga$ is a path with length $r$ we often write $q^{(r)}$ instead of $q^{\ga}$.

\subsubsection{Iteration of the algorithm}\label{s2ss1sss2}

When $T\in\De(\cR)$ is such that the $r$-th iterated of $Q$ is defined, we have an explicit formula for $T^{(r)}:=Q^{r}(T)$.

\begin{lem}\label{lem1s2ss1sss2}
Let $\ga\in\Pi(\cR)$ be a path in the Rauzy diagram with length $r$ and let $B_{\ga}$ and $\De_{\ga}$ be respectively the matrix and the simplicial cone defined in \S \ref{sss1ss1s2}. Then for any $T\in\De_{\ga}$ the $r$-th iterated of $Q$ is defined and the length datum $\la^{(r)}$ of $Q^{r}(T)$ is given by the formula
$$
\la^{(r)}=^{t}B_{\ga}^{-1}\la.
$$
\end{lem}

\begin{demo}
Let us first consider an arrow $\ga$. Call $\pi$ its starting point, $\al$ its winner and $\be$ its loser. Suppose that $\ga$ is of type top, the other case being identic. Consider $T=(\pi,\la)$ in $\De_{\ga}$. By definition of the matrix $B_{\ga}$ we have
$$
B_{\ga}=id +E_{\be,\al}
$$
where $E_{\be,\al}$ is the matrix whose entry in the column $\al$ and row $\be$ is $1$ and all the others are $0$. Hence $^{t}B_{\ga}\RR_{+}^{d}$ is the open half-cone of those $\la$ in $\RR_{+}^{d}$ with $\la_{\al}>\la_{\be}$. Since $\pi^{t}(\al)=\pi^{b}(\be)=d$ this last condition is equivalent to condition (\ref{eq1ss1s2(algorithmR-V}) and therefore $Q$ is defined on $T$. Moreover the length datum $\la^{(1)}$ of $Q(T)$ is given by $\la^{(1)}=^{t}B_{\ga}^{-1}\la$, according to equation (\ref{eq3ss1s2(algorithmR-V)}).

The proof goes on by induction on $r$. Suppose that the lemma is proved for any concatenation $\ga_{1}...\ga_{r-1}$ of $r-1$ arrows. Consider a path $\ga=\ga_{1}\dots\ga_{r-1}\ga_{r}$ starting at $\pi$, where $\ga_{r}$ is an arrow which can be concatenated to $\ga_{r-1}$. Consider any $T=(\pi,\la)$ in $\De_{\ga}=\{\pi\}\times^{t}B_{\ga}\RR_{+}^{d}$. Since $B_{\ga}=B_{\ga_{r}}B_{\ga_{1}...\ga_{r-1}}$ it follows that $\De_{\ga}$ is a sub-cone of $\De_{\ga_{1}...\ga_{r-1}}$, therefore the inductive hypothesis applies, that is $T^{(r-1)}=Q^{r-1}(T)$ is defined and its length datum is given by $\la^{(r-1)}= ^{t}B^{-1}_{\ga_{1}...\ga_{r-1}}\la$. Moreover we have $T^{(r-1)}\in\De_{\ga_{r}}$, thus applying the argument in the first part of the lemma we get that $T^{(r)}=Q(T^{(r-1)})$ is defined an its length datum is $\la^{(r)}=^{t}B^{-1}_{\ga_{r}}\la^{(r-1)}=^{t}B_{\ga}^{-1}\la$. The lemma is proved.
\end{demo}

For $r\in\NN$ denote $\De_{r}(\cR)$ the domain of $Q^{r}$. According to Lemma \ref{lem1s2ss1sss2} the connected components of $\De_{r}(\cR)$ are labeled by paths $\ga$ in $\Pi(\cR)$ of length $r$. If $\ga$ is such a path ending at $\pi'$ then $Q^{r}:\De_{\ga}\to\De_{\pi'}$ is a homeomorphism. The intersection $\De_{\infty}(\cR):=\bigcap_{r\in\NN}\De_{r}(\cR)$ is the set of those i.e.t. $T$ such that the map $Q$ can be applied infinitely many times. It is a set with full Lebesgue measure, since is the intersection of countably many sets of full Lebesgue measure. The complement of $\Delta_{\infty}(\cR)$ is the set of those i.e.t.s $T$ such that $T^{(r)}=Q^{r}(T)$ eventually does not satisfy condition $(\ref{eq1ss1s2(algorithmR-V})$, that is the algorithm stops. For the complement the following characterization holds (see Corollary 2 at page 37 of \cite{ytre}).

\begin{lemnonnum}
When applied to an i.e.t. $T$ the Rauzy algorithm $Q$ eventually stops if and only if $T$ has a connection.
\end{lemnonnum}

From now on, when applying the map $Q$, we will not worry about its domain, keeping in mind that it is defined on $\De^{(1)}(\cR)$ modulo a set of measure zero.

\subsubsection{Return times}\label{S2ss1sss3}

Let $\ga$ be Rauzy path with length $r$ and consider $T$ in $\De_{\ga}$. Let $I$ be the interval where $T$ acts. According to the definition of the Rauzy-Veech algorithm the i.e.t. $T^{(r)}=Q^{r}(T)$ is the first return of $T$ to the subinterval $I^{(r)}$ of $I$ whose left endpoint is the left endpoint of $I$ and whose right endpoint is the rightmost singularity of $T^{(r-1)}$ (top or bottom). For any $\al$ and $\be$ denote $u^{(r),t}_{\al}$ and $u^{(r),b}_{\be}$ the corresponding singularities respectively for $T^{(r)}$ and for $(T^{(r)})^{-1}$. Since $T^{(r)}$ is the first return map of $T$ to the subinterval $I^{(r)}$ then there exist two non-negative integers $l(\be,r)$ and $h(\al,r)$ such that
$$
u^{(r),b}_{\be}=T^{l(\be,r)}u^{b}_{\be}
\textrm{ and }
u^{(r),t}_{\al}=T^{-h(\al,r)}u^{t}_{\al}.
$$
Moreover $l(\be,r)$ and $h(\al,r)$ are the smallest non-negative integers such that respectively $T^{l(\be,r)}u^{b}_{\be}\in I^{(r)}$ and $T^{-h(\al,r)}u^{t}_{\al}\in I^{(r)}$. It is also evident that $l(\be,r)$ and $h(\al,r)$ just depend on $\ga$ (and of course on $\al$ or $\be$).

Let $\la^{(r)}$ be the length datum of $T^{(r)}$ and for any $\al$ and $\be$ define the intervals
$$
I^{(r),t}_{\al}:=
(u^{(r),t}_{\al},u^{(r),t}_{\al}+\la^{(r)}_{\al})
\textrm{ and }
I^{(r),b}_{\be}:=
(u^{(r),b}_{\be},u^{(r),b}_{\be}+\la^{(r)}_{\be}).
$$
They are the intervals where respectively $T^{(r)}$ and $(T^{(r)})^{-1}$ act as a translation. For any $\al$ we have $T^{(r)}(I^{(r),t}_{\al})=I^{(r),b}_{\al}\subset I^{(r)}$ and since $T^{(r)}$ is the first return map of $T$ to $I^{(r)}$, then there exists a positive integer $R=R(\al,r)$, called \emph{return time}, such that
$$
T^{(r)}(I^{(r),t}_{\al})=
T^{R}(I^{(r),t}_{\al})
\textrm{ and }
T^{k}(I^{(r),t}_{\al})\cap I^{(r)}=
\emptyset
$$
for all $k\in\{0,\dots,R(\al,r)-1\}$. Moreover $T^{(r)}$ acts as a translation on $I^{(r),t}_{\al}$, therefore for these values of $k$ the image $T^{k}(I^{(r),t}_{\al})$ does not contain in its interior any singularity $u^{t}_{\xi}$ for $T$ or any singularity $u^{b}_{\xi}$ for $(T)^{-1}$. In particular for $k\in\{0,\dots,R(\al,r)-1\}$ the intervals $T^{k}(I^{(r),t}_{\al})$ are disjoint each other and have length equal to $\la^{(r)}_{\al}$. Recall the vector $q^{(r)}=B_{\ga}\vec{1}$ defined in \S \ref{sss1ss1s2}. Lemma \label{lem1s2ss1sss2} implies that the length datum $\la^{(r)}$ satisfy $^{t}B_{\ga}\la^{(r)}=\la$, where $\la$ is the length datum of $T$. Therefore according to the discussion above, for any $\al$ and any $r$ we have
$$
R(\al,r)=q^{\ga}_{\al}.
$$

Finally we observe that for any $\al$ and $r$ the integers $h(\al,r)$ and $l(\al,r)$ satisfy
$
T^{h(\al,r)}(I^{(r),t}_{\al})=
(u^{t}_{\al},u^{t}_{\al}+\la^{(r)}_{\al})
\subset
I^{t}_{\al}
$
and
$
T^{-l(\al,r)}(I^{(r),b}_{\al})=
(u^{b}_{\al},u^{b}_{\al}+\la^{(r)}_{\al})
\subset
I^{b}_{\al}
$.
Moreover we have trivially $T(I^{t}_{\al})=I^{b}_{\al}$, thus
$
T^{h(r,\al)+1+l(r,\al)}(I^{(r),t}_{\al})=
I^{(r),b}_{\al}
$.
Therefore any $\al\in\cA$ and any $r$ we have
$$
h(r,\al)+1+l(r,\al)=q^{\ga}_{\al}.
$$

\subsubsection{Normalized Rauzy-Veech algorithm and Zorich's acceleration}\label{Normalized Rauzy Veech algorithm and Zorich's acceleration}

The Rauzy-Veech algorithm has interesting recurrence properties on the spaces of rays in $\De(\cR)$ rather that on $\De(\cR)$ itself. Therefore it is convenient to introduce a normalization on the sum of the lengths of the intervals. For $\la\in\RR^{\cA}_{+}$ recall the notation $\|\la\|:=\sum_{\xi\in\cA}\la_{\xi}$ and $\widehat{\la}:=\|\la\|^{-1}\la$. For any combinatorial datum $\pi$ in some $\cR$ we write
$$
\De^{(1)}_{\pi}:=
\{(\pi,\la)\in\De_{\pi};\|\la\|=1\}.
$$
The normalized length datum of an i.e.t. $T\in\De^{(1)}_{\pi}$ will be denoted $\widehat{\la}$. For any Rauzy class $\cR$ the set of all normalized i.e.t.s with combinatorial datum in $\cR$ is denoted $\De^{(1)}(\cR):=\bigsqcup_{\pi\in \cR}\De^{(1)}_{\pi}$.

\begin{defi}\label{defalgoritmorauzyveechnormalizzato}
Let $\cR$ be a Rauzy class over an alphabet $\cA$. The \emph{normalized Rauzy-Veech algorithm} is the map $\widehat{Q}:\De^{(1)}(\cR)\to \De^{(1)}(\cR)$ defined by
$$
\widehat{Q}(\pi,\la):=
(\widetilde{\pi},\frac{\widetilde{\la}}{\|\widetilde{\la}\|}),
$$
where $(\widetilde{\pi},\widetilde{\la})=Q(\pi,\la)$ is the Rauzy-Veech algorithm introduced in paragraph \ref{ss1s2(algorithmR-V)}.
\end{defi}

If $T\in\De^{(1)}(\cR)$ is an i.e.t. without connections, for any $r\in\NN$ we denote $\widehat{T}^{(r)}:=\widehat{Q}^{r}(T)$. For any $r$ let $\ga_{r}$ be the simple arrow associated to the step $\widehat{T}^{(r)}=\widehat{Q}(\widehat{T}^{(r-1)})$ of the algorithm. We obtain a sequence $\ga_{1},\ga_{2},..,\ga_{r},...$ of simple arrows. We denote $\ga(T,r)$ the concatenation $\ga_{1}...\ga_{r}$ of the first $r$ arrows in the sequence. We have $\ga(T,r)\prec\ga(T,r+1)$ with respect to the ordering $\prec$ in Definition \ref{def1ss1s2(algorithmR-V)}. Then we define $\ga(T,\infty)$ as the half infinite path in $\Pi(\cR)$ such that $\ga(T,r)\prec\ga(T,\infty)$ for all $r>0$.

Veech proved that $\widehat{Q}$ has an unique invariant measure which is absolutely continuous with respect to the Lebesgue measure, nevertheless this measure is not finite (see \cite{veech}). Zorich introduced an \emph{acceleration} of $\widehat{Q}$ with a finite invariant measure (see \cite{zorich}). For an i.e.t. $T$ without connections we define the integer $N(T)$ as the minimum of those $r\in\NN$ such that the type of $T$ is different from the type of $\widehat{Q}^{r}(T)$.

\begin{defi}\label{defaccelerazionezorich}
The \emph{Zorich's acceleration} is the map $\cZ:\De^{(1)}(\cR)\to \De^{(1)}(\cR)$ defined by $\cZ(T):=\widehat{Q}^{N(T)}(T)$.
\end{defi}

We recall the following important result in the ergodic theory of i.e.t.s (\cite{zorich}).

\begin{thmnonnum}[Zorich]
The map $\cZ$ in Definition \ref{defaccelerazionezorich} has an unique invariant measure $\mu$ which is absolutely continuous with respect to the Lebesgue measure on $\De^{(1)}(\cR)$. Moreover $\mu$ is ergodic.
\end{thmnonnum}

\subsection{Reduction of Rauzy classes}\label{Reduction of rauzy classes}

In this paragraph we describe \emph{reduction of Rauzy classes}, a combinatorial operation on Rauzy classes introduced in \cite{agy} generalizing a previous simpler version in \cite{av}. We follow closely \S 5 of \cite{agy}.

\subsubsection{Decorated Rauzy classes}\label{decoratedrauzyclasses}

Let $\cR$ be a Rauzy class with alphabet $\cA$ and $\cA'$ be a proper subset of $\cA$. An arrow is called $\cA'$-\emph{colored} if its winner belongs to $\cA'$. A path $\ga\in \Pi(\cR)$ is $\cA'$-colored if it is a concatenation of $\cA'$-colored arrows.

For an element $\pi$ in $\cR$ we say that $\pi$ is $\cA'$-\emph{trivial} if the last letters on both the top and the bottom rows of $\pi$ do not belong to $\cA'$, $\pi$ is $\cA'$-\emph{intermediate} if exactly one of those letters belongs to $\cA'$ and finally $\pi$ is $\cA'$-\emph{essential} if both letters belong to $\cA'$.
An $\cA'$-\emph{decorated Rauzy class} is a maximal subset $\cR_{\ast}$ of $\cR$ whose elements can be joined by an $\cA'$-colored path. Let $\Pi_{\ast}(\cR_{\ast})$ be the set of $\cA'$-colored paths starting and ending at permutations in $\cR_{\ast}$.

A decorated Rauzy class is called \emph{trivial} if it contains a trivial element $\pi$, in this case $\cR_{\ast}=\{\pi\}$ and $\Pi_{\ast}(\cR_{\ast})=\{\pi\}$, recalling that vertices are identified with zero-length paths. A decorated Rauzy class is called \emph{essential} if it contains an essential element. Admissibility implies that any essential decorated Rauzy class contains intermediate elements.

Let $\cR_{\ast}$ be an essential decorated Rauzy class and let $\cR_{\ast}^{ess}\subset\cR_{\ast}$ be the subset of essential elements. Let $\Pi^{ess}_{\ast}(\cR_{\ast})$ be the set of paths in $\Pi_{\ast}(\cR_{\ast})$ starting and ending at elements of $\cR_{\ast}^{ess}$. An \emph{arc} is a minimal non-trivial path in $\cR_{\ast}^{ess}$. In general an arc is concatenation of more that one arrow, anyway all arrows in the same arc are of the same type and have the same winner, so winner and type of an arc are well defined. The losers in an arc are all distinct, moreover the first loser is in $\cA'$ and the others are not. Any element in $\cR_{\ast}^{ess}$ is the starting point of a top and of a bottom arc and also the ending point of a top and a bottom arc.

If $\ga\in \Pi_{\ast}(\cR_{\ast})$ is an arrow then there exist unique paths $\ga_{s}$ and $\ga_{e}$ in $\Pi_{\ast}(\cR_{\ast})$ such that $\ga_{s}\ga\ga_{e}$ is an arc, called the \emph{completion} of $\ga$. If $\pi$ is intermediate the completion of the $\cA'$-colored arrow starting at $\pi$ is the only arc passing through $\pi$.

If $\pi\in \cR_{\ast}$ we define $\pi^{ess}$ as follows. If $\pi$ is essential then $\pi^{ess}=\pi$, if $\pi$ is intermediate let $\pi^{ess}$ be the end of the arc passing through $\pi$.

To $\ga\in\Pi_{\ast}(\cR_{\ast})$ we associate an element $\ga^{ess}\in \Pi^{ess}_{\ast}(\cR_{\ast})$ as follows. For a trivial path $\pi\in \cR_{\ast}$ we use the previous definition of $\pi^{ess}$. Assuming that $\ga$ is an arrow we distinguish two cases:
\begin{enumerate}
\item If $\ga$ starts at an essential element, we let $\ga^{ess}$ be the completion of $\ga$.
\item Otherwise, we let $\ga^{ess}$ be the endpoint of the completion of $\ga$.
\end{enumerate}
We extend the definition to paths $\ga\in \Pi_{\ast}(\cR_{\ast})$ by concatenation. Notice that if $\ga \in \Pi^{ess}_{\ast}(\cR_{\ast})$ then $\ga^{ess}=\ga$.

\subsubsection{Reduction of Rauzy classes}\label{reductionofrauzyclasses}

Given a permutation $\pi$ on the alphabet $\cA$, even not admissible, whose top and bottom rows end with different letters, we obtain the \emph{admissible end} of $\pi$ by deleting as many letters from the top and bottom rows of $\pi$ as necessary to obtain an admissible permutation. The resulting permutation belongs to some Rauzy class $\cR''$ on some alphabet $\cA''\subset \cA$.

Let $\cR_{\ast}$ be an essential decorated Rauzy class and $\pi\in \cR_{\ast}^{ess}$. Delete all the letters not belonging to $\cA'$ from the top and bottom rows of $\pi$. The resulting permutation $\pi'$ is not necessarily admissible, but since $\pi$ is essential the letters in the end of the top and bottom rows of $\pi'$ are distinct. Let $\pi^{red}$ be the admissible end of $\pi'$. We call $\pi^{red}$ the \emph{reduction} of $\pi$. We extend the operation of reduction from $\cR_{\ast}^{ess}$ to $\cR_{\ast}$ defining the reduction of $\pi\in \cR_{\ast}$ as the reduction of $\pi^{ess}$.

If $\ga\in\Pi^{ess}_{\ast}(\cR_{\ast})$ is an arc starting at $\pi_{s}$ and ending in $\pi_{e}$, then the reductions of $\pi_{s}$ and $\pi_{e}$ belong to the same Rauzy class and we define $\ga^{red}$ as the arrow which joins $\pi_{e}^{red}$ with $\pi_{s}^{red}$. The reduced arrow $\ga^{red}$ has the same type and the same winner of the arc $\ga$ and its loser is the first loser of $\ga$. It follows that the set of reductions of all elements $\pi$ in $\cR_{\ast}$ is a Rauzy class $\cR^{red}$ on some alphabet $\cA''\subset \cA'\subset\cA$. We define the reduction of a path $\ga\in \Pi_{\ast}(\cR_{\ast})$ as follows. If $\ga$ is a trivial (zero-length) path or an arc, it is defined as above. We extend the definition to the case $\ga\in \Pi^{ess}_{\ast}(\cR_{\ast})$ by concatenation. In general we let the reduction of $\ga$ to be equal to the reduction of $\ga^{ess}$. Restricted to essential elements the operation of reduction give a bijection $red:\cR^{ess}_{\ast}\to\cR^{red}$. If we think to elements $\pi\in\cR$ as trivial paths we can extend the previous operation to a bijection $red:\Pi^{ess}_{\ast}(\cR_{\ast})\to\Pi(\cR^{red})$ compatible with concatenation on the set of arcs.

\subsubsection{Drift in essential decorated Rauzy classes}\label{driftinessentialdecoratedrauzyclasses}
Let $\cR_{\ast}\subset \cR$ be an essential $\cA'$-decorated Rauzy class. For $\pi\in \cR_{\ast}$ let $\al_{t}(\pi)$ (respectively $\al_{b}(\pi)$) be the rightmost letter in the top (respectively in the bottom) row of $\pi$ that belongs to $\cA\setminus\cA'$. Let $d_{t}(\pi)$ (respectively $d_{b}(\pi)$) be the position of $\al_{t}(\pi)$ (respectively of $\al_{b}(\pi)$). Let $d(\pi):=d_{t}(\pi)+d_{b}(\pi)$. An essential element of $\cR_{\ast}$ is thus some $\pi$ such that $d_{t}(\pi)<d$ and $d_{b}(\pi)<d$. If $\pi_{s}$ is an essential element of $\cR_{\ast}$ and $\ga$ is an arrow starting at $\pi_{s}$ and ending at $\pi_{e}$ then
\begin{enumerate}
\item
$d_{t}(\pi_{e})=d_{t}(\pi_{s})$ or $d_{t}(\pi_{e})=d_{t}(\pi_{s})+1$, the second possibility happening iff $\ga$ is a bottom arrow whose winner precedes $\al_{t}(\pi_{s})$ in the top of $\pi_{s}$.
\item
$d_{b}(\pi_{e})=d_{b}(\pi_{s})$ or $d_{b}(\pi_{e})=d_{b}(\pi_{s})+1$, the second possibility happening iff $\ga$ is a top arrow whose winner precedes $\al_{b}(\pi_{s})$ in the bottom of $\pi_{s}$.
\end{enumerate}
In particular $d(\pi_{e})=d(\pi_{s})$ or $d(\pi_{e})=d(\pi_{s})+1$. In the second case we say that $\ga$ is \emph{drifting}. Let $\cR^{red}$ be the reduction of $\cR_{\ast}$ and let $\cA''\subset \cA'\subset \cA$ be the alphabet of $\cR^{red}$. If $\pi\in \cR_{\ast}$ is essential, then there exists some $\al\in\cA''$ which either precedes $\al_{t}(\pi)$ in the top row of $\pi$ or precedes $\al_{b}(\pi)$ in the bottom row of $\pi$, we call such an $\al$ \emph{good}. Indeed, if $\ga\in \Pi_{\ast}(\cR_{\ast})$ is a path starting at $\pi$, ending with a drifting arrow and minimal with this property, then the winner of the last arrow of $\ga$ belongs to $\cA''$ and either precedes $\al_{t}(\pi)$ in the top of $\pi$ (if the drifting arrow is a bottom) or precedes $\al_{b}(\pi)$ in the bottom of $\pi$ (if the drifting arrow is a top).

Note that if $\ga\in \Pi^{ess}_{\ast}(\cR_{\ast})$ is an arrow starting and ending at essential elements $\pi_{s},\pi_{e}$ then a good letter for $\pi_{s}$ is also a good letter for $\pi_{e}$. Moreover, if $\ga$ is not drifting then the winner of $\ga$ is not a good letter for $\pi_{s}$.

\subsubsection{Standard decomposition of separated paths}\label{standarddecompositionofseparatedpaths}
An arrow is called $(\cA\setminus\cA')$-\emph{separated} if both its winner and its loser belong to $\cA'$. A path $\ga\in \Pi_{\ast}(\cR_{\ast})$ is $(\cA\setminus\cA')$-separated if it is a concatenation of $(\cA\setminus\cA')$-separated arrows. We also say that a Rauzy path $\ga$ is \emph{complete} (or $\cA$-\emph{complete}) if for any letter $\al\in\cA$ there exists an arrow composing $\ga$ having $\al$ as winner.

If $\ga\in \Pi(\cR)$ is a non-trivial maximal $(\cA\setminus\cA')$-separated path then there exists an essential $\cA'$-decorated Rauzy class $\cR_{\ast}\subset \cR$ such that $\ga\in \Pi_{\ast}(\cR_{\ast})$. Moreover, if $\ga=\ga_{1}...\ga_{n}$ then any arrow $\ga_{i}$ starts at an essential element $\pi_{i}\in \cR^{ess}_{\ast}$ (and $\ga_{n}$ ends at an intermediate element of $\cR_{\ast}$ by maximality).

\begin{rem}\label{rem1standarddecompositionofseparatedpaths}
Let $r:=d(\pi_{n})-d(\pi_{1})$. Let $\ga=\ga^{(1)}\ga^{1}...\ga^{(r)}\ga^{r}$, where the $\ga^{i}$ are drifting arrows and $\ga^{(i)}$ are (possibly trivial) concatenation of non drifting arrows. If $\al$ is a good letter for $\pi_{1}$, then it follows that $\al$ is not the winner of any arrow in any $\ga^{(i)}$. The reduction of any $\ga^{(i)}$ are therefore non-complete paths in $\Pi(\cR^{red})$.
\end{rem}

\subsection{Lebesgue measure and distortion estimate}\label{thedistortionestimate}

Let $\pi$ be an admissible combinatorial datum in some Rauzy class $\cR$ over the alphabet $\cA$ and consider the associated simplex $\De^{(1)}_{\pi}$ of those $T\in\De_{\pi}$ with $\|\la\|=1$. We call $\leb_{d-1}$ the Lebesgue measure on $\De^{(1)}_{\pi}$ normalized in order to give measure one to it. For any finite path $\ga\in\Pi(\cR)$ starting at $\pi$ we define a sub-simplex of $\De^{(1)}_{\pi}$ by
$$
\De^{(1)}_{\ga}:=\De_{\ga}\cap \De^{(1)}_{\pi}.
$$
Modulo identifying $\De^{(1)}_{\pi}$ with the standard simplex $\De^{(1)}:=\{\la\in\RR^{d}_{+};\|\la\|=1\}$, the vertices of $\De^{(1)}_{\ga}$ are the vectors $(1/q^{\ga}_{\xi})^{t}B_{\ga}e_{\xi}$ with $\xi\in\cA$. Since $B_{\ga}$ belongs to $\sldz$ for any $\ga\in\Pi(\cR)$, we have the formula (see equation 5.5 in \cite{veech2})
\begin{equation}\label{eqmisurasimplessi}
\leb_{d-1}(\De^{(1)}_{\ga})=
\prod_{\xi\in\cA}(q^{\ga}_{\xi})^{-1}.
\end{equation}
Let $\Ga$ be a disjoint family (see Definition \ref{def1ss1s2(algorithmR-V)}) of Rauzy paths starting at $\pi$. Disjointness of $\Ga$ means that the simplices $\De^{(1)}_{\ga}$ with $\ga\in\Ga$ are disjoint. In this case we have
$$
\leb_{d-1}
(\bigcup_{\ga\in\Ga}\De^{(1)}_{\ga})=
\sum_{\ga\in\Ga}
\leb_{d-1}(\De^{(1)}_{\ga}).
$$

\subsubsection{Distortion estimate}\label{sss1ss3s2(distortionestimate)}

Fix $\pi$ in some Rauzy class $\cR$ and consider a finite path $\ga$ in $\Pi(\cR)$ starting at $\pi$. In view of Lemma \ref{lem1s2ss1sss2} we can interpret $\leb_{d-1}(\De^{(1)}_{\ga})$ as the probability that $\ga\prec\ga(T,\infty)$ for $T\in\De^{(1)}_{\pi}$. We introduce the notation
$$
\PP(\ga):=
\leb_{d-1}(\De^{(1)}_{\ga}).
$$
Fix a Rauzy path $\nu$ starting at $\pi'$ and ending in $\pi$. Let $r$ be the length of $\nu$ and set $\widehat{T}^{(r)}=\widehat{Q}^{r}(T)$ for $T\in\De^{(1)}_{\nu}$. For any $\ga$ starting at $\pi$ the concatenation $\nu\ga$ is defined. As before we can think to $\leb_{d-1}(\De^{(1)}_{\nu\ga})$ as the probability that $\nu\ga\prec\ga(T,\infty)$ for $T\in\De^{(1)}_{\pi'}$. Therefore the ratio
\begin{equation}\label{eqprobabilitacondizionatasimplessi}
P_{\nu}(\De^{(1)}_{\ga}):=
\frac
{\leb_{d-1}(\De^{(1)}_{\nu\ga})}
{\leb_{d-1}(\De^{(1)}_{\nu})}=
\frac
{\prod_{\xi\in\cA}q^{\nu}_{\xi}}
{\prod_{\xi\in\cA}q^{\nu\ga}_{\xi}}
\end{equation}
can be considered as the probability that $\ga\prec\ga(\widehat{T}^{(r)},\infty)$ for $\widehat{T}^{(r)}\in\De^{(1)}_{\pi}$ given that $\nu\prec\ga(T,\infty)$. When $\ga$ varies among all Rauzy path $\ga$ starting at $\pi$, equation (\ref{eqprobabilitacondizionatasimplessi}) defines a probability measure $P_{\nu}$ on Borel sets of $\De^{(1)}_{\pi}$. We also introduce the notation
$$
\PP_{\nu}(\ga):=P_{\nu}(\De^{(1)}_{\ga}).
$$
If $\Ga$ is a family of Rauzy paths starting at $\pi$ we write $\PP_{\nu}(\Ga):=P_{\nu}\big(\bigcup_{\ga\in\Ga}\De^{(1)}_{\ga}\big)$. In particular, if $\Ga$ is a disjoint family, we have
$$
\PP_{\nu}(\Ga)=\sum_{\ga\in\Ga}\PP_{\nu}(\ga).
$$
Observe that $B_{\nu\ga}=B_{\ga}B_{\nu}$ and thus $q^{\nu\ga}=B_{\ga}q^{\nu}$. For any vector $q$ in $\RR^{\cA}_{+}$ set $N(q):=\prod _{\xi\in\cA}q_{\xi}$ and define
$\PP_{q}(\ga):=
\frac{N(q)}
{N(B_{\ga}q)}
$.
For $\cA'\subset\cA$ and $q\in\RR^{\cA}_{+}$ set $M_{\cA'}(q):=\max_{\xi\in\cA'}q_{\xi}$. In the trivial case $\cA'=\cA$ we simply write $M(q):=M_{\cA}(q)$. Denote $\Pi_{\pi}(\cR)$ the set of those $\ga\in\Pi(\cR)$ starting at $\pi$. Theorem 5.4 in \cite{agy} gives the following distortion estimate.

\begin{thmnonnum}[Avila-Gouezel-Yoccoz]
There are two constants $C>0$ and $\theta>1$, depending only on the number of intervals $d$, with the following property. Let $\cA'\subset\cA$ be a non-empty proper subset, $m$ and $M$ be integers with $0\leq m\leq M$ and $q$ be any vector in $\RR^{\cA}_{+}$. Then we have
$$
\PP_{q}
\{\ga\in\Pi_{\pi}(\cR);
M(B_{\ga}q)>2^{M}M(q),
M_{\cA'}(B_{\ga}q)<2^{M-m}M(q)\}
\leq C(m+1)^{\theta}2^{-m}.
$$
\end{thmnonnum}
We also recall the following estimate (Proposition 5.9 in \cite{agy}).

\begin{propnonnum}
There are two constants $C>0$ and $\theta>1$, depending only on the number of intervals $d$, such that for any proper subset $\cA'\subset\cA$, any $M\in\NN$ and any $q$ in $\RR^{\cA}_{+}$ we have
$$
\PP_{q}
\{\ga\in\Pi_{\pi}(\cR);
\ga\textrm{ is not complete };
M(B_{\ga}q)>2^{M}M(q)\}
\leq C(M+1)^{\theta}2^{-M}.
$$
\end{propnonnum}

Fix any finite Rauzy path $\nu$ ending in $\pi$. Setting $q=q^{\nu}$ in the theorem and the proposition above, we get

\begin{equation}\label{eq1teo1thedistorsionestimate}
\PP_{\nu}
\{\ga\in\Pi_{\pi}(\cR);
M(B_{\ga}q^{\nu})>2^{M}M(q^{\nu}),
M_{\cA'}(B_{\ga}q^{\nu})<2^{M-m}M(q^{\nu})\}
\leq C(m+1)^{\theta}2^{-m}.
\end{equation}

\begin{equation}\label{eq2teo1thedistorsionestimate}
\PP_{\nu}
\{\ga\in\Pi_{\pi}(\cR);
\ga\textrm{ is not complete };
M(B_{\ga}q^{\nu})>2^{M}M(q^{\nu})\}
\leq C(M+1)^{\theta}2^{-M}.
\end{equation}

\section{Reduced Triples}\label{s3}

This section is devoted to the proof of Theorem \ref{teoremaf}. We consider the non-normalized version $Q$ of the Rauzy-Veech algorithm and for $T$ without connections we write $T^{(r)}=Q^{r}(T)$. Following the notation introduced in \S \ref{Normalized Rauzy Veech algorithm and Zorich's acceleration}, we denote $\ga(T,\infty)$ the half-infinite Rauzy path in the Rauzy diagram generated by $T$ and $\ga(T,r)$ the concatenation of the first $r$ arrows of $\ga(T,\infty)$.

\begin{defi}\label{defdetectingproducing}
Consider $T$ without connections and a triple $(\be,\al,n)$ reduced for $T$. We say that the triple $(\be,\al,n)$ \emph{is detected} by $Q$ if there exists $r$ such that the singularities $u^{(r),b}_{\be}$ and $u^{(r),t}_{\al}$ satisfy
$$
|T^{n}u^{b}_{\be}-u^{t}_{\al}|=
|u^{(r),b}_{\be}-u^{(r),t}_{\al}|.
$$

We say that $(\be,\al,n)$ \emph{is produced} by $Q$ if there exists some $r$ and a letter $\xi$ in $\cA$ such that
$$
|T^{n}u^{b}_{\be}-u^{t}_{\al}|=\la^{(r)}_{\xi}.
$$
\end{defi}

\subsection{Detection of reduced triples with the algorithm}\label{s3ss2}

In this subsection we prove the first part of Theorem \ref{teoremaf}.

\begin{lem}\label{lem1s3ss1sss0}
Consider $T$ without connections. If the triple $(\be,\al,n)$ is reduced for $T$, then is detected by the algorithm at some step $T^{(r)}$.
\end{lem}

\begin{demo}
Suppose that $T^{n}u^{b}_{\be}<u^{t}_{\al}$, the other case being symmetric. Consider $m\in\{0,\dots,n\}$ such that $T^{n-m}u^{b}_{\be}=
\min\{T^{n-i}u^{b}_{\be};i\in\{0,\dots,n\}\}$, that is the leftmost among the points $T^{n-i}u^{b}_{\be}$ with $0\leq i\leq n$. Since the triple $(\be,\al,n)$ is reduced for $T$, then for the same $m$ we have $T^{-m}u^{t}_{\al}=\min\{T^{i}u^{t}_{\al};i\in\{0,..,n\}\}$. Therefore
$$
|T^{n}u^{b}_{\be}-u^{t}_{\al}|=
|T^{n-m}u^{b}_{\be}-T^{-m}u^{t}_{\al}|
$$
and in particular $T^{n-m}u^{b}_{\be}<T^{-m}u^{t}_{\al}$.

For $j\in\NN$ let $I^{(j)}$ be the interval where $T^{(j)}=Q^{j}(T)$ acts. The left endpoint of any $I^{(j)}$ is the left endpoint of $I$, the right endpoint of $I^{(j)}$ is the rightmost singularity of $T^{(j-1)}$. Consider $r$ defined by
$$
r:=\max
\{j\in\NN;
T^{n-m}u^{b}_{\be}\in I^{(j)}
\textrm{ and }
T^{-m}u^{t}_{\al}\in I^{(j)}\}.
$$
Since $T^{n-m}u^{b}_{\be}<T^{-m}u^{t}_{\al}$ then maximality of $r$ implies $T^{-m}u^{t}_{\al}\not\in I^{(r+1)}$, therefore $u^{(r),b}_{\xi}\leq T^{-m}u^{t}_{\al}$ and $u^{(r),t}_{\xi}\leq T^{-m}u^{t}_{\al}$ for any $\xi\in\cA$. On the other hand $T^{(r)}$ is the first return map of $T$ to $I^{(r)}$, thus as it is explained in \S \ref{S2ss1sss3}, there are non-negative integers $l$ and $h$ such that
$$
T^{n-m}u^{b}_{\be}=(T^{(r)})^{l}(u^{(r),b}_{\be})
\textrm{ and }
T^{-m}u^{t}_{\al}=(T^{(r)})^{-h}(u^{(r),t}_{\al}).
$$

We prove that $h=0$. If $h>0$ then $u^{(r),t}_{\al}=T^{-h'}u^{t}_{\al}$ with $0\leq h'<m$. Since $T^{-m}u^{t}_{\al}$ is the leftmost among the points $T^{-i}u^{t}_{\al}$ with $0\leq i\leq n$ then the last condition implies $T^{-m}u^{t}_{\al}<u^{(r),t}_{\al}$, which is absurd by the discussion above. In particular $u^{(r),t}_{\al}$ is the rightmost singularity of $T^{(r)}$.

Now we prove that $l=0$. If $l>0$ then $u^{(r),b}_{\be}=T^{l'}u^{b}_{\be}$ with $0\leq l'<n-m$. Since $T^{n-m}u^{b}_{\be}$ is the leftmost among the points $T^{n-i}u^{b}_{\be}$ with $0\leq i\leq n$ then the last condition implies $T^{n-m}u^{b}_{\be}<u^{(r),b}_{\be}$. Moreover we proved that $u^{(r),t}_{\al}$ is the rightmost singularity of $T^{(r)}$, thus $u^{(r),b}_{\be}<u^{(r),t}_{\al}$. Combining these two conditions we get
$$
T^{n-m}u^{b}_{\be}<
T^{l'}u^{b}_{\be}<
T^{-m}u^{t}_{\al}
$$
with $l'-(-m)=l'+m<n$, which is absurd because $(\be,\al,n)$ is reduced for $T$. The lemma is proved.
\end{demo}

\subsubsection{Minimal detecting paths}\label{s3ss1sss1}

We introduce some notation. Let $\pi$ be an admissible combinatorial datum in the same Rauzy class of $\pi_{0}$. For any $\xi$ in $\cA$ denote $e_{\xi}$ the corresponding vector of the standard basis of $\RR^{d}$, in order to have $\la=\sum_{\xi\in\cA}\la_{\xi}e_{\xi}$ for any $\la\in\RR^{d}$. For a pair $(\be,\al)$ such that $\pi^{b}(\be)>1$ and $\pi^{t}(\al)>1$ define the integer vectors
$$
w^{b}_{\be}(\pi):=
\sum_{\pi^{b}(\xi)<\pi^{b}(\be)}e_{\xi}
\textrm{ and }
w^{t}_{\al}(\pi):=
\sum_{\pi^{t}(\xi)<\pi^{t}(\al)}e_{\xi},
$$
then set $w_{\be,\al}(\pi):=w^{b}_{\be}(\pi)-w^{t}_{\al}(\pi)$. Thus the singularities of an i.e.t. $T=(\pi,\la)$ in $\De_{\pi}$ are given by $u^{b}_{\be}=\langle w^{b}_{\be}(\pi),\la\rangle$ and $u^{t}_{\al}=\langle w^{t}_{\al}(\pi),\la\rangle$ and therefore
$$
u^{b}_{\be}-u^{t}_{\al}=
\langle w_{\be,\al}(\pi),\la\rangle.
$$
In particular the set of those $T$ in $\De^{(1)}_{\pi}$ with $u^{b}_{\be}=u^{t}_{\al}$ coincides with the intersection $\De_{\pi}\cap(w_{\be,\al}(\pi))^{\perp}$, where $(w_{\be,\al}(\pi))^{\perp}$ denotes the hyperplane normal to $w_{\be,\al}(\pi)$.\\

Consider $T$ without connections and a triple $(\be,\al,n)$ reduced for $T$. According to Lemma \ref{lem1s3ss1sss0}, let $r$ in $\NN$ such that the singularities $u^{(r),t}_{\al}$ and $u^{(r),b}_{\be}$ of $T^{(r)}$ satisfy
$$
|T^{n}u^{b}_{\be}-u^{t}_{\al}|=
|u^{(r),t}_{\al}-u^{(r),b}_{\be}|.
$$
In general there exist more than one value of $r$ satisfying the last condition and we denote $r_{min}$ the minimal one. Then we denote $\ga(\be,\al,n,T):=\ga(T,r_{min})$ the \emph{minimal detecting path} for $T$ and the triple $(\be,\al,n)$.

\begin{defi}\label{def1s3ss1sss1}
Fix $n$ in $\NN$ and a pair $(\be,\al)$ with $\pi^{b}(\be)>1$ and $\pi^{t}(\al)>1$. Denote $\Ga(\be,\al,n)$ the family of minimal detecting paths $\ga=\ga(\be,\al,n,T)$, where $T$ varies among the elements without connections in $\De_{\pi_{0}}$ such that the triple $(\be,\al,n)$ is reduced for $T$. For any $\ga$ in $\Ga(\be,\al,n)$ denote $\De^{\ast}_{\ga}$ the set of those $T$ in $\De_{\ga}$ such that $(\be,\al,n)$ is reduced for $T$.
\end{defi}

\begin{rem}\label{rem1s3ss1sss1}
For any Rauzy path $\ga$ with length $r$ and any $T\in\De_{\ga}$, the integers $l(\be,r)$ and $h(\al,r)$ defined in \S \ref{S2ss1sss3} with the property $u^{(r),b}_{\be}=T^{l(\be,r)}u^{b}_{\be}$ and $u^{(r),t}_{\al}=T^{-h(\al,r)}u^{t}_{\al}$ depend only on $\ga$. Thus $\Ga(\pi_{0},\be,\al,n)$ is a disjoint family, since its elements are minimal by definition.
\end{rem}

\begin{prop}\label{prop1s3ss1sss1}
Fix $n$ in $\NN$ and a pair $(\be,\al)$ with $\pi^{b}(\be)>1$ and $\pi^{t}(\al)>1$. For any $T$ in $\De_{\pi_{0}}$ such that the triple $(\be,\al,n)$ is reduced for $T$ there exists an unique path $\ga$ in $\Ga(\be,\al,n)$ with $T\in\De_{\ga}$ and such that
$$
|T^{n}u^{b}_{\be}-u^{t}_{\al}|=
|\langle\la,B_{\ga}^{-1}w_{\be,\al}(\pi)\rangle|,
$$
where $\pi$ is the ending point of $\ga$. On the other hand for any $\ga$ in $\Ga(\be,\al,n)$ the set $\De^{\ast}_{\ga}$ is an open subset of $\De_{\ga}$ whose closure contains $\De_{\ga}\cap(w_{\be,\al}(\pi))^{\perp}$.
\end{prop}

\begin{demo}
For $T$ as in the statement we set $\ga:=\ga(\be,\al,n,T)$. Such $\ga$ is an element of $\Ga(\be,\al,n)$ and we trivially have $T\in\De_{\ga}$, moreover $\ga$ is unique because $\Ga(\be,\al,n)$ is a disjoint family. Let $r$ be the length of $\ga$ and $\pi$ be its ending point. According to Lemma \ref{lem1s3ss1sss0} we have $|T^{n}u^{b}_{\be}-u^{t}_{\al}|=|u^{(r),t}_{\al}-u^{(r),b}_{\be}|$, where $u^{(r),t}_{\al}$ and $u^{(r),b}_{\be}$ are the singularities of $T^{(r)}$. Moreover $T^{(r)}$ belongs to $\De_{\pi}$, thus for $\la^{(r)}=^{t}B_{\ga}^{-1}\la$ we have
$$
|u^{(r),t}_{\al}-u^{(r),b}_{\be}|=
|\langle\la^{(r)},w_{\be,\al}(\pi)\rangle|,
$$
which implies the first part of the proposition. Consider $\ga$ in $\Ga(\be,\al,n)$ with length $r$ and ending in $\pi$. The map $T\mapsto T^{(r)}$ is an homeomorphism between $\De_{\ga}$ and $\De_{\pi}$. Let $I^{(r)}$ be the interval where $T^{(r)}$ acts. The triple $(\be,\al,n)$ is reduced for $T$ if and only if the singularities $u^{(r),t}_{\al}$ and $u^{(r),b}_{\be}$ are in some specified order in $I^{(r)}$ with respect to the other singularities $u^{(r),t}_{\xi}$ and $u^{(r),b}_{\xi}$ with $\xi\in\cA$. We do not compute explicitly the required order relations, anyway they depend uniquely on $\ga$ and $\pi$ and obviously determine an open condition for $T^{(r)}\in\De_{\pi}$. Hence having $(\be,\al,n)$ as reduced triples is an open condition for $T\in\De_{\ga}$. The last part of the proposition follows observing that for $T\in\De_{\ga}$ condition $T\to (w_{\be,\al}(\pi))^{\perp}$ is equivalent to $|T^{n}u^{b}_{\be}-u^{t}_{\al}|\to 0$. The proposition is proved.
\end{demo}

\subsection{Measure estimate for $\Ga(\be,\al,n)$}\label{s3ss2}

\subsubsection{Combinatorial properties of the family $\Ga(\be,\al,n)$}\label{s3ss2sss1}

\begin{lem}\label{lem1s3ss2sss1}
For any $\ga$ in $\Ga(\be,\al,n)$ the last arrow $\ga_{last}$ of $\ga$ is either a top arrow with loser $\be$ or a bottom arrow with loser $\al$. Moreover we have $\max\{q^{\ga}_{\al},q^{\ga}_{\be}\}>n/2$. Finally if $W=W(\ga)$ is the winner of $\ga_{last}$, we have $q^{\ga}_{W}<n$.
\end{lem}

\begin{demo}
Consider $\ga$ in $\Ga(\be,\al,n)$ and $T$ in $\De_{\ga}$. Let $r$ be the length of $\ga$. For the singularities of $T^{(r)}$ we have $|T^{n}u^{b}_{\be}-u^{t}_{\al}|=|u^{(r),t}_{\al}-u^{(r),b}_{\be}|$ and $r$ is the smallest integer such that this last condition holds, since $\ga$ is minimal by definition. Hence we have either $u^{(r),b}_{\be}=T^{(r-1)}(u^{(r-1),b}_{\be})$, that is $\ga_{last}$ is a top arrow with loser $\be$, or $u^{(r),t}_{\al}=(T^{(r-1)})^{-1}(u^{(r-1),t}_{\al})$, that is $\ga_{last}$ is a bottom arrow with winner $\al$. Recall from \S \ref{S2ss1sss3} that for any $\xi$ in $\cA$ and any $r$ we have integers $l(\xi,r)$ and $h(\xi,r)$ satisfying $q^{(r)}_{\xi}=l(\xi,r)+h(\xi,r)+1$ and such that $u^{(r),b}_{\xi}=T^{l(\xi,r)}u^{b}_{\be}$ and $u^{(r),t}_{\xi}=T^{-h(\xi,r)}u^{t}_{\al}$. Since $l(\be,r)+h(\al,r)=n$ then the second part of the statement follows.
Finally observe that if $\be$ loses against $W$ in $\ga_{last}$ then $u^{(r-1),b}_{\be}\in I^{(r-1),t}_{W}$ and $l(\be,r)=l(\be,r-1)+q^{(r-1)}_{W}$. Similarly if $\al$ loses in $\ga_{last}$ then $u^{(r-1),t}_{\al}\in I^{(r-1),b}_{W}$ and $h(\al,r)=h(\al,r-1)+q^{(r-1)}_{W}$. Moreover in both cases we have $q^{(r)}_{W}=q^{(r-1)}_{W}$, thus the third part of the statement follows. The lemma is proved.
\end{demo}

\begin{lem}\label{lem2s3ss2sss1}
Consider $\ga$ in $\Ga(\be,\al,n)$ and let $\pi$ be its ending point. There exist a pair of letters $(W_{+},W_{-})$ such that $\langle w_{\be,\al}(\pi),e_{W_{+}}\rangle=1$ and $\langle w_{\be,\al}(\pi),e_{W_{-}}\rangle=-1$, that is
$$
\pi=
\left(
\begin{array}{cccccc}
 \dots &  W_{-} & \dots  &  \al  &  \dots  &   W_{+}  \\
 \dots &  W_{+} & \dots  &  \be  &  \dots  &   W_{-}  \\
\end{array}
\right),
$$
where the cases $W_{+}=\al$, $W_{-}=\be$ and $\be=\al$ are (separately) possible. In particular if $W$ is the last winner of $\ga$ we have $\langle w_{\be,\al}(\pi),e_{W}\rangle=\pm1$.
\end{lem}

\begin{demo}
Let $W$ be the last winner of $\ga$. According to Lemma \ref{lem1s3ss2sss1} the last loser in $\ga$ is either $\be$ or $\al$. We suppose that it is $\be$, the other case being symmetric. Thus we have
$$
\pi=
\left(
\begin{array}{cccc}
       &     & \dots &   W   \\
 \dots &  W  & \be   &  \dots\\
\end{array}
\right)
$$
and for $W_{+}:=W$ we have $\langle w_{\be,\al}(\pi),e_{W_{+}}\rangle=1$. In particular the last part of the lemma follows (if the last loser in $\ga$ is $\al$ we get $\langle w_{\be,\al}(\pi),e_{W}\rangle=-1$). Let $r$ be the length of $\ga$ and consider any $T$ in $\De_{\ga}$ having $(\be,\al,n)$ as reduced triple. For such $T$ we have $u^{(r),b}_{\be}=T^{l(\be,r)}u^{b}_{\be}$, $u^{(r),t}_{\al}=T^{-h(\al,r)}u^{t}_{\al}$ and $u^{(r),b}_{W}=T^{l(W,r)}u^{b}_{W}$, where $l(\be,r)$, $h(\al,r)$ and $l(W,r)$ are the integers introduced in \S \ref{S2ss1sss3}. Moreover the combinatorics of $\pi$ implies $u^{(r),b}_{W}<u^{(r),b}_{\be}$. Since $\ga_{last}$ is a top arrow with winner $W$ and loser $\be$, then we have $l(\be,r)=l(\be,r-1)+q^{(r-1)}_{W}$ and $l(W,r)=l(W,r-1)<q^{(r-1)}_{W}$, thus $l(W,r)+h(\al,r)<l(\be,r)+h(\al,r)=n$. It follows that $u^{(r),t}_{\al}>u^{(r),b}_{W}$, otherwise the triple $(\be,\al,n)$ is not reduced for $T$. Finally observe that $\pi^{t}(W)=d\geq \pi^{t}(\al)$ and $\pi^{b}(W)<\pi^{b}(\be)$, thus condition $u^{(r),t}_{\al}>u^{(r),b}_{W}$ implies that there exists a letter $W_{-}$ such that $\pi^{t}(W_{-})<\pi^{t}(\al)$ and $\pi^{b}(\be)\geq\pi^{b}(W_{-})$, which is equivalent to $\langle e_{W_{-}},w_{\be,\al}(\pi)\rangle=-1$. The lemma is proved.
\end{demo}

\subsubsection{Local description of condition $|T^{n}u^{b}_{\be}-u^{t}_{\al}|<\rho$}\label{s3ss2sss2}

Here we prove a local estimate for those $T$ such that $(\be,\al,n)$ is reduced for $T$ and satisfies $|T^{n}u^{b}_{\be}-u^{t}_{\al}|<\rho$. We need a finite measure setting, thus for $\ga\in \Ga(\be,\al,\epsilon)$ we consider the $(d-1)$-simplex $\De^{(1)}_{\ga}$ in $\De^{(1)}_{\pi_{0}}$. Lemma \ref{lem2s3ss2sss1} implies the following corollary.

\begin{cor}\label{cor1s3ss2sss2}
Consider $\ga$ in $\Ga(\be,\al,n)$ and let $\pi$ be its ending point. Then the intersection $\De^{(1)}_{\pi}\cap(w_{\be,\al}(\pi))^{\perp}$ is a non-empty convex subset of $\De^{(1)}_{\pi}$ with dimension $d-2$. Moreover, if $W$ is the last winner of $\ga$, then $e_{W}$ is a vertex of $\De^{(1)}_{\pi}$ which does not belong to $(w_{\be,\al}(\pi))^{\perp}$.
\end{cor}

\begin{demo}
According to Lemma \ref{lem2s3ss2sss1}, consider a pair of letters $(W_{+},W_{-})$ such that $\langle w_{\be,\al}(\pi),e_{W_{+}}\rangle=1$ and $\langle w_{\be,\al}(\pi),e_{W_{-}}\rangle=-1$. The first part of the statement follows observing that any $\xi\in\cA$ satisfies one of the following three conditions
\begin{itemize}
\item
$\langle w_{\be,\al}(\pi),e_{\xi}\rangle=0$, that is $e_{\xi}$ belongs to $(w_{\be,\al}(\pi))^{\perp}\cap\De^{(1)}_{\pi}$.
\item
$\langle w_{\be,\al}(\pi),e_{\xi}\rangle=1$, thus $\langle w_{\be,\al}(\pi),e_{\xi}+e_{W_{-}}\rangle=0$, that is $2^{-1}(e_{\xi}+e_{W_{-}})$ belongs to $(w_{\be,\al}(\pi))^{\perp}\cap\De^{(1)}_{\pi}$.
\item
$\langle w_{\be,\al}(\pi),e_{\xi}\rangle=-1$, thus $\langle w_{\be,\al}(\pi),e_{\xi}+e_{W_{+}}\rangle=0$, that is $2^{-1}(e_{\xi}+e_{W_{+}})$ belongs to $(w_{\be,\al}(\pi))^{\perp}\cap\De^{(1)}_{\pi}$.
\end{itemize}
In particular by Lemma \ref{lem1s3ss2sss1} we have $\langle w_{\be,\al}(\pi),e_{W}\rangle=\pm1$, thus $e_{W}$ does not belong to $(w_{\be,\al}(\pi))^{\perp}$. The corollary is proved.
\end{demo}

\begin{lem}\label{lem1s3ss2sss2}
Consider $\ga$ in $\Ga(\be,\al,n)$ and let $\pi$ be its ending point. Then the intersection $\De^{(1)}_{\ga}\cap(B_{\ga}^{-1}w_{\be,\al}(\pi))^{\perp}$ is a non-empty convex subset of $\De^{(1)}_{\ga}$ with dimension $d-2$. Moreover, if $W$ is the last winner of $\ga$, for any $\rho>0$ we have
$$
\frac
{\leb
\{T\in\De^{(1)}_{\ga};
|T^{n}u^{b}_{\be}-u^{t}_{\al}|<\rho\}}
{\leb(\De^{(1)}_{\ga})}
\leq
\rho q^{\ga}_{W}.
$$
\end{lem}

\begin{demo}
Recall that the vertexes of $\De^{(1)}_{\ga}$ are the vectors $v_{\xi}:=(1/q^{\ga}_{\xi})^{t}B_{\ga}e_{\xi}$ with $\xi\in\cA$. Observe that for any $\xi$ we have
$$
\langle B_{\ga}^{-1}w_{\be,\al}(\pi),^{t}B_{\ga}e_{\xi}\rangle=
\langle w_{\be,\al}(\pi),e_{\xi}\rangle,
$$
thus the first part of the statement follows from Corollary \ref{cor1s3ss2sss2}. To get the second part consider the linear form $\la\mapsto f_{\ga}(\la):=q^{\ga}_{W}\langle \la,B_{\ga}^{-1}w_{\be,\al}(\pi)\rangle$. The first part of the statement says that $\ker(f_{\ga})$ intersects $\De^{(1)}_{\ga}$ in its interior. Moreover we know from Corollary \ref{cor1s3ss2sss2} that $\langle w_{\be,\al}(\pi),e_{W}\rangle=\pm1$, thus $f_{\ga}(v_{W})=\pm1$. On the other hand the first part of Proposition \ref{prop1s3ss1sss1} implies that for any $T$ in $\De^{(1)}_{\ga}$ we have $
|T^{n}u^{b}_{\be}-u^{t}_{\al}|=
|\langle\la,B_{\ga}^{-1}w_{\be,\al}(\pi)\rangle|
$, thus $|T^{n}u^{b}_{\be}-u^{t}_{\al}|<\rho$ if and only if $|f_{\ga}(\la)|<\rho q^{\ga}_{W}$, therefore the second part of the statement follows. The lemma is proved.
\end{demo}

\subsubsection{Global estimate}\label{s3ss2sss3}

\begin{prop}\label{props3ss2}
There exists a positive constant $C$, depending only on the number of intervals $d$, such that for any pair of letters $(\be,\al)$ with $\pi^{b}_{0}(\be)>1$ and $\pi^{t}_{0}(\al)>1$ and any positive integer $N$ we have
$$
\sum_{2^{N-1}\leq n<2^{N}}
\bigg(
\sum_{\ga\in\Ga(\be,\al,n)}
q^{\ga}_{W(\ga)}\leb(\De^{(1)}_{\ga})
\bigg)
\leq C2^{N}.
$$
\end{prop}

\begin{demo}
For any $W\in\cA$ we call $\Ga(\be,\al,n;W)$ the sub-family of those paths $\ga\in\Ga(\be,\al,n)$ whose last winner is $W$. Since $\Ga(\be,\al,n)=\bigsqcup_{W\in\cA}\Ga(\be,\al,n;W)$ and $\cA$ is finite with $d$ elements, then it is enough to prove the statement replacing $\Ga(\be,\al,n)$ by $\Ga(\be,\al,n;W)$. If $C(W)$ is the constant that we get for $\Ga(\be,\al,n;W)$, then $C=\sum_{W\in\cA}C(W)$ works for $\Ga(\be,\al,n)$. To simplify notation, from now on we write $\Ga_{n}$ instead of $\Ga(\be,\al,n;W)$. Denote $\log$ the logarithm in base $2$. For any $k\in\{0,\dots,\log n\}$ we denote $\Ga_{n,k}$ the sub-family of those $\ga\in\Ga_{n}$ such that $2^{k}\leq q^{\ga}_{W}<2^{k+1}$. Observe that $\Ga_{n,k}=\emptyset$ for $k>\log n$, according to Lemma \ref{lem1s3ss2sss2}. The main argument in the proof is contained in the following lemma.

\begin{lem}\label{lem1s3ss2sss3}
For any positive integer $n$, any $k\in\{0,\dots,\log n\}$ and any $i\in\{0,\dots,2^{k}-1\}$ the families $\Ga_{(n+i),k}$ are each other disjoint.
\end{lem}

\begin{demo}
Let $m>n$ and consider $\ga\in\Ga_{n,k}$ and $\ga'\in\Ga_{m,k}$ such that $\ga\prec\ga'$. By definition $W$ is the last winner both of $\ga$ and $\ga'$ with $2^{k}\leq q^{\ga}<2^{k+1}$ and $2^{k}\leq q^{\ga'}<2^{k+1}$. Let $r$ and $r+j$ be the length respectively of $\ga$ and $\ga'$, where $j\geq1$. Recall from \S \ref{S2ss1sss3} that there are integers $l(\be,r)$ and $h(\al,r)$, depending only on $\ga$, such that for any $T\in\De^{(1)}_{\ga}$ we have $u^{(r),b}_{\be}=T^{l(\be,r)}u^{b}_{\be}$ and $u^{(r),t}_{\al}=T^{-h(\al,r)}u^{t}_{\al}$. Similarly let $l(\be,r+j)$ and $h(\al,r+j)$ be the integers corresponding to $\ga'$. Since $\ga\in\Ga_{n,k}$ then
$$
n=
l(\be,r)+h(\al,r)=
l(\be,r-1)+h(\al,r-1)+q^{\ga}_{W}.
$$
On the other hand $\ga'\in\Ga_{m,k}$, thus
$$
m=
l(\be,r+j)+h(\al,r+j)=
l(\be,r+j-1)+h(\al,r+j-1)+q^{\ga'}_{W}\geq
$$
$$
l(\be,r)+h(\al,r)+q^{\ga}_{W}=n+2^{k},
$$
where the last inequality holds since $\ga\prec\ga'$. Therefore condition $\ga\prec\ga'$ implies $m\geq n+2^{k}$ and since this is newer satisfied by $m=n+i$ with $0\leq i<2^{k}$ the lemma is proved.
\end{demo}

We recall that $\Ga(\be,\al,n)$ is a disjoint family, so any $\Ga_{n,k}$ is disjoint (since it is a sub-family of $\Ga(\be,\al,n)$), that is $\sum_{\ga\in\Ga_{n,k}}\PP(\ga)=\PP(\Ga_{n,k})$. By definition of $\Ga_{n,k}$ we have trivially
$$
\sum_{2^{N-1}\leq n<2^{N}}
\sum_{\ga\in\Ga_{n}}
q^{\ga}_{W(\ga)}\PP(\ga)
\leq
2\sum_{2^{N-1}\leq n<2^{N}}
\sum_{k=0}^{\log n}
2^{k}\PP(\Ga_{n,k}).
$$
Now we observe that if $2^{N-1}\leq n<2^{N}$ then any $\ga\in\Ga_{n,k}$ satisfies $2^{k}\leq q^{\ga}_{W}<n<2^{N}$ and hence $k<N$. We have the identity
$$
\sum_{2^{N-1}\leq n<2^{N}}
\sum_{k=0}^{\log n}
2^{k}\PP(\Ga_{n,k})=
\sum_{k=0}^{N-1}
2^{k}\sum_{2^{N-1}\leq n<2^{N}}
\PP(\Ga_{n,k}).
$$

Fix any $k\in\{0,\dots,N-1\}$ and any $i\in\{0,\dots,2^{N-k-1}-1\}$. According to Lemma \ref{lem1s3ss2sss3}, the families $\Ga_{n,k}$ with $2^{N-1}+i2^{k}\leq n<2^{N-1}+(i+1)2^{k}$ are disjoint. We set $\cG_{N,k,i}:=\bigsqcup_{n}\Ga_{n,k}$, where $n$ varies in $\{2^{N-1}+i2^{k},\dots,2^{N-1}+(i+1)2^{k}\}$. By disjointness $\PP(\cG_{N,k,i})=\sum_{n}\PP(\Ga_{n,k})$, where $n$ varies in the same interval, therefore we have the identity
$$
\sum_{2^{N-1}\leq n<2^{N}}
\PP(\Ga_{n,k})=
\sum_{i=0}^{2^{N-k-1}-1}
\PP(\cG_{N,k,i}).
$$
Observe that $\max\{q^{\ga}_{\be},q^{\ga}_{\al}\}>n/2\geq 2^{N-2}$ for any $\ga\in\Ga_{n,k}$ with $2^{N-1}\leq n<2^{N}$, according to Lemma \ref{lem1s3ss2sss2}. On the other hand we have $q^{\ga}_{W}<2^{k+1}$, thus the family $\cG_{N,k,i}$ is contained in $\{\ga;M(q^{\ga})>2^{N-k+3}m(q^{\ga})\}$. Equation  (\ref{eq1teo1thedistorsionestimate}) in the background implies that there exists two positive constants $C$ and $\te$, depending only on the number of intervals $d$, such that $\PP(\cG_{N,k,i})\leq C(N-k+3)^{\te}2^{-(N-k+3)}$. Modulo changing the constant $C$, we write
$$
\PP(\cG_{N,k,i})
\leq
C\frac{(N-k)^{\te}}{2^{N-k}}
$$
Applying the two identities above and the last inequality we get
$$
\sum_{2^{N-1}\leq n<2^{N}}
\sum_{k=0}^{\log n}
2^{k}\PP(\Ga_{n,k})
\leq
C
\sum_{k=0}^{N-1}
2^{k}\cdot2^{N-k-1}
\frac{(N-k)^{\te}}{2^{N-k}}=
2^{N-1}
C\sum_{m=1}^{N}
\frac{m^{\te}}{2^{m}}.
$$
It follows that
$$
\sum_{2^{N-1}\leq n<2^{N}}
\sum_{\ga\in\Ga_{n}}
q^{\ga}_{W(\ga)}\PP(\ga)
\leq
2^{N}C\sum_{m=1}^{+\infty}
\frac{m^{\te}}{2^{m}}
$$
and the statement in the proposition follows since $\sum_{m=1}^{+\infty}m^{\te}2^{-m}<+\infty$.
\end{demo}

\subsection{Producing reduced triples with the algorithm}\label{s3ss3}

The results in this paragraph hold for any $T$ in $\De_{\pi_{0}}$. We consider the non-normalized version $Q$ of the Rauzy-Veech algorithm. For any $T=(\pi_{0},\la)$ in $\De_{\pi_{0}}$ without connections we call $(\pi^{(r)},\la^{(r)})$ the pair of combinatorial and length data of $T^{(r)}=Q^{r}(T)$.

\subsubsection{Properties A and B}\label{s3ss3sss1}

We introduce two combinatorial properties for pairs $(\be,\al)$ with $\pi_{0}^{t}(\al)>1$ and $\pi_{0}^{b}(\be)>1$, depending only on the Rauzy class $\cR$ of $\pi_{0}$.

\begin{defi}\label{defproprietaAeB}
Let $(\be,\al)$ be an ordered pair of letters with $\pi_{0}^{t}(\al)>1$, $\pi_{0}^{b}(\be)>1$.
\begin{itemize}
\item
We say that $(\be,\al)$ satisfies \emph{property A} if there exists a combinatorial datum $\pi=\pi(\be,\al)$ in $\cR$ such that
$$
\pi^{t}(\al)=\pi^{b}(\be)=d
$$
that is we have
$$
\pi=  \left(
\begin{array}{cc}
\dots &   \al   \\
\dots &   \be   \\
\end{array}
\right).
$$
\item
We say that $(\be,\al)$ satisfies \emph{property B} if there exists a combinatorial datum $\pi=\pi(\be,\al)$ in $\cR$ and a letter $V$ in $\cA$ such that
$$
\begin{array}{ll}
\{\xi\in \cA\,;\,\pi^{t}(\xi)<\pi^{t}(\al)\}\cup \{V\}=
\{\xi\in \cA\,;\,\pi^{b}(\xi)<\pi^{b}(\be)\}\\
\pi^{t}(V)=\pi^{b}(\al)=d.
\end{array}
$$
When $(\be,\al)$ satisfies property B and $\pi$ is a combinatorial datum as above we call $L$ the letter such that $\pi^{b}(L)=\pi^{b}(\be)-1$ and $\pi^{t}(L)<\pi^{t}(\al)$. Admissibility of $\pi$ implies $L\neq V$, thus we have
$$
\pi= \left(
\begin{array}{ccccccc}
\dots &   L    &  \dots  &   \al  &  \dots  &        &  V  \\
\dots &   V    &  \dots  &    L   &   \be   &  \dots & \al \\
\end{array}
\right).
$$
\end{itemize}
\end{defi}

Theorem \ref{propcombinatoria} in \S \ref{Main combinatorial property} establishes a combinatorial property of Rauzy classes which implies that any pair $(\be,\al)$ with $\pi_{0}^{t}(\al)>1$ and $\pi_{0}^{b}(\be)>1$ satisfies either property A or property B in Definition \ref{defproprietaAeB} (or both).

\subsubsection{Steps of the algorithm producing reduced triples}\label{s3ss3sss2}

Consider a fixed Rauzy path $\eta$. For any other path $\ga$ we say that $\ga$ \emph{ends with $\eta$} if there exists $\nu$ such that $\ga=\nu\eta$.

\begin{lem}\label{lem1As3ss3sss2}
Consider a pair $(\be,\al)$ satisfying property A in Definition \ref{defproprietaAeB} and an element $\pi=\pi(\be,\al)$ in $\cR$ such that $\pi^{b}(\be)=\pi^{t}(\al)=d$. Then there exists a finite path $\eta$ in $\Pi(\cR)$ with the following properties:
\begin{itemize}
\item
The last arrow of $\eta$ is a top arrow with loser $\al$.
\item
$\pi$ is in third-to-last position in $\eta$
\item
If $\ga$ is a path with length $r$ and ending with $\eta$, then there exists a positive integer $n=n(r)$ with $n\leq \|q^{\ga}\|$ such that for any $T$ in $\De_{\ga}$ the triple $(\be,\al,n)$ is reduced for $T$ and
$
\la^{(r)}_{\al}=
|T^{n}u^{b}_{\be}-u^{t}_{\al}|
$.
\end{itemize}
\end{lem}

\begin{demo}
Let $\ga^{t}_{\al}$ be the top arrow with winner $\al$ connecting $\pi$ to $R^{t}(\pi)$. Then let $\ga^{b}_{W}$ be the bottom arrow starting from $R^{t}(\pi)$ with loser $\al$, where $W$ is the winner of this second arrow. Consider the concatenation $\ga^{t}_{\al}\ga^{b}_{W}$ of these two arrows and let $\eta$ be any path ending with $\ga^{t}_{\al}\ga^{b}_{W}$. The first two properties of $\eta$ are therefore evident.

Let $\ga$ be a path with length $r$ and ending with $\eta$ and consider $T$ in $\De_{\ga}$. According to \S \ref{S2ss1sss3} in the background let $l(\be,r-2)$ and $h(\al,r-2)$ be the non-negative integers such that for the singularities of $T^{(r-2)}$ we have $u^{(r-2),b}_{\be}=T^{l(\be,r-2)}u^{b}_{\be}$ and $u^{(r-2),t}_{\al}=T^{-h(\al,r-2)}u^{t}_{\al}$. We set $n:=l(\be,r-2)+h(\al,r-2)$. In \S \ref{S2ss1sss3} we prove that $q^{(r)}_{\xi}=l(\xi,r)+h(\xi,r)+1$ for any $\xi$ and any $r$, therefore we have obviously $n<\|q^{(r-2)}\|$ and thus $n<\|q^{(r)}\|$. Observe that $\ga$ ends with $\ga^{t}_{\al}\ga^{b}_{W}$. According to Lemma \ref{lem1s2ss1sss2} we have $T^{(r-2)}\in \De_{\ga^{t}_{\al}\ga^{b}_{W}}$ and in particular $T^{(r-2)}\in \De_{\ga^{t}_{\al}}$. Since $\ga^{t}_{\al}$ is a top arrow with winner $\al$ and loser $\be$, then Equation (\ref{eq3ss1s2(algorithmR-V)}) in the background implies that for the step $T^{(r-1)}=Q(T^{(r-2)})$ of the algorithm we have
$$
\la^{(r-1)}_{\al}=
|u^{(r-2),b}_{\be}-u^{(r-2),t}_{\al}|
$$
(observe that $u^{(r-2),b}_{\be}$ and $u^{(r-2),t}_{\al}$ are the rightmost singularities of $T^{(r-2)}$ respectively in the top and bottom row). Furthermore $T^{(r-1)}\in\De_{\ga^{b}_{W}}$, thus the winner of $T^{(r-1)}$ is $W$ and we have $\la^{(r-1)}_{\al}=\la^{(r)}_{\al}$ because $\al\not=W$. Summing up, we proved that
$$
\la^{(r)}_{\al}=
|T^{l(r-2,\be)}u^{b}_{\be}-T^{-h(r-2,\al)}u^{t}_{\al}|.
$$

Denote $J$ the subinterval of $I^{(r-2)}$ whose endpoints are $u^{(r-2),b}_{\be}$ and $u^{(r-2),t}_{\al}$. The triple $(\be,\al,n)$ is reduced for $T$ if and only if for any $k\in\{-l(\be,r-2),\dots,h(\al,r-2)\}$ the image $T^{k}(J)$ of $J$ does not intersect in its interior any singularity $u^{t}_{\xi}$ of $T$ or any singularity $u^{b}_{\xi}$ of $T^{-1}$. Let $I^{(r-2),t}_{\al}$ and $I^{(r-2),b}_{W}$ be the subintervals of $I^{(r-2)}$ where respectively $T^{(r-2)}$ and $(T^{(r-2)})^{-1}$ act as a translation. The return times to $I^{(r-2)}$ under iteration of $T$ for these intervals are respectively $q^{(r-2)}_{\al}$ and $-q^{(r-2)}_{W}$ (see \S \ref{S2ss1sss3}). Since $J\subset I^{(r-2),t}_{\al}\cap I^{(r-2),b}_{W}$, then the required property on $T^{k}(J)$ holds for $k\in\{-q^{(r-2)}_{W},\dots,q^{(r-2)}_{\al}\}$. It is enough to prove that $h(\al,r-2)\leq q^{(r-2)}_{\al}$ and $l(\be,r-2)\leq q^{(r-2)}_{W}$. Condition $h(\al,r-2)\leq q^{(r-2)}_{\al}$ follows directly from the relation $q^{(r-2)}_{\al}=h(\al,r-2)+l(\al,r-2)+1$. On the other hand observe that $u^{(r-2),b}_{\be}$ is an endpoint of $I^{(r-2),b}_{W}$ and we have $u^{b}_{\be}=T^{-l(\be,r-2)}(u^{(r-2),b}_{\be})$. If $l(\be,r-2)> q^{(r-2)}_{W}$ then $m=l(\be,r-2)-q^{(r-2)}_{W}$ is a positive integer with $T^{m}u^{b}_{\be}\in I^{(r-2)}$ and this is absurd because in \S \ref{S2ss1sss3} we prove that $l(\be,r-2)$ is the smallest integer such that $T^{l(\be,r-2)}u^{b}_{\be}\in I^{(r-2)}$. The lemma is proved.
\end{demo}

\begin{lem}\label{lem1Bs3ss3sss2}
Consider a pair $(\be,\al)$ satisfying property B, an element $\pi=\pi(\be,\al)$ in $\cR$ as in Definition \ref{defproprietaAeB} and the associated letters $V$ and $L$. Then there exists a finite path $\eta$ in $\Pi(\cR)$ with the following properties:
\begin{itemize}
\item
The last arrow of $\eta$ is a bottom arrow with winner $\al$.
\item
$\pi$ in second-to-last position in $\eta$.
\item
If $\ga$ is a path with length $r$ and ending with $\eta$, then we have an integer $n=n(r)$ with $n\leq \|q^{\ga}\|$ such that for $T$ in $\De_{\ga}$ we have
$
\la^{(r)}_{V}=
|T^{n}u^{b}_{\be}-u^{t}_{\al}|
$.
Furthermore if $\la^{(r)}_{V}<\la^{(r)}_{L}$, then the triple $(\be,\al,n)$ is reduced for $T$.
\end{itemize}
\end{lem}

\begin{demo}
Consider the bottom arrow $\ga^{b}_{\al}$ with winner $\al$ starting at $\pi$ and take any path $\eta$ having $\ga^{b}_{\al}$ as last arrow. The first two statements are therefore evident.

Let $\ga$ be a path with length $r$ and ending with $\eta$ and consider $T$ in $\De_{\ga(r)}$. Since $\pi$ is second to last in $\eta$ and $\ga(r)$ ends with $\eta$, then $T^{(r-1)}\in\De_{\pi}$ and the combinatorics of $\pi$ implies
$$
\la^{(r-1)}_{V}=
|u^{(r-1),b}_{\be}-u^{(r-1),t}_{\al}|.
$$
As in Lemma \ref{lem1As3ss3sss2} consider the non-negative integers $l(\be,r-1)$ and $h(\al,r-1)$ such that the singularities of $T^{(r-1)}$ satisfy $u^{(r-1),t}_{\al}=T^{-h(\al,r-1)}u^{t}_{\al}$ and $u^{(r-1),b}_{\be}=T^{l(\be,r-1)}u^{b}_{\be}$. We set $n:=l(r-1,\be)+h(r-1,\al)$ and as in Lemma \ref{lem1As3ss3sss2} we have $n<\|q^{(r-1)}\|<\|q^{(r)}\|$. Since $T^{(r-1)}\in\De_{\ga^{b}_{\al}}$ then $\al$ is the winner of $T^{(r-1)}$ and we have $\la^{(r)}_{V}=\la^{(r-1)}_{V}$, because $V\not=\al$. Summing up we get
$$
\la^{(r)}_{V}=|T^{l(r-1,\be)}u^{b}_{\be}-T^{-h(r-1,\al)}u^{t}_{\al}|.
$$

Assume that $\la^{(r)}_{V}<\la^{(r)}_{L}$ and denote $J$ the subinterval of $I^{(r-1)}$ whose endpoints are $u^{(r-1),b}_{\be}$ and $u^{(r-1),t}_{\al}$. The triple $(\be,\al,n)$ is reduced for $T$ if and only if the iterates $T^{k}(J)$ do not contain any singularity $u^{t}_{\xi}$ of $T$ or any singularity $u^{b}_{\xi}$ of $T^{-1}$ for $k\in\{-l(\be,r-1),\dots,h(\al,r-1)\}$.

Condition $T^{(r-1)}\in\De_{\ga^{b}_{\al}}$ implies $\la^{(r-1)}_{V}<\la^{(r-1)}_{\al}$ and it follows from the combinatorics of $\pi$ that $u^{(r-1),t}_{\al}<u^{(r-1),b}_{\be}<u^{(r-1),t}_{\al}+\la^{(r-1)}_{\al}$, that is $J\subset I^{(r-1),t}_{\al}$. On the other hand we are assuming that $\la^{(r)}_{V}<\la^{(r)}_{L}$, which is equivalent to $\la^{(r-1)}_{V}<\la^{(r-1)}_{L}$ because $L\not=\al$. Therefore the combinatorics of $\pi$ implies $J\subset I^{(r-1),b}_{L}$. The required condition on $T^{k}(J)$ is satisfied for $k\in\{-q^{(r-1)}_{L},\dots,q^{(r-1)}_{V}\}$ and arguing as in Lemma \ref{lem1As3ss3sss2} we get $l(\be,r-1)\leq q^{(r-1)}_{L}$ and $h(\al,r-1)\leq q^{(r-1)}_{V}$. The lemma is proved.
\end{demo}

Theorem \ref{propcombinatoria} implies that any pair $(\be,\al)$ with $\pi^{t}_{0}(\al)>1$ and $\pi^{b}_{0}(\be)>1$ satisfies either property A or property B in Definition \ref{defproprietaAeB}. According to the two cases we can associate to $(\be,\al)$ a path $\eta$ respectively as in Lemma \ref{lem1As3ss3sss2} or as in Lemma \ref{lem1Bs3ss3sss2}. We resume the results in this section with the following proposition.

\begin{prop}\label{props3ss3}
Let $(\be,\al)$ be a pair with $\pi^{t}_{0}(\al)>1$ and $\pi^{b}_{0}(\be)>1$.

If $(\be,\al)$ satisfies property A, let $\eta$ be a path as in Lemma \ref{lem1As3ss3sss2}. For $T$ in $\De_{\pi_{0}}$ without connections let $(r_{k})_{k\in\NN}$ be instants such that $\ga(T,r_{k})$ ends with $\eta$. Then to any such $r_{k}$ it corresponds an integer $n(k)<\|q^{\ga(T,r_{k})}\|$ such that $(\be,\al,n(k))$ is reduced for $T$ and
$$
|T^{n(k)}u^{b}_{\be}-u^{t}_{\al}|=
\la^{(r_{k})}_{\al}.
$$

If $(\be,\al)$ satisfies property B, let $V$ and $L$ be the associated letters as in Definition \ref{defproprietaAeB} and let $\eta$ be a path given by Lemma \ref{lem1Bs3ss3sss2}. For $T$ in $\De_{\pi_{0}}$ without connections let $(r_{k})_{k\in\NN}$ be instants such that $\ga(T,r_{k})$ ends with $\eta$. Then to any such $r_{k}$ it corresponds an integer $n(k)<\|q^{\ga(T,r_{k})}\|$ such that
$$
|T^{n(k)}u^{b}_{\be}-u^{t}_{\al}|=
\la^{(r_{k})}_{V}.
$$
Moreover, if $\la^{(r_{k})}_{L}<\la^{(r_{k})}_{V}$ then $(\be,\al,n(k))$ is reduced for $T$.
\end{prop}

\subsubsection{Proof of Theorem \ref{teoremaf}}\label{s3ss3sss3}

The first part of Theorem \ref{teoremaf} is proved by Lemma \ref{lem1s3ss1sss0}, thus it only remains to prove the second part of the Theorem. Let $(\be,\al)$ be a pair with $\pi^{t}_{0}(\al)>1$ and $\pi^{b}_{0}(\be)>1$. According to Theorem \ref{propcombinatoria} it satisfies either property A or property B in Definition \ref{defproprietaAeB}.

If $(\be,\al)$ satisfies property A, let $\eta$ be a path as in Lemma \ref{lem1As3ss3sss2}. Recall that the Zorich's map $\cZ$ is an ergodic acceleration of the normalized map $\widehat{Q}$, therefore the $\widehat{Q}$-orbit of almost any $\widehat{T}$ in $\De^{(1)}_{\pi_{0}}$ enters in $\De^{(1)}_{\eta}$ infinitely many times. According to Lemma \ref{lem1s2ss1sss2} in the background this is equivalent to have infinitely many instants $(r_{k})_{k\in\NN}$ such that $\ga(\widehat{T},r_{k})$ ends with $\eta$. Observe that for $T=(\pi_{0},\la)$ and $\widehat{T}=(\pi_{0},\widehat{\la})$ with $\widehat{\la}=\|\la\|^{-1}\la$ we have $\ga(T,r_{k})=\ga(\widehat{T},r_{k})$. Therefore for almost any $T$ in $\De_{\pi_{0}}$ there exist infinitely many instants $(r_{k})_{k\in\NN}$ such that $\ga(T,r_{k})$ ends with $\eta$. The first part of Proposition \ref{props3ss3} implies the statement for the pair $(\be,\al)$.

On the other hand, if $(\be,\al)$ satisfies property B, let $V$ and $L$ be the associated letters as in Definition \ref{defproprietaAeB} and let $\eta$ be a path as in Lemma \ref{lem1Bs3ss3sss2}. Let $\pi_{1}$ be ending point of $\eta$ and $l$ be the length of $\eta$. For $\widehat{T}\in\De^{(1)}_{\eta}$ write $\widehat{T}^{(l)}=(\pi_{1},\widehat{\la}^{(l)})$. The map $\widehat{T}\mapsto\widehat{T}^{(l)}$ is an homeomorphism between $\De^{(1)}_{\eta}$ and $\De^{(1)}_{\pi_{1}}$, thus
$
\{\widehat{T}\in\De^{(1)}_{\eta};
\widehat{\la}^{(l)}_{V}<\widehat{\la}^{(l)}_{L}\}
$
is an open subset of $\De^{(1)}_{\eta}$ with positive measure. Therefore the $\widehat{Q}$-orbit of almost any $\widehat{T}$ in $\De^{(1)}_{\pi_{0}}$ enters in
$
\{\widehat{T}\in\De^{(1)}_{\eta};
\widehat{\la}^{(l)}_{V}<\widehat{\la}^{(l)}_{L}\}
$
infinitely many times. Arguing as in the previous case we get that for almost any $T$ in $\De_{\pi_{0}}$ there exist infinitely many instants $(r_{k})_{k\in\NN}$ such that $\ga(T,r_{k})$ ends with $\eta$ and $\la^{(r_{k})}_{V}<\la^{(r_{k})}_{L}$. The second part of Proposition \ref{props3ss3} implies the statement for the pair $(\be,\al)$. Theorem \ref{teoremaf} is proved.

\subsubsection{Proof of Proposition \ref{proposizioneg}}\label{s3ss3sss4}

Fix a pair $(\be,\al)$ with $\pi^{b}_{0}(\be)>1$ and $\pi^{t}_{0}(\al)>1$. Recall from \S \ref{s3ss1sss1} that for any $\ga$ in $\Ga(\be,\al,n)$ we denote $\De^{\ast}_{\ga}$ the subset of those $T$ in $\De_{\ga}$ such that the triple $(\be,\al,n)$ is reduced for $T$. Denote $\De^{\ast\ast}_{\ga}$ the subset of those $T$ in $\De^{\ast}_{\ga}$ such that $|T^{n}u^{b}_{\be}-u^{t}_{\al}|<\varphi(n)$. The union
$$
\De^{\ast\ast}(\Ga(\be,\al,n)):=
\bigsqcup_{\ga\in\Ga(\be,\al,n)}
\De^{\ast\ast}_{\ga}
$$
is the set of all $T$ such that $(\be,\al,n)$ is a reduced triple for $T$ satisfying (\ref{eqKhinTSI}).

\begin{lem}\label{lem1s3ss3sss4}
If $\{\varphi(n)\}_{n\in\NN}$ be a positive sequence, then for any $N\in\NN$ the union $\bigcup_{n\geq N}\De^{\ast\ast}(\Ga(\be,\al,n))$ is an open and dense subset of $\De_{\pi_{0}}$.
\end{lem}

\begin{demo}
Since $\varphi(n)>0$, then for any $\ga$ in $\Ga(\be,\al,n)$ the set of those $T=(\pi_{0},\la)$ in $\De_{\ga}$ with $|\langle\la,B_{\ga}^{-1}w_{\be,\al}(\pi)\rangle|<\varphi(n)$ is open in $\De_{\ga}$. According to Proposition \ref{prop1s3ss1sss1} this last condition is equivalent to $|T^{n}u^{b}_{\be}-u^{t}_{\al}|<\varphi(n)$, moreover the set $\De^{\ast}_{\ga}$ is open. Therefore $\De^{\ast\ast}_{\ga}$ is an open subset of $\De_{\ga}$ for any $\ga$ in $\Ga(\be,\al,n)$. It just remains to prove density of $\bigcup_{n\geq N}\De^{\ast\ast}(\Ga(\be,\al,n))$.

Fix any open set $X$ in $\De_{\pi_{0}}$. We prove that there exists $n\geq N$ and $\ga$ in $\Ga(\be,\al,n)$ such that $X\cap\De^{\ast\ast}_{\ga}\not=\emptyset$. Since $X$ has positive Lebesgue measure, according to the second part of Theorem \ref{teoremaf} there exists $T$ in $X$, a letter $\xi$ and infinitely many integers $r_{k}$ and $n(k)$ such that $(\be,\al,n(k))$ is a reduced triple for $T$ with
$$
\la^{(r_{k})}_{\xi}=
|T^{n(k)}u^{b}_{\be}-u^{t}_{\al}|,
$$
in particular $\xi=\al$ if the pair $(\be,\al)$ satisfies property A and $\xi=V$ if $(\be,\al)$ satisfies property B. Moreover we can also suppose that $T$ is uniquely ergodic.
It is well-known that $\la^{(r)}_{\xi}\to 0$ as $r\to\infty$ for any $\xi\in\cA$ (see Corollary 1 at page 37 of \cite{ytre}). Therefore there exists $k$ with $n(k)\geq N$, otherwise we have $\la^{n(k)}_{\xi}>\epsilon:=\min_{n<N}|T^{n}u^{b}_{\be}-u^{t}_{\al}|$ for all $k\in\NN$, which is absurd. Since $T$ is uniquely ergodic, then the intersection of the infinitely many cones
$$
\De_{\ga(T,r_{1})}\supset\dots\supset
\De_{\ga(T,r_{k})}\supset
\De_{\ga(T,r_{k+1})}\cdots
$$
is the half-line in $\De_{\pi_{0}}$ spanned by $T$ (see \S 8.1 in \cite{ytre}). Therefore for $k$ big enough we have $\proiettivo(\De_{\ga(T,r_{k})})\subset\proiettivo(X)$, where $\proiettivo(X)$ denotes the space of rays spanned by the elements of $X$. Fix $k$ with $n(k)\geq N$ and $\proiettivo(\De_{\ga(T,r_{k})})\subset\proiettivo(X)$, set $n:=n(k)$ and let $\ga$ be the unique path in $\Ga(\be,\al,n)$ such that $\ga\prec\ga(T,r_{k})$. Since $\varphi(n)>0$ there exists $T_{\ast}=(\pi_{0},\la_{\ast})$ in $\De_{\ga(T,r_{k})}\cap X$ such that either $\la^{(r_{k})}_{\ast\al}<\varphi(n)$ if $(\be,\al)$ has property A, or $\la^{(r_{k})}_{\ast V}<\min\{\la^{(r_{k})}_{\ast L},\varphi(n)\}$ if $(\be,\al)$ has property B. Proposition \ref{props3ss3} implies that the triple $(\be,\al,n)$ is reduced for such $T_{\ast}$ and satisfies (\ref{eqKhinTSI}). Since $\De_{\ga(T,r_{k})}\subset\De_{\ga}$, we have $T_{\ast}\in\De^{\ast\ast}_{\ga}$, thus $X\cap\De^{\ast\ast}_{\ga}\not=\emptyset$. The lemma is proved.
\end{demo}

Lemma \ref{lem1s3ss3sss4} implies that for any pair $(\be,\al)$ with $\pi^{b}_{0}(\be)>1$ and $\pi^{t}_{0}(\al)>1$ the intersection
$
\bigcap_{N=1}^{+\infty}
\bigcup_{n\geq N}
\De^{\ast\ast}(\Ga(\be,\al,n))
$
is a residual set. The elements of the set
$$
\bigcap_{\pi^{b}_{0}(\be)>1,\pi^{t}_{0}(\al)>1}
\bigg(
\bigcap_{N=1}^{+\infty}
\bigcup_{n\geq N}
\De^{\ast\ast}(\Ga(\be,\al,n))
\bigg)
$$
are $\varphi$-Liouville i.e.t.s and since the intersection of finitely many residual sets is still a residual set, then Proposition \ref{proposizioneg} follows.

\section{Dichotomy for Khinchin-type condition}\label{s4}

This section is devoted to the proof of Theorem \ref{teoremaa}. We first prove that Proposition \ref{propcasoconvergente} and Proposition \ref{propcasodivergente} imply the theorem.

\begin{lem}\label{lem1s4ss0sss0}
Proposition \ref{propcasoconvergente} implies the convergent case of Theorem \ref{teoremaa}. On the other hand Proposition \ref{propcasodivergente} implies the divergent case of Theorem \ref{teoremaa}.
\end{lem}

\begin{demo}
Let $\pi_{0}$ and $(\be,\al)$ be as in Theorem \ref{teoremaa}. For any $\rho>0$ define the set $\De_{\pi_{0},\rho,+}$ of those $T=(\pi_{0},\la)$ in $\De_{\pi_{0}}$ with $\|\la\|>\rho$. Similarly define the set $\De_{\pi_{0},\rho,-}$ of those $T=(\pi_{0},\la)$ in $\De_{\pi_{0}}$ with $\|\la\|<\rho$. Finally, for any positive sequence $\varphi=\{\varphi(n)\}_{n\in\NN}$ let $\varphi_{\rho}=\{\varphi_{\rho}(n)\}_{n\in\NN}$ be the sequence given by $\varphi_{\rho}(n):=\rho\varphi(n)$. For any $T=(\pi_{0},\la)$ in $\De_{\pi_{0}}$ let $\widehat{T}=(\pi_{0},\widehat{\la})$ be the corresponding normalized i.e.t. in $\De^{(1)}_{\pi_{0}}$, where $\widehat{\la}:=\|\la\|^{-1}\la$. If $u^{b}_{\be}$ and $u^{t}_{\al}$ are the singularities for $T$, denote $\widehat{u}^{b}_{\be}$ and $\widehat{u}^{t}_{\al}$ the corresponding singularities for $\widehat{T}$. Observe that a triple $(\be,\al,n)$ is reduced for $T$ if and only if it is reduced for $\widehat{T}$, moreover we have
$$
|T^{n}u^{b}_{\be}-u^{t}_{\al}|=
\|\la\|\cdot
|\widehat{T}^{n}\widehat{u}^{b}_{\be}-\widehat{u}^{t}_{\al}|.
$$

We first prove that Proposition \ref{propcasoconvergente} implies the convergent case of Theorem \ref{teoremaa}. Consider $\varphi$ such that $\varphi(n)$ is decreasing monotone and $\sum_{n=1}^{+\infty}\varphi(n)<+\infty$. Fix $\rho>0$. The sequence $\varphi_{\rho^{-1}}$ satisfies the assumption in Proposition \ref{propcasoconvergente}. For any $T$ in $\De_{\pi_{0},\rho,+}$ and any $n$ we have
$$
|\widehat{T}^{n}\widehat{u}^{b}_{\be}-\widehat{u}^{t}_{\al}|<
\rho^{-1}|T^{n}u^{b}_{\be}-u^{t}_{\al}|.
$$
Consider $T$ in $\De_{\pi_{0},\rho,+}$ such that there exist infinitely many triples $\{(\be,\al,n_{k})\}_{k\in\NN}$ satisfying (\ref{eqKhinTSI}) with respect to the sequence $\varphi$. Any triple $(\be,\al,n_{k})$ satisfies (\ref{eqKhinTSI}) for the corresponding normalized $\widehat{T}$ with respect to the sequence $\varphi_{\rho^{-1}}$. Suppose that $\cS$ is a positive-measure subset of $\De_{\pi_{0},\rho,+}$ such that for any $T$ in $\cS$ there exist infinitely many triples $(\be,\al,n_{k})$ as above. Then $\widehat{\cS}:=\{\widehat{T};T\in\cS\}$ is a subset of  $\De^{(1)}_{\pi_{0}}$ with positive measure, which is absurd by Proposition \ref{propcasoconvergente}. Summing up it follows that for almost any $T$ in $\De_{\pi_{0},\rho,+}$ there exist only finitely many triples $(\be,\al,n)$ reduced for $T$ and such that $|T^{n}u^{b}_{\be}-u^{t}_{\al}|<\varphi(n)$. We observe that we can repeat the argument for arbitrary small $\rho$ and for any pair of letters $(\be,\al)$ with $\pi^{b}_{0}(\be)>1$ and $\pi^{t}_{0}(\al)>1$, thus the convergent part of Theorem \ref{teoremaa} follows.

Now we prove that Proposition \ref{propcasodivergente} implies the divergent case of Theorem \ref{teoremaa}. Consider $\varphi$ such that $n\varphi(n)$ is decreasing monotone and $\sum_{n=1}^{+\infty}\varphi(n)=+\infty$. Fix $\rho>0$. For any $T$ in $\De_{\pi_{0},\rho,-}$ and any $n$ we have
$$
|T^{n}u^{b}_{\be}-u^{t}_{\al}|<
\rho
|\widehat{T}^{n}\widehat{u}^{b}_{\be}-\widehat{u}^{t}_{\al}|.
$$
The sequence $\varphi_{\rho^{-1}}$ satisfies the assumption in  Proposition \ref{propcasodivergente}, therefore for almost any $T_{\ast}$ in $\De^{(1)}_{\pi_{0}}$ there exist infinitely many triples $(\be,\al,n)$ reduced for $T_{\ast}$ and satisfying (\ref{eqKhinTSI}) with respect to the sequence $\varphi_{\rho^{-1}}$. Fix such a $T_{\ast}$ and let $\{(\be,\al,n_{k})\}_{k\in\NN}$ be the infinite family of triples reduced for $T_{\ast}$ and satisfying (\ref{eqKhinTSI}) with respect to the sequence $\varphi_{\rho^{-1}}$. For any $T$ in $\De_{\pi_{0},\rho,-}$ with $\widehat{T}=T_{\ast}$ any triple $(\be,\al,n_{k})$ is also reduced for $T$ and satisfies condition (\ref{eqKhinTSI}) with respect to the sequence $\varphi$. Therefore for almost any $T$ in $\De_{\pi_{0},\rho,-}$ there exists infinitely many triples $(\be,\al,n)$ reduced for $T$ and such that $|T^{n}u^{b}_{\be}-u^{t}_{\al}|<\varphi(n)$. Finally we observe that we can repeat the argument for arbitrary big $\rho$, thus the divergent part of Theorem \ref{teoremaa} follows.
\end{demo}

\subsection{Proof of convergent case}\label{s4ss1}

According to Theorem \ref{teoremaf} reduced triples are detected by the Rauzy-Veech algorithm. Here we show that Proposition \ref{propcasoconvergente} follows if we prove it for reduced triples.

\begin{lem}\label{lem1s4ss1sss0}
Consider $T$ in $\De_{\pi_{0}}$ (or in $\De^{(1)}_{\pi_{0}}$) such that there exist a sequence of triples $(\be_{k},\al_{k},m_{k})$ with $k\in\NN$ satisfying condition (\ref{eqKhinTSI}) but not necessarily reduced for $T$. Then there exists a pair $(\be,\al)$ and infinitely many integers $n_{k}$ such that $(\be,\al,n_{k})$ is a reduced triple for $T$ which satisfies (\ref{eqKhinTSI}).
\end{lem}

\begin{demo}
Let $(\be,\al,m)$ be a triple such that $|T^{m}u^{b}_{\be}-u^{t}_{\al}|<\varphi(m)$ but which is not reduced for $T$. Then there exists $k$ in $\{0,\dots,m\}$ and a letter $\xi$ in $\cA$ such that $T^{-k}I(\be,\al,m)$ contains in its interior either $u^{t}_{\xi}$ or $u^{b}_{\xi}$, moreover we can suppose that $k$ is minimal with one of the two properties. If $u^{t}_{\xi}\in T^{-k}I(\be,\al,m)$ then $|T^{m-k}u^{b}_{\be}-u^{t}_{\xi}|<|T^{m}u^{b}_{\be}-u^{t}_{\al}|$, thus $|T^{m-k}u^{b}_{\be}-u^{t}_{\xi}|<\varphi(m-k)$, since $\varphi$ is monotone. Similarly, if $u^{b}_{\xi}\in T^{-k}I(\be,\al,m)$ then $|u^{b}_{\xi}-T^{-k}u^{t}_{\al}|<|T^{m}u^{b}_{\be}-u^{t}_{\al}|$, so $|T^{k}u^{b}_{\xi}-u^{t}_{\al}|<\varphi(k)$ by minimality of $k$. In both cases we pass from $(\be,\al,m)$ to an other triple $(\be',\al',n)$ satisfying equation (\ref{eqKhinTSI}) with $n<m$. Applying iteratively the argument we get a triple reduced for $T$ which still satisfies equation (\ref{eqKhinTSI}).

Now let us suppose that $T$ admits infinitely many triples $(\be_{k},\al_{k},m_{k})_{k\in\NN}$ which are solutions of equation (\ref{eqKhinTSI}), but not necessarily reduced. With the argument above we get a sequence of reduced solutions $\{(\be'_{k},\al'_{k},n_{k})\}_{k\in\NN}$ for $T$. Finally there exist at least one pair $(\be,\al)$ appearing infinitely many times in the sequence. The lemma is proved.
\end{demo}

\subsubsection{Proof of Proposition \ref{propcasoconvergente}}\label{s4ss1sss1}

Fix a combinatorial datum $\pi_{0}$ and a pair of letters $(\be,\al)$ as in Theorem \ref{teoremaa}. For any $n\in\NN$ denote $\cI(\be,\al,n)$ the set of those $T$ in $\De^{(1)}_{\pi_{0}}$ such that the triple $(\be,\al,n)$ is reduced for $T$ and satisfies condition (\ref{eqKhinTSI}).

According to Proposition \ref{prop1s3ss1sss1} the set of those $T$ in $\De^{(1)}_{\pi_{0}}$ such that the triple $(\be,\al,n)$ is reduced for $T$ is contained in $\bigsqcup_{\ga\in\Ga(\be,\al,n)}\De^{(1)}_{\ga}$, therefore we have
$$
\cI(\be,\al,n)=
\bigsqcup_{\ga\in\Ga(\be,\al,n)}
\cI(\be,\al,n)\cap\De^{(1)}_{\ga},
$$
where the union is disjoint since $\Ga(\be,\al,n)$ is a disjoint family. Moreover Lemma \ref{lem1s3ss2sss2} implies that for any $\ga\in\Ga(\be,\al,n)$, denoting $W(\ga)$ the last winner of $\ga$, we have
$$
\leb_{d-1}(\cI(\be,\al,n)\cap\De^{(1)}_{\ga})<
\varphi(n)q^{\ga}_{W(\ga)}\leb_{d-1}(\De^{(1)}_{\ga}).
$$
Recalling that we denote $\PP(\ga)=\leb_{d-1}(\De^{(1)}_{\ga})$ we have
$$
\sum_{n\in\NN}\leb(\cI(\be,\al,n))
\leq
\sum_{n\in\NN}\varphi(n)
\sum_{\ga\in\Ga(\be,\al,n)}
q^{\ga}_{W(\ga)}\PP(\ga)\leq
$$
$$
\sum_{N\in\NN}
\varphi(2^{N-1})
\sum_{2^{N-1}\leq n<2^{N}}
\bigg(
\sum_{\ga\in\Ga(\be,\al,n)}
q^{\ga}_{W(\ga)}\PP(\ga)
\bigg)\leq
2C\sum_{N\in\NN}
2^{N-1}\varphi(2^{N-1}),
$$
where the second inequality follows since $\varphi(n)$ is decreasing monotone and the third is consequence of Proposition \ref{props3ss2}, where $C$ is the constant in the proposition. We recall the following classic result in calculus

\begin{lem}\label{lem1s4ss1sss1}
Let $\varphi:[1,+\infty)\to(0,+\infty)$ be a positive and monotone decreasing function and consider any real number $\te>1$. Then $\sum_{n=1}^{+\infty}\varphi(n)<+\infty$ if and only if $\sum_{N=1}^{+\infty}\te^{N}\varphi(\te^{N})$.
\end{lem}

\begin{demo}
Since $\varphi$ is decreasing monotone, for $\te^{N-1}\leq n<\te^{N}$ we have
$$
(\te^{N}-\te^{N-1})\varphi(\te^{N})\leq
\sum_{\te^{N-1}\leq n<\te^{N}}\varphi(n)<
(\te^{N}-\te^{N-1})\varphi(\te^{N-1}),
$$
hence
$
(1-\te^{-1})\sum_{N=1}^{+\infty}
\te^{N}\varphi(\te^{N})\leq
\sum_{n=1}^{+\infty}\varphi(n)<
(\te-1)\sum_{N=1}^{+\infty}
\te^{N-1}\varphi(\te^{N-1}).
$
The lemma is proved.
\end{demo}

Since we are assuming that $\varphi(n)$ is decreasing monotone with $\sum_{n=1}^{+\infty}\varphi(n)<+\infty$, then Lemma \ref{lem1s4ss1sss1} implies $\sum_{N=1}^{+\infty}2^{N-1}\varphi(2^{N-1})<+\infty$, that is
$$
\sum_{n=1}^{+\infty}
\leb_{d-1}(\cI(\be,\al,n))
<+\infty.
$$
According to the first part of Borel-Cantelli Lemma almost any $T\in\De^{(1)}_{\pi_{0}}$ is contained just in finitely many sets $\cI(\be,\al,n)$, that is there are just finitely many integers $n$ to which correspond a triple $(\be,\al,n)$ reduced for $T$ and solution of (\ref{eqKhinTSI}). This last condition holds for any pair $(\be,\al)$, then Lemma \ref{lem1s4ss1sss0} implies that Proposition \ref{propcasoconvergente} holds. The convergent case of Theorem \ref{teoremaa} is proved.

\subsection{Shrinking target criterion for the divergent case}\label{s4ss2}

We consider the normalized version $\widehat{Q}$ of the Rauzy-Veech algorithm. For $T$ in $\De^{(1)}_{\pi_{0}}$ without connections we call $(\pi^{(r)},\widehat{\la}^{(r)})$ the pair of combinatorial and length data of $\widehat{T}^{(r)}=\widehat{Q}^{r}(T)$. Recall that $\ga(T,\infty)$ denotes the half-infinite path in the Rauzy diagram generated by $T$ and $\ga(T,r)$ is the concatenation of the first $r$ arrows of $\ga(T,\infty)$.

\subsubsection{First return to a neat path}\label{s4ss2sss1}

Let $\eta$ be any finite path of length $l$, that is $\eta$ is the concatenation $\eta_{1}\dots\eta_{l}$ of $l$ elementary arrows. Consider $T$ in $\De^{(1)}_{\pi_{0}}$ without connections. For $r\in\NN$ the path $\ga(T,r)$ ends with $\eta$ if and only if $\widehat{T}^{(r-l)}\in\De^{(1)}_{\eta}$, according to Lemma \ref{lem1s2ss1sss2} in the background. Therefore, motivated by Proposition \ref{props3ss3}, we look for instants $r$ such that the iterates $\widehat{T}^{(r)}=\widehat{Q}^{r}(T)$ belong to $\De^{(1)}_{\eta}$. Since $\widehat{Q}$ has an ergodic acceleration (the Zorich's map) there are infinitely many such $r$.\\

We say that a finite Rauzy path $\eta$ is \emph{neat} if any time that we can write $\eta=\eta_{1}\eta_{2}=\eta_{3}\eta_{1}$ either $\eta=\eta_{1}$ or $\eta_{1}$ is trivial.

\begin{lem}\label{lem1s4ss2sss1}
Let us consider any finite Rauzy path $\eta$ and the associated simplicial cone $\De_{\eta}$. Let $l=l(\eta)$ be the number of elementary arrows which compose $\eta$. Then $\eta$ is neat if and only if for any $T\in\De_{\eta}$ we have $T^{(i)}\not\in\De_{\eta}$ for all $i\in\{1,\dots,l-1\}$.
\end{lem}

\begin{demo}
Let us first suppose that $\eta$ is not neat, that is there exist three non trivial paths $\eta_{1},\eta_{2},\eta_{3}$ such that $\eta=\eta_{1}\eta_{2}=\eta_{3}\eta_{1}$. We consider the sub-cone $\De_{\eta\eta_{2}}$ of $\De_{\eta}$ and any $T\in\De_{\eta\eta_{2}}$. Let $i$ be the length of $\eta_{3}$. Since $\eta_{1}$ is not trivial then $1\leq i\leq l-1$ and we have $T^{(i)}\in \De_{\eta_{1}\eta_{2}}=\De_{\eta}$.

On the other hand we have $T\in\De_{\eta}$ if and only if the first $l$ steps of the algorithm $Q$ applied to $T$ are given by the $l$ arrows $\eta_{1},\dots,\eta_{l}$ composing $\eta$. Let us suppose that for some $1\leq i\leq l-1$ we have $T^{(i)}\in\De_{\eta}$. This means that $\eta$ begins with $\eta_{i+1}\dots\eta_{l}$, on the other hand $\eta_{i+1}\dots\eta_{l}$ is also the ending part of $\eta$, so $\eta$ is not neat. The lemma is proved.
\end{demo}

Let $\eta:\pi_{0}\to\pi_{1}$ be a neat path starting at $\pi_{0}$ and ending in $\pi_{1}$. Consider the sub-simplex $\De^{(1)}_{\eta}$ of $\De^{(1)}_{\pi_{0}}$ and the \emph{first entering map} of $\widehat{Q}$ in $\De^{(1)}_{\eta}$, that is the map $R_{\eta}:\Delta^{(1)}(\cR)\to\Delta^{(1)}_{\eta}$ defined by $R_{\eta}(T):=\widehat{Q}^{E(T)}(T)$, where
$$
E(T):=\min \{k\in \NN^{\ast};\widehat{Q}^{k}(T)\in \De^{(1)}_{\eta}\}.
$$
As we argued in the beginning of \S \ref{s4ss2sss1}, $R_{\eta}$ is defined almost everywhere on $\Delta^{(1)}(\cR)$. We denote $R_{\eta,\pi_{1}}:=R_{\eta}|_{\De^{(1)}_{\pi_{1}}}$ its restriction to the simplex $\De^{(1)}_{\pi_{1}}$. Consider the homeomorphism $\widehat{Q}_{\eta}:\Delta^{(1)}_{\eta}\to\Delta^{(1)}_{\pi_{1}}$ defined by
$$
\widehat{Q}_{\eta}(\pi_{0},\la)=
(\pi_{1},\frac{^{t}B^{-1}_{\eta}\la}{\|^{t}B^{-1}_{\eta}\la\|}).
$$
$\widehat{Q}_{\eta}$ can be concatenated with $R_{\eta,\pi_{1}}$. We define the map
$\cF_{\eta}:\De^{(1)}_{\pi_{1}}\to \De^{(1)}_{\pi_{1}}$ by
\begin{equation}\label{eqdefmappaeffe}
\cF_{\eta}(T):=\widehat{Q}_{\eta}\circ R_{\eta,\pi_{1}}(T).
\end{equation}
Let $\Ga^{\eta}$ be the set of Rauzy paths $\ga$ starting and ending in $\pi_{1}$ which contain $\eta$ and are minimal with this property with respect to the ordering $\prec$. In other words the elements of $\Ga^{\eta}$ are the paths $\ga$ which admit a decomposition
$$
\ga=\nu\eta
$$
with $\nu\in\Pi(\cR)$ and have the property that for no proper sub-path $\ga'$ with $\ga'\prec\ga$ the same decomposition is possible. The simplices $\De^{(1)}_{\ga}$ with $\ga\in\Ga^{\eta}$ are exactly the connected components of the domain of the map $\cF_{\eta}$, on each of them $\cF_{\eta}$ acts as a projective linear map, that is, if we write $T=(\pi_{1},\la)$ then
$$
T\in \De^{(1)}_{\ga}
\Leftrightarrow
\cF_{\eta}(T)=
(\pi_{1},\frac{^{t}B_{\ga}^{-1}\la}{\|^{t}B_{\ga}^{-1}\la\|}).
$$
For $k\in\NN$ let us introduce the set $\Ga^{(k),\eta}$ of those finite paths $\ga_{k}$ starting and ending in $\pi_{1}$ which contain exactly $k$ distinct copies of the path $\eta$ and are minimal with this property with respect to the ordering $\prec$. Observe that by minimality all these paths end with $\eta$. The connected components of the $k$-th iterated $\cF_{\eta}^{k}$ of $\cF_{\eta}$ are exactly the simplices $\Delta^{(1)}_{\ga_{k}}$ with $\ga_{k}\in \Ga^{(k),\eta}$. For any $\ga_{k}$ in $\Ga^{(k),\eta}$ and for any $T$ in $\De^{(1)}_{\ga_{k}}$ we have $\cF_{\eta}^{k}(T)=(\pi_{1},\frac{^{t}B_{\ga_{k}}^{-1}\la}{\|^{t}B_{\ga_{k}}^{-1}\la\|})$. Since $\cF_{\eta}$ is defined almost everywhere, for all $k\in \NN$ we have
$$
\De^{(1)}_{\pi_{1}}=
\bigsqcup_{\ga_{k}\in \Ga^{(k),\eta}}
\De^{(1)}_{\ga_{k}}\mod 0.
$$

\subsubsection{Uniform control of the speed of shrinking}\label{s4ss2sss2}

Let $\eta$ be a neath path starting at $\pi_{0}$ and ending at $\pi_{1}$. For a generic $T$ in $\De^{(1)}_{\pi_{0}}$ consider the sequence of instants $r_{k}=r_{k}(T)$ such that the initial segment $\ga(T,r_{k})$ of $\ga(T,\infty)$ ends with $\eta$ for any $k\in\NN$. As we argued in the beginning of \S \ref{s4ss2sss1} for almost any $T$ there are infinitely many such $r_{k}$. In particular $\widehat{T}^{(r_{k})}$ belongs to $\De^{(1)}_{\pi_{1}}$ for any $k$, thus we can write $\widehat{T}^{(r_{k})}=(\pi_{1},\widehat{\la}^{(r_{k})})$. In order to establish a shrinking target criterion (Proposition \ref{props4ss2}) we need to compare the length of an interval of $\widehat{T}^{(r_{k})}$ with the quantity
$$
\|\la^{(r_{k})}\|^{-1}
\varphi(\|q^{\ga(T,r_{k})}\|),
$$
where $\la^{(r_{k})}$ is the length datum of the non-normalized i.e.t. $T^{(r_{k})}=Q^{r_{k}}(T)$. The expression above does not depend just on $k$, but also on $T$ itself. In this section we get a lower bound for it uniform in $T$. We use well known results on the ergodic theory of the Zorich's map $\cZ$ introduced in Definition \ref{defaccelerazionezorich}.

\begin{lem}\label{lem1s4ss2sss2}
For any neat path $\eta$ there exists a constant $\theta >1$ (depending on $\eta$) such that for almost any $T\in \De^{(1)}_{\pi_{0}}$ and for any $k$ big enough we have
$$
\|q^{\ga(T,r_{k})}\|\leq \te^{k}.
$$
\end{lem}

\begin{demo}
We first recall a basic property of the map $\cZ$. Let us denote $\widetilde{\ga}$ paths in the Rauzy diagram corresponding to iterations of $\cZ$. If $\widetilde{\ga}$ corresponds to $N$ iterations of $\cZ$ then it has a decomposition $\widetilde{\ga}=\widetilde{\ga}^{1}\dots\widetilde{\ga}^{N}$, where any $\widetilde{\ga}^{i}$ is concatenation of simple arrows with the same winner and it is maximal with this property. The product of matrices $\widetilde{B}_{N}:=B_{\widetilde{\ga}^{N}}\cdots B_{\widetilde{\ga}^{1}}$ defines a cocycle $T\mapsto (\cZ^{N}(T),\widetilde{B}_{N})$ over $\cZ$, known as \emph{Zorich's cocycle}. Oseledet's Theorem, together with ergodicity of $\cZ$ (Theorem 2 in \cite{zorich}), implies that there exists a constant $\nu>0$ such that for a generic $T$ we have
$$
\lim_{N\to\infty}
\frac{1}{N}
\log(\|B_{\widetilde{\ga}^{N}}\cdots B_{\widetilde{\ga}^{1}}\|)
=\nu,
$$
where $\nu$ is the maximal \emph{Lyapunov exponent} of the cocycle.

For any path $\ga_{k}\in\Ga^{(k),\eta}$ and any $T\in\De^{(1)}_{\ga_{k}}$ without connections let us denote $\widetilde{\ga}_{k}(T)$ the shortest segment of the Zorich's path generated by $T$ which satisfies $\ga_{k}\prec\widetilde{\ga}_{k}(T)$. We obviously have
$$
\|q^{\ga_{k}}\|\leq \|q^{\widetilde{\ga}_{k}(T)}\|.
$$

Fix any (small) $\epsilon>0$. For a generic $T$ and for any $k\in\NN$ consider the decomposition of $\widetilde{\ga}_{k}(T)$ in elementary Zorich's steps: $\widetilde{\ga}_{k}(T)=\widetilde{\ga}^{1}\dots\widetilde{\ga}^{N}$, where $N=N(T,k)$. Applying the result recalled above on Zorich's cocycle, for any $k$ big enough we get
$$
\|q^{\widetilde{\ga}_{k}(T)}\|=
\|B_{\widetilde{\ga}_{1}\dots\widetilde{\ga}_{N}}\|=
\|B_{\widetilde{\ga}^{N}}\cdots B_{\widetilde{\ga}^{1}}\|\leq
e^{N(\nu+\epsilon)}
$$

To complete the proof it is enough to show that there exists a constant $C>0$ such that $N(T,k)\leq Ck$ for almost any $T$ and for any $k$ big enough (then the constant in the statement is $\te=e^{C(\nu+\epsilon)}$). We denote $\chi$ the characteristic function of the set $\De^{(1)}_{\eta}$. Then we denote $S_{n}\chi$ and $\widetilde{S}_{n}\chi$ the $n$-th \emph{Birkhoff sum} of the function $\chi$ respectively under iterations of the map $\widehat{Q}$ and of the map $\cZ$. In formulae
$S_{n}\chi(T):=
\chi(\widehat{Q}^{n-1}(T))
+\dots+
\chi(T)
$
and
$
\widetilde{S}_{n}\chi(T):=
\chi(\cZ^{n-1}(T))
+\dots+
\chi(T).
$
Since $\cZ$ is ergodic then $(1/n)\widetilde{S}_{n}\chi(T)\to \mu(\De^{(1)}_{\eta})$ as $n\to\infty$ for almost any $T$ (where $\mu$ is the smooth ergodic measure for $\cZ$). We set $C:=1/\mu(\De^{(1)}_{\eta})$ and we have $n\leq (C+\epsilon)\widetilde{S}_{n}\chi(T)$ for almost any $T$ and any $n$ big enough.

Now we recall that $\eta$ is neat, thus for any $\ga_{k}\in\Ga^{(k),\eta}$ and any $T\in\De^{(1)}_{\eta}$, recalling the definition of the entering times $r_{k}=r_{k}(T)$ in the statement, we have
$$
S_{r_{k}}\chi(T)=k.
$$
For the same $\ga_{k}$ and $T$ we consider the Zorich path $\widetilde{\ga}_{k}(T)$ introduced above and its decomposition in elementary Zorich's steps $\widetilde{\ga}_{k}(T)=\widetilde{\ga}^{1}\dots\widetilde{\ga}^{N}$, where $N=N(T,k)$. Since $\cZ$ is an acceleration of $\widehat{Q}$ then $\widetilde{S}_{N}\chi(T)\leq S_{r_{k}}\chi(T)=k$, thus, recalling the estimate above $N\leq (C+\epsilon)\widetilde{S}_{N}\chi(T)\leq (C+\epsilon)k$. The estimate for $N$ is complete and therefore the proof of the lemma too.
\end{demo}

We say that a finite Rauzy path $\eta$ is \emph{positive} if the matrix $B_{\eta}$ (introduced in \S \ref{sss1ss1s2}) has all its entries positive.

\begin{lem}\label{lem2s4ss2sss2}
Let $\eta$ be a positive path and set $M:=\|B_{\eta}\|$. Then for any path $\ga$ ending with $\eta$, that is admitting a decomposition $\ga=\nu\eta$, and for any $\al,\be\in\cA$ we have
$$
q^{\ga}_{\al}\leq
Mq^{\ga}_{\be}.
$$
\end{lem}

\begin{demo}
Just observe that $q^{\ga}=B_{\ga}\vec{1}=B_{\eta}B_{\nu}\vec{1}$.
\end{demo}

Define a positive function $\varphi:[1,+\infty)\to(0,+\infty)$ setting $\varphi(t):=n\varphi(n)/t$ for any positive integer $n$ and any $t$ with $n\leq t<n+1$. Restricted to positive integers the function $\varphi(t)$ equals the sequence $\varphi(n)$. Moreover $t\varphi(t)$ is decreasing monotone and in particular $\varphi(t)$ is decreasing monotone. Let $\eta$ be a positive and neat path. Set $M:=\|B_{\eta}\|$ and let $\theta >1$ be the constant appearing in Lemma \ref{lem1s4ss3sss1}. We get a monotone decreasing sequence setting
\begin{equation}\label{eqpsifunzionediphi}
\psi_{k}:=
\frac{\te^{k}\varphi (\te^{k})}
{dM}.
\end{equation}

Since $\varphi(t)$ is decreasing monotone, Lemma \ref{lem1s4ss1sss1} implies that $\sum_{k=1}^{\infty}\psi_{k}=+\infty$. We resume the results in Lemma \ref{lem1s4ss2sss2} and Lemma \ref{lem2s4ss2sss2} with the following lemma.

\begin{lem}\label{lem3s4ss2sss2}
Let $\eta$ be a neat and positive path and consider the map $\cF_{\eta}$ defined by equation (\ref{eqdefmappaeffe}). Then for almost any $T$ and any $k$ big enough we have
$$
\psi_{k}
\leq
\frac{1}
{\|\la^{(r_{k})}\|}\varphi(\|q^{\ga(T,r_{k})}\|),
$$
where $r_{k}$ are the instants defined at the beginning of \S \ref{s4ss2sss2} and $\psi_{k}$ is the sequence defined in (\ref{eqpsifunzionediphi}) with parameter $\te$ given by Lemma \ref{lem1s4ss2sss2}.
\end{lem}

\begin{demo}
Observe that the statement in Lemma \ref{lem2s4ss2sss2} is equivalent to $q^{\ga}_{\xi}\geq\|q^{\ga}\|/dM$ for any $\xi\in\cA$, where $M:=\|B_{\eta}\|$. Thus we have
$$
\langle q^{\ga(T,r_{k})},\la^{(r_{k})} \rangle=
\sum_{\xi\in\cA}\la^{(r_{k}}_{\xi}
\langle q^{\ga(T,r_{k})},e_{\xi}\rangle\geq
\frac
{\|\la^{(r_{k})}\|\cdot\|q^{\ga(T,r_{k})}\|}
{dM}.
$$
On the other hand, observe that
$$
\langle q^{\ga(T,r_{k})},\la^{(r_{k})} \rangle=
\langle B_{\ga(T,r_{k})}\vec{1},^{t}B_{\ga(T,r_{k})}^{-1}\la \rangle=
\langle \vec{1},\la\rangle=1.
$$
Therefore we have
$$
\frac{1}{\|\la^{(r_{k})}\|}
\varphi(\|q^{\ga(T,r_{k})}\|)
\geq
\frac{\|q^{\ga(T,r_{k})}\|}{dM}
\varphi(\|q^{\ga(T,r_{k})}\|)
$$
Lemma \ref{lem1s4ss2sss2} implies that $\|q^{\ga(T,r_{k})}\|\leq \te^{k}$ for almost any $T$ and any $k$ bigger than some $k_{0}=k_{0}(T)$. Since $t\varphi(t)$ is decreasing monotone then the lemma follows.
\end{demo}

\subsubsection{Shrinking target formulation}\label{s4ss2sss3}

In this paragraph we state a shrinking-target criterion for the map $\cF_{\eta}$ which implies Proposition \ref{propcasodivergente}.

\begin{defi}\label{defreferencepath}
Let $(\be,\al)$ be a pair with $\pi_{0}^{t}(\al)>1$ and $\pi_{0}^{b}(\be)>1$. A \emph{reference path} for $(\be,\al)$ is a neat and positive path $\eta:\pi_{0}\to \pi_{1}$ starting at $\pi_{0}$ and ending in $\pi_{1}$ with the following property:
\begin{itemize}
\item
If $(\be,\al)$ satisfies property A then $\eta$ is chosen according to Lemma \ref{lem1As3ss3sss2} and contains at least $2$ arrows with winner $\al$.
\item
If $(\be,\al)$ satisfies property B and $V$ is the letter appearing in Definition \ref{defproprietaAeB}, then $\eta$ is chosen according to Lemma \ref{lem1As3ss3sss2} and contains at least $d$ arrows with winner $V$.
\end{itemize}
\end{defi}

\begin{rem}\label{rem1s4ss2sss3}
Lemmas \ref{lem1As3ss3sss2} and \ref{lem1Bs3ss3sss2} just specify the ending part (the last arrow or the last two) of the path $\eta$ that they provide, whereas they leave complete freedom in the choice of its beginning. This makes it possible to choose an appropriate $\eta$ which satisfying all the required properties.
\end{rem}

Let $\pi_{0}$ and $(\be,\al)$ be respectively a combinatorial datum and a pair as in Proposition \ref{propcasodivergente}. Let $\eta:\pi_{0}\to\pi_{1}$ be a reference path for $(\be,\al)$ as in Definition \ref{defreferencepath} and let $\cF_{\eta}:\De^{(1)}_{\pi_{1}}\to\De^{(1)}_{\pi_{1}}$ be the map defined in equation \ref{eqdefmappaeffe}.

\begin{prop}\label{props4ss2}
Let $(\be,\al)$ be a pair satisfying property A and $\eta$ be a reference path for $(\be,\al)$. If for almost any $T\in\De^{(1)}_{\pi_{1}}$ there exist infinitely many $k\in\NN$ such that
$$
\cF^{k}_{\eta}(T)\in
\{(\pi_{1},\la)\in\De^{(1)}_{\pi_{1}};\la_{\al}<\psi_{k}\}
$$
then Proposition \ref{propcasodivergente} holds for the pair $(\be,\al)$.

Let $(\be,\al)$ be a pair with property B, let $V$ and $L$ be the associated letters, and let $\eta$ be a reference path for $(\be,\al)$. If for almost any $T\in\De^{(1)}_{\pi_{1}}$ there are infinitely many $k\in\NN$ such that
$$
\cF^{k}_{\eta}(T)\in
\{
(\pi_{1},\la)\in\De^{(1)}_{\pi_{1}};
\la_{V}<\min\{\la_{L},\psi_{k}\}
\}
$$
then Proposition \ref{propcasodivergente} holds for the pair $(\be,\al)$.
\end{prop}

\begin{demo}
For a generic $T$ in $\De^{(1)}_{\pi_{0}}$ consider the sequence of instants $r_{k}=r_{k}(T)$ with $r_{1}<r_{2}<\dots$ such that $\ga(T,r_{k})$ ends with $\eta$ for any $k\in\NN$. As we argued in the beginning of \S \ref{s4ss2sss1} for almost any $T$ there are infinitely many such $r_{k}$. Fix such a $T$ and write $\widehat{T}^{(r_{k})}=\widehat{Q}^{(r_{k})}(T)$. In particular $\ga(T,r_{1})$ ends with $\eta$, thus $\widehat{T}^{(r_{1})}\in\De^{(1)}_{\pi_{1}}$. Therefore the definition of the map $\cF_{\eta}$ implies that for any $r_{k}$ as above we have
$$
\widehat{T}^{(r_{k})}=
\cF_{\eta}^{k-1}(\widehat{T}^{(r_{1})}).
$$
In particular $\widehat{T}^{(r_{k})}$ belongs to $\De^{(1)}_{\pi_{1}}$ for any $k$ and we write
$
\widehat{T}^{(r_{k})}=
(\pi_{1},\widehat{\la}^{(r_{k})})
$.

If $(\be,\al)$ satisfies property A then the reference path $\eta$ is chosen according to Lemma \ref{lem1As3ss3sss2}. Any path $\ga(T,r_{k})$ ends with $\eta$ by construction, thus the first part of Proposition \ref{props3ss3} implies that to any $r_{k}$ it corresponds an integer $n(k)<\|q^{\ga(T,r_{k})}\|$ such that the triple $(\be,\al,n(k))$ is reduced for $T$ and
$
|T^{n(k)}u^{b}_{\be}-u^{t}_{\al}|=
\la^{(r_{k})}_{\al}
$,
where $\la^{(r_{k})}$ is the non-normalized length datum. Condition (\ref{eqKhinTSI}) therefore is equivalent to $\la^{(r_{k})}_{\al}<\varphi(n(k))$. Since $n\varphi(n)$ is decreasing monotone then also $\varphi(n)$ is, hence $\varphi(\|q^{\ga(T,r_{k})}\|)\leq \varphi(n(k))$. Moreover the estimate in Lemma \ref{lem3s4ss2sss2} holds for almost any $T$ and any $k$ big enough, it follows that if there are infinitely many instants $r_{k}$ such that $\widehat{\la}^{(r_{k})}_{\al}<\psi_{k}$ then there are also infinitely many reduced triples $(\be,\al,n(k))$ satisfying (\ref{eqKhinTSI}).

If $(\be,\al)$ satisfies property B, then the reference path $\eta$ is chosen according to Lemma \ref{lem1Bs3ss3sss2}. Let $V$ and $L$ be the letters appearing in Definition \ref{defproprietaAeB}. Since any $\ga(T,r_{k})$ ends with $\eta$, then the second part of Proposition \ref{props3ss3} implies that to any $r_{k}$ it corresponds an integer $n(k)<\|q^{\ga(T,r_{k})}\|$ such that
$
|T^{n(k)}u^{b}_{\be}-u^{t}_{\al}|=
\la^{(r_{k})}_{V}
$,
moreover $(\be,\al,n(k))$ is a reduced triple for $T$ if $\la^{(r_{k})}_{V}<\la^{(r_{k})}_{L}$. A reduced solution of (\ref{eqKhinTSI}) therefore corresponds to $\la^{(r_{k})}_{V}<\min\{\la^{(r_{k})}_{V},\varphi(n(k))\}$. Arguing as in the previous case, if there are infinitely many instants $r_{k}$ such that $\widehat{\la}^{(r_{k})}_{V}<\min\{\widehat{L}^{(r_{k})}_{\al},\psi_{k}\}$ then there are also infinitely many reduced triples $(\be,\al,n(k))$ satisfying (\ref{eqKhinTSI}). The proposition is proved.
\end{demo}

\subsection{Refined shrinking targets}\label{s4ss3}

In this paragraph we treat a technical issue which is necessary to translate the shrinking target property in Proposition \ref{props4ss2} into the setting of the Borel-Cantelli lemma. We write $T=(\pi_{1},\la)$ and we parameterize the family of targets in Proposition \ref{props4ss2} with a parameter $\epsilon>0$ as follows. For a pair $(\be,\al)$ satisfying property A the general target is the set $\{T\in\De^{(1)}_{\pi_{1}};\la_{\al}<\epsilon\}$ and for a pair $(\be,\al)$ satisfying property B it is the set $\{T\in\De^{(1)}_{\pi_{1}};\la_{V}<\min\{\la_{L},\epsilon\}\}$. In both cases the $(d-1)$-volume of these sets is obviously proportional to $\epsilon$. For any $\epsilon$ we define a subset of these targets in order to satisfy these two properties.
\begin{itemize}
\item
The refined targets belong to sigma-algebra generated by the connected components of the domain of the map $\cF_{\eta}$.
\item
The $(d-1)$-volume of the refined targets is proportional to $\epsilon$.
\end{itemize}

\subsubsection{General construction}\label{s4ss3sss1}

Let $W$ be any letter in the alphabet $\cA$ and let us denote $\cA_{W}$ the sub-alphabet $\cA\setminus \{W\}$. Fix any $\pi\in\cR$ and any $\epsilon>0$.

Let $E(\pi,W,\epsilon)$ be the family of those $\cA_{W}$-colored paths $\ga$ which start at $\pi$, satisfy $q^{\ga}_{W}>1/\epsilon$ and are minimal with this property with respect to the ordering $\prec$. By minimality of its paths $E(\pi,W,\epsilon)$ is a disjoint family. Let $N(\pi,W,\epsilon)$ be the family of paths $\nu$ starting at $\pi$, which satisfy $q^{\nu}_{W}<1/\epsilon$, end with an arrow with winner $W$ and are maximal with these two properties with respect to $\prec$. By maximality of its paths, $N(\pi,W,\epsilon)$ is a disjoint family. Moreover $E(\pi,W,\epsilon)$ and $N(\pi,W,\epsilon)$ are disjoint each other and satisfy $\PP(E(\pi,W,\epsilon))+\PP(E(\pi,W,\epsilon))=1$. We define the sets
$$
\De^{(1)}(E(\pi,W,\epsilon)):=
\bigsqcup_{\ga\in E(\pi,W,\epsilon)}
\De^{(1)}_{\ga},
$$
$$
\De^{(1)}(N(\pi,W,\epsilon)):=
\bigsqcup_{\nu \in N(\pi,W,\epsilon)}
\De^{(1)}_{\ga},
$$
where the unions are disjoint because of disjointness of the corresponding families of paths. Condition
$
\PP(E(\pi,W,\epsilon))+\PP(E(\pi,W,\epsilon))=1
$
is equivalent to
$$
\De^{(1)}_{\pi}=
\De^{(1)}(E(\pi,W,\epsilon))
\sqcup
\De^{(1)}(N(\pi,W,\epsilon))
\mod 0.
$$
\begin{lem}\label{lem1s4ss3sss1}
For any letter $W \in \cA$, any $\pi\in \cR$ and any $\epsilon >0$ we have
$$
\De^{(1)}(E(\pi,W,\epsilon))
\subset
\{T\in\De^{(1)}_{\pi};\la_{W}<\epsilon\}.
$$
\end{lem}

\begin{demo}
It is enough to prove that for any $\ga$ in $E(\pi,W,\epsilon)$ and for any $T$ in $\De^{(1)}_{\ga}$ we have $\la_{W}<\epsilon$, where we write $T=(\pi,\la)$. We observe that the length datum $\la$ of any $T$ in $\De^{(1)}_{\ga}$ is a convex combination of the vectors $v_{\xi}:=(1/q^{\ga}_{\xi})^{t}B_{\ga}e_{\xi}$ with $\xi\in\cA$. Any of these vectors has $W$-coordinate equal to
$$
\langle v_{\xi},e_{W}\rangle=
\frac
{\langle ^{t}B_{\ga}e_{\xi},e_{W}\rangle}
{q^{\ga}_{\xi}}=
\frac
{\langle e_{\xi},B_{\ga}e_{W}\rangle}
{q^{\ga}_{\xi}}.
$$
The scalar product above is maximum for the vertex $v_{W}$ of the simplex $\De^{(1)}_{\ga}$, with value $\langle B_{\ga}e_{W},e_{W}\rangle (q^{\ga}_{W})^{-1}$. Moreover the letter $W$ newer wins in $\ga$ by definition of $E(\pi,W,\epsilon)$, so $B_{\ga}e_{W}=e_{W}$ and the maximum value of the $W$-coordinate for length data in $\De^{(1)}_{\ga}$ is $(q^{\ga}_{W})^{-1}$, which is smaller than $\epsilon$ for $\ga$ in $E(\pi,W,\epsilon)$. The lemma is proved.
\end{demo}

In Theorem \ref{propstimalocalegenerale} we prove that there exists a positive constant $C$, depending only on the number of intervals $d$, such that for any $\epsilon>0$ we have
$$
\leb_{d-1}
(\De^{(1)}(E(\pi,W,\epsilon)))
\geq C\epsilon.
$$

\subsubsection{$\eta$-measurability}\label{s4ss3sss2}

Let $\eta$ be any neat path starting at $\pi_{0}$ and ending in $\pi_{1}$ and consider the map $\cF_{\eta}$ defined in equation (\ref{eqdefmappaeffe}). We say that a finite path $\nu$ starting at $\pi_{1}$ is $\eta$\emph{-measurable} if the simplex $\De^{(1)}_{\nu}$ is measurable with respect to the sigma-algebra generated by the connected components of the domain of the map $\cF_{\eta}$, that is the sigma-algebra whose atoms are the simplices $\De^{(1)}_{\ga}$ for $\ga\in\Ga^{\eta}$. It is easy to see that a finite path $\nu$ starting at $\pi_{1}$ is $\eta$-measurable if and only if it does not contain $\eta$ as a sub-path. Then we say that a subset $X$ of $\De^{(1)}_{\pi}$ is $\eta$-measurable if there exists a family $\Ga=\Ga(X)$ of $\eta$-measurable pats such that $X=\bigsqcup_{\nu\in\Ga}\De^{(1)}_{\nu}$.

\subsubsection{Pairs with property A}\label{s4ss3sss3}

Consider a pair $(\be,\al)$ satisfying property A and let $\eta$ be a reference path for $(\be,\al)$ starting at $\pi_{0}$ and ending at $\pi_{1}$, as in Definition \ref{defreferencepath}. For any $\epsilon>0$ we define two families of paths just setting $\cE^{A}(\pi_{1},\epsilon):=E(\pi_{1},\al,\epsilon)$ and $\cN^{A}(\pi_{1},\epsilon):=N(\pi_{1},\al,\epsilon)$, where $E(\pi_{1},\al,\epsilon)$ and $N(\pi_{1},\al,\epsilon)$ are defined in \S \ref{s4ss3sss1}. For $(\be,\al)$ with property A the refined shrinking target is the set
$$
\De^{(1)}(\cE^{A}(\pi_{1},\epsilon)):=
\bigsqcup_{\ga\in\cE^{A}(\pi_{1},\epsilon)}
\De^{(1)}_{\ga},
$$
where the union above is disjoint because $\cE^{A}(\pi_{1},\epsilon)$ is a disjoint family. The complement of the target is
$$
\De^{(1)}(\cN^{A}(\pi_{1},\epsilon)):=
\bigsqcup_{\nu\in\cN^{A}(\pi_{1},\epsilon)}
\De^{(1)}_{\nu}.
$$
It follows directly from the definition that $\cE^{A}(\pi_{1},\epsilon)$ and $\cN^{A}(\pi_{1},\epsilon)$ are disjoint with $\PP(\cE^{A}(\pi_{1},\epsilon))+\PP(\cE^{A}(\pi_{1},\epsilon))=1$, thus
$$
\De^{(1)}_{\pi_{1}}=
\De^{(1)}(\cE^{A}(\pi_{1},\epsilon))
\sqcup
\De^{(1)}(\cN^{A}(\pi_{1},\epsilon))
\mod 0.
$$

\begin{prop}\label{propAs4ss3}
Let $(\be,\al)$ be a pair satisfying property A and $\eta$ be a reference path for $(\be,\al)$. For any $\epsilon>0$ all paths in the families $\cE^{A}(\pi_{1},\epsilon)$ and $\cN^{A}(\pi_{1},\epsilon)$ are $\eta$-measurable. Moreover the refined target $\De^{(1)}(\cE^{A}(\pi_{1},\epsilon))$ is contained in the set $\{T\in\De^{(1)}_{\pi_{1}};\la_{\al}<\epsilon\}$ and satisfies $\leb_{d-1}(\De^{(1)}(\cE^{A}(\pi_{1},\epsilon))\geq C\epsilon$, where $C$ is the constant in Theorem \ref{propstimalocalegenerale}.
\end{prop}

\begin{demo}
The last part of the proposition follows directly from Lemma \ref{lem1s4ss3sss1} and Theorem \ref{propstimalocalegenerale}, thus we just have to check $\eta$-measurability. As we argue in \S \ref{s4ss3sss2}, $\eta$-measurable paths are those which do not contain $\eta$ as sub-path. We recall from Definition \ref{defreferencepath} that for $(\be,\al)$ satisfying property A any reference path $\eta$ contains at least $2$ arrows with winner $\al$. The required property follows observing that $\al$ never wins in paths $\ga\in \cE^{A}(\pi_{1},\epsilon)$ and wins just once in paths $\nu\in \cN^{A}(\pi_{1},\epsilon)$. The proposition is proved.
\end{demo}

\subsubsection{Pairs with property B}\label{s4ss3sss4}

Let $(\be,\al)$ be a pair of letters satisfying property B and let $V$ and $L$ be the associated letters as in Definition \ref{defproprietaAeB}. Let $\eta$ be a reference for $(\be,\al)$ starting in $\pi_{0}$ and call $\pi_{1}$ its ending point. According to Lemma \ref{lem1Bs3ss3sss2} the element $\pi(\be,\al)$ associated to $(\be,\al)$ in Definition \ref{defproprietaAeB} is in second to last position in $\eta$ and $\al$ wins in the last arrow of $\eta$. It follows that $\pi_{1}$ satisfies
$\{\xi\in \cA\,;\, \pi_{1}^{t}(\xi)<\pi_{1}^{t}(\al) \}\cup \{V\}= \{\xi\in \cA\,;\, \pi_{1}^{b}(\xi)<\pi_{1}^{b}(\be) \}$ and
$\pi_{1}^{t}(V)=\pi_{1}^{b}(\be)$, that is
$$
\pi_{1}=
\left(
\begin{array}{ccccccc}
\dots &   L    &  \dots  &   \al  &    V    &  \dots &     \\
\dots &   V    &  \dots  &    L   &   \be   &  \dots & \al \\
\end{array}
\right).
$$
Define the sub-alphabet $\cA':=\{\xi\in \cA\,;\,\pi^{t}_{1}(\xi)<\pi^{t}_{1}(\al)\}$ and call $a$ the number of elements of $\cA'$. Let $\cR_{\ast}$ be the essential $(\cA\setminus\cA')$-decorated Rauzy class which contains $\pi_{1}$ (see \S \ref{Reduction of rauzy classes} for a description of the formalism for reduction of Rauzy classes).

\begin{rem}\label{rem1s4ss3sss4}
Since $\pi_{1}$ is an essential element of $\cR_{\ast}$ then there exists a good letter for $\pi_{1}$, furthermore there exists only one good letter and it is evident that it is $V$. Therefore any path $\ga$ in the set $E(\pi_{1},V,\epsilon)$ (see Definition \ref{s4ss3sss1}) is $\cA'$-separated, indeed by definition $V$ never wins in $\ga$.
\end{rem}

\begin{lem}\label{lem1s4ss3sss4}
Let $\hat{\ga}:\pi_{1}\to\hat{\pi}$ be any $\cA'$-separated path starting at $\pi_{1}$ and ending in $\hat{\pi}$. Then for any letter $\xi\in \cA'\cup \{\al\}$ we have:
$$
\hat{\pi}^{t}(\xi)=\pi_{1}^{t}(\xi).
$$
Moreover if $\hat{\ga}:\pi_{1}\to\hat{\pi}$ is maximal $\cA'$-separated, then its ending point $\hat{\pi}$ satisfies $\hat{\pi}^{b}(L)=d$ and we have:
$$
B_{\hat{\ga}}e_{V}=\sum_{\xi\in\cA\setminus\cA'}e_{\xi}.
$$
\end{lem}

\begin{demo}
Suppose that there exists an $\cA'$-separated path $\hat{\ga}$ starting at $\pi_{1}$ and ending in $\hat{\pi}$ which does not satisfy the first part of the lemma. We can suppose that $\hat{\ga}$ is minimal with this property, that is $\hat{\pi}$ is the first element in $\hat{\ga}$ where the condition does not hold. By minimality there exists a letter $\xi\in \cA'\cup\{\al\}$ such that $\hat{\pi}^{t}(\xi)=\pi_{1}^{t}(\xi)+1$. Let $\ga_{last}$ be the last arrow in $\hat{\ga}$, call $W$ its winner and $\pi$ its starting point, i.e. $\pi$ is in second-to-last position in $\hat{\ga}$. Observe that $\ga_{last}$ has to be a bottom arrow and its starting point $\pi$ has to satisfy $\pi^{t}(W)<\pi^{t}(\xi)$. Since $\hat{\ga}$ is $\cA'$-separated then $W\in\cA\backslash\cA'$, therefore $\pi$ still doesn't satisfy condition in the first part of the lemma, which is absurd by minimality of $\hat{\ga}$.

Now let us consider a maximal $\cA'$-separated path $\hat{\ga}$ starting at $\pi_{1}$ and ending in $\hat{\pi}$. By maximality of $\hat{\ga}$ there exists a letter $\xi\in\cA'$ such that $\hat{\pi}^{t}(\xi)=d$ or $\hat{\pi}^{b}(\xi)=d$. By the first part of the lemma the only possibility is $\hat{\pi}^{b}(\xi)=d$. Moreover $L$ is the rightmost letter of $\cA'$ in the permutation $\pi_{1}$ and in order to invert its position with respect to any other letter of $\xi\in \cA'$ it has to arrive in last position in the bottom line at least one time. Since $\hat{\ga}$ is $\cA'$-separated this can happen only at its ending point $\hat{\pi}$, therefore we have $\hat{\pi}^{b}(L)=d$. To prove the second part of the lemma let us decompose $\hat{\ga}$ as $\hat{\ga}=\ga^{(1)}\ga_{1},\dots,\ga^{(m)}\ga_{m}$, where $m=d-a-1$ and for any $i=1,\dots,m$ the sub-path $\ga^{(i)}$ is not drifting and $\ga_{i}$ is a drifting arrow. Let us write $\hat{\ga}^{(i)}:=\ga^{(1)}\ga_{1},\dots,\ga^{(i)}$. For any $i=1,\dots,m$ call $\al_{i}$ and $\be_{i}$ respectively the winner and the loser of the arrow $\ga_{i}$, then call $\pi_{s}^{(i)}$ and $\pi_{e}^{(i)}$ respectively the starting and ending point of $\ga_{i}$. Since $V$ is the only good letter for $\pi$ we have $\al_{1}=V$. Then we have
$$
B_{\ga^{(1)}}e_{V}=e_{V}
\texttt{ and }
B_{\ga_{1}}e_{V}=e_{V}+e_{\be_{1}}
$$
and the only good letters for $\pi_{e}^{(1)}$ are $V$ and $\be_{1}$. Let us put $\cI_{0}:=\{V\}$. Let us fix $k\leq m$ and suppose by induction that for any $1\leq i<k$ we have that there exists a subset $\cI_{i}\subset \cA\setminus\cA'$ with $i+1$ elements and such that
$$
B_{\hat{\ga}^{(i)}}(e_{V})=
\sum_{\xi\in \cI_{i-1}}e_{\xi}
\texttt{ and }
B_{\hat{\ga}^{(i)}\ga_{i}}(e_{V})=
\sum_{\xi\in \cI_{i}}e_{\xi}
$$
and such that the good letters for $\pi_{e}^{(i)}$ are exactly the letters of $\cI_{i}$. We observe that the first step of the induction is satisfied by $\cI_{0}$. Let us consider the path $\hat{\ga}^{(k)}=\hat{\ga}^{(k-1)}\ga_{k-1}\ga^{(k)}$. By the induction hypothesis $B_{\hat{\ga}^{(k-1)}\ga_{k-1}}(e_{V})=\sum_{\xi\in \cI_{k-1}}e_{\xi}$ and the good letters for $\pi_{e}^{(k-1)}$ are exactly the letters of $\cI_{k-1}$. Since there is no drift for any arrow in $\ga^{(k)}$, then the winner of such arrows is never in $\cI_{k-1}$, therefore we have  $B_{\hat{\ga}^{(k)}}(e_{V})=\sum_{\xi\in \cI_{k-1}}e_{\xi}$. Now let us consider the $k$-th drifting arrow $\ga_{k}$, its winner $\al_{k}$ and its loser $\be_{k}$. The first part of the lemma says that all the drifting arrows of $\hat{\ga}$ are top arrows, therefore $\pi_{s}^{(k),b}(\be_{k})=d$ and $\be_{k}$ is not an element of $\cI_{k-1}$. Moreover since $\ga_{k}$ is drifting we have $\pi_{s}^{(k),t}(\al_{k})=d$ and $\pi_{s}^{(k),b}(\al_{k})<d_{b}(\pi_{s}^{(k)})$ (see paragraph \ref{driftinessentialdecoratedrauzyclasses} for the notation), therefore $\pi_{e}^{(k),b}(\be_{k})=\pi_{s}^{(k),b}(\al_{k})+1$, that is $\be_{k}$ moves in good position for $\pi_{e}^{(k)}$. Putting $\cI_{k}:=\cI_{k-1}\cup \{\be_{k}\}$ the inductive step follows. The lemma is proved.
\end{demo}

Fix any path $\ga$ in $E(\pi_{1},V,\epsilon)$ and observe that it is $\cA'$-separated, according to remark \ref{rem1s4ss3sss4}. Let $\cE_{\ga}$ be the family of those maximal $\cA'$-separated paths $\hat{\ga}$ starting at $\pi_{1}$ and such that $\ga\prec\hat{\ga}$ (that is $\hat{\ga}$ begins with $\ga$). Any $\hat{\ga}$ in $\cE_{\ga}$ has its ending point $\hat{\pi}$ in the essential $(\cA\setminus\cA')$-decorated class $\cR_{\ast}$. Lemma \ref{lem1s4ss3sss4} implies $\hat{\pi}^{t}(L)=\pi^{t}_{1}(L)\leq a$ and $\hat{\pi}^{b}(L)=d$. Therefore there exists a path $\eta(\hat{\ga})$ starting at $\hat{\pi}$ which is concatenation of $d-a$ bottom arrows $\eta_{1},\dots,\eta_{d-a}$, each one with winner $L$ and such that any $\xi$ in $\cA\setminus\cA'$ is the loser of some $\eta_{i}$. For $i\in\{1,\dots,d-a-1\}$ let $\sigma_{i}$ be the path starting at the point where $\eta_{i}$ ends and with looser $L$ and define $\nu_{i}:=\eta_{1}\dots\eta_{i}\sigma_{i}$, that is the concatenation of the first $i$ arrows with winner $L$ followed by the arrow where $L$ loses. Finally denote $\cN(\hat{\ga}):=\{\hat{\ga}\nu_{1},\dots,\hat{\ga}\nu_{d-a}\}$. We define the family of paths

$$
\cE^{B}(\pi_{1},\epsilon):=
\bigsqcup_{\ga\in E(\pi_{1},V,\epsilon)}
\bigsqcup_{\hat{\ga}\in\cE_{\ga}}
\hat{\ga}\eta(\hat{\ga}),
$$
$$
\cN^{B}(\pi_{1},\epsilon):=
N(\pi_{1},V,\epsilon)\sqcup
\bigg(
\bigsqcup_{\ga\in E(\pi_{1},V,\epsilon)}
\bigsqcup_{\hat{\ga}\in \cE_{\ga}}
\cN(\hat{\ga})
\bigg).
$$
For a pair $(\be,\al)$ satisfying property B the refined target is the set
$$
\De^{(1)}(\cE^{B}(\pi_{1},\epsilon)):=
\bigsqcup_{\ga' \in \cE^{B}(\pi_{1},\epsilon)}
\De^{(1)}_{\ga'},
$$
where the union is disjoint since $\cE^{B}(\pi_{1},\epsilon)$ is a disjoint family. Similarly the complement of the target is
$$
\De^{(1)}(\cN^{B}(\pi_{1},\epsilon)):=
\bigsqcup_{\nu'\in\cN^{B}(\pi_{1},\epsilon)}
\De^{(1)}_{\nu'}.
$$
$\cE^{B}(\pi_{1},\epsilon)$ and $\cE^{B}(\pi_{1},\epsilon)$ are disjoint each other with $\PP(\cE^{B}(\pi_{1},\epsilon))+\PP(\cE^{B}(\pi_{1},\epsilon))=1$, hence
$$
\De^{(1)}_{\pi_{1}}=
\De^{(1)}(\cE^{B}(\pi_{1},\epsilon))
\sqcup
\De^{(1)}(\cN^{B}(\pi_{1},\epsilon))
\mod 0.
$$

\begin{lem}\label{lem2s4ss3sss4}
Let $(\be,\al)$ be a pair of letters satisfying property B, let $V$ and $L$ be the associated letters as in Definition \ref{defproprietaAeB} and let $\eta:\pi_{0}\to\pi_{1}$ be a reference path for $(\be,\al)$. Then we have
$$
\De^{(1)}(\cE^{B}(\pi_{1},\epsilon))
\subset
\{
T\in\De^{(1)}_{\pi_{1}};
\la_{V}<\min \{\la_{L},\epsilon\}
\}.
$$
\end{lem}

\begin{demo}
We want to prove that for any $\ga'$ in $\cE^{B}(\pi_{1},\epsilon)$ and any $T=(\pi_{1},\la)$ in $\De^{(1)}_{\ga'}$ we have $\la_{V}<\min \{\la_{L},\epsilon\}$. Any $\ga'$ in $\cE^{B}(\pi_{1},\epsilon)$ begins with a path $\ga\in E(\pi_{1},V,\epsilon)$ and Lemma \ref{lem1s4ss3sss1} implies that for such $\la$ we have $\la_{V}<\epsilon$, therefore it is enough to prove that $\la_{V}<\la_{L}$. Consider $\ga'$ in $\cE^{B}(\pi_{1},\epsilon)$ and decompose it as $\ga'=\hat{\ga}\eta(\hat{\ga})$, where $\hat{\ga}\in\cE_{\ga}$ for some $\ga$ in $E(\pi_{1},V,\epsilon)$. We have
$$
B_{\ga'}e_{V}=
B_{\eta(\hat{\ga})}
\bigg(
\sum_{\xi\in \cA\setminus\cA'}e_{\xi}
\bigg)=
\sum_{\xi\in\cA\setminus\cA'}e_{\xi},
$$
where the first equality is consequence of the second part of Lemma \ref{lem1s4ss3sss4} and the second equality holds because the winner of any arrow composing $\eta(\hat{\ga})$ is $L$, which does not belong to $\cA\setminus\cA'$. On the other hand we have
$$
B_{\ga'}e_{L}=
B_{\eta(\hat{\ga})}e_{L}=
e_{L}+\bigg(\sum_{\xi\in\cA\setminus\cA'}e_{\xi}\bigg).
$$
Here the first equality follows since $\hat{\ga}$ is $\cA'$-separated and thus it does not contain arrows with winner $L$. The second equality follows because the ending point $\hat{\pi}$ of $\hat{\ga}$ satisfies $\{\xi\in\cA;\hat{\pi}^{t}(\xi)>a\}=\cA\setminus\cA'$, according to the first part of Lemma \ref{lem1s4ss3sss4}, and on the other hand any of these letters loses against $L$ in some arrow composing $\eta(\hat{\ga})$. Summing up, we proved that any $\ga'$ in $\cE^{B}(\pi_{1},\epsilon)$ satisfies $B_{\ga'}(e_{L}-e_{V})=e_{L}$. Consider the vertices $v_{\xi}$ of the simplex $\De^{(1)}_{\ga'}$. They satisfy $(q^{\ga'}_{\xi})v_{\xi}=^{t}B_{\ga'}e_{\xi}$ for any $\xi$ in $\cA$, therefore
$$
\langle v_{\xi},e_{L}-e_{V}\rangle=
(q^{\ga'}_{\xi})^{-1}
\langle ^{t}B_{\ga}e_{\xi},e_{L}-e_{V}\rangle =
(q^{\ga'}_{\xi})^{-1}
\langle e_{\xi},e_{L}\rangle\geq 0.
$$
In particular $\langle v_{L},e_{L}-e_{V}\rangle>0$. Since any $\la$ in $\De^{(1)}_{\ga'}$ is convex combination of the vertices $v_{\xi}$, then $\langle\la,e_{L}-e_{V}\rangle>0$. The lemma is proved.
\end{demo}

\begin{lem}\label{lem3s4ss3sss4}
There exists a positive constant $C'>0$, depending only on the number of intervals $d$, such that for any pair of letters $(\be,\al)$ satisfying property B and for any $\epsilon >0$ we have
$$
\PP(\cE^{B}(\pi_{1},\epsilon))\geq C'\epsilon.
$$
\end{lem}

\begin{demo}
Consider any $\ga$ in $E(\pi_{1},V,\epsilon)$ and denote $\cE^{B}(\pi_{1},\epsilon|\ga)$ the set of those paths $\ga'$ in $\cE^{B}(\pi_{1},\epsilon)$ which begins with $\ga$. We have
$$
\PP(\cE^{B}(\pi_{1},\epsilon))=
\sum_{\ga\in E(\pi_{1},V,\epsilon)}
\PP(\ga)\PP_{\ga}(\cE^{B}(\pi_{1},\epsilon|\ga)).
$$
Then it is enough to prove that for any $\ga$ in $E(\pi_{1},V,\epsilon)$ we have
$$
\PP_{\ga}(\cE^{B}(\pi_{1},\epsilon|\ga))
\geq
\frac{1}{2^{d-a}},
$$
combining this estimate with Theorem \ref{propstimalocalegenerale} we get the required estimate with $C':=C2^{-(d-a)}$, where $C$ is the constant appearing in Theorem \ref{propstimalocalegenerale}. We recall that for a fixed $\ga$ in $E(\pi_{1},V,\epsilon)$ any $\ga'$ in $\cE^{B}(\pi_{1},\epsilon|\ga)$ is decomposed as $\ga'=\hat{\ga}\eta(\hat{\ga})$ with $\hat{\ga}$ in $\cE_{\ga}$. We have
$$
\PP_{\ga}(\cE^{B}(\pi_{1},\epsilon|\ga))=
\sum_{\hat{\ga}\in\cE_{\ga}}
\PP_{\ga}(\hat{\ga}\eta(\hat{\ga}))=
\sum_{\hat{\ga}\in\cE_{\ga}}
\PP_{\ga}(\hat{\ga})\PP_{\hat{\ga}}(\eta(\hat{\ga}))
\geq
\inf_{\hat{\ga}\in\cE_{\ga}}
\PP_{\hat{\ga}}(\eta(\hat{\ga}))
$$
because $\{\De^{(1)}_{\hat{\ga}};\hat{\ga}\in\cE_{\ga}\}$ form a partition $\mod 0$ of $\De^{(1)}_{\ga}$. For any $\hat{\ga}\in\cE_{\ga}$ the path $\eta(\hat{\ga})$ is concatenation of $d-a$ bottom arrows with winner $L$, moreover any letter in $\cA\setminus \cA'$ loses in exactly one of these arrows. Therefore for the concatenation $\ga'=\hat{\ga}\eta(\hat{\ga})$ we have
$$
q^{\ga'}_{\xi}=
q^{\hat{\ga}}_{\xi}+q^{\hat{\ga}}_{L}
\textrm{ for }
\xi\in \cA\setminus\cA'
\textrm{ and }
q^{\ga'}_{\xi}=
q^{\hat{\ga}}_{\xi}
\textrm{ for }
\xi\in \cA'.
$$
Any $\hat{\ga}$ in $\cE_{\ga}$ is $\{L\}$-separated, thus we have $q^{\hat{\ga}}_{L}=1$ and it follows that $q^{\ga'}_{\xi}=q^{\hat{\ga}}_{\xi}+1\leq 2q^{\hat{\ga}}_{\xi}$ for any $\xi\in\cA\setminus\cA'$. Applying Equation (\ref{eqprobabilitacondizionatasimplessi}) we get
$$
\PP_{\hat{\ga}}(\eta(\hat{\ga}))=
\prod_{\xi\in\cA}
\frac{q^{\hat{\ga}}_{\xi}}{q^{\ga'}_{\xi}}>
2^{-(d-a)},
$$
which implies $\PP_{\ga}(\cE^{B}(\pi_{1},\epsilon|\ga))\geq 2^{-(d-a)}$. The proposition is proved.
\end{demo}

\begin{prop}\label{propBs4ss3sss4}
Let $(\be,\al)$ be a pair satisfying property B and $\eta$ be a reference path for $(\be,\al)$. For any $\epsilon>0$ all paths in the families $\cE^{B}(\pi_{1},\epsilon)$ and $\cN^{B}(\pi_{1},\epsilon)$ are $\eta$-measurable. Moreover the refined target $\De^{(1)}(\cE^{B}(\pi_{1},\epsilon))$ is contained in the set $\{T\in\De^{(1)}_{\pi_{1}};\la_{V}<\min\{\la_{L},\epsilon\}\}$ and satisfies $\leb_{d-1}(\De^{(1)}(\cE^{B}(\pi_{1},\epsilon))\geq C'\epsilon$, where $C'$ is the constant in Lemma \ref{lem3s4ss3sss4}.
\end{prop}

\begin{demo}
The letter $V$ wins at most $d-a+1$ times in paths $\ga'\in\cE(\pi_{1},V,\epsilon)$ and $\nu'\in\cN(\pi_{1},\al,\epsilon)$, thus such $\ga'$ and $\nu'$ cannot contain a sub-path $\eta$ as in the statement. Recalling from \S \ref{s4ss3sss2} that $\eta$-measurable paths are those which do not contain $\eta$ as sub-path we get that $\cE^{B}(\pi_{1},\epsilon)$ and $\cN^{B}(\pi_{1},\epsilon)$ are $\eta$-measurable. The last part of the proposition follows from Lemma \ref{lem2s4ss3sss4} and Lemma \ref{lem3s4ss3sss4}. The proposition is proved.
\end{demo}

\subsection{Proof of the divergent case}\label{s4ss4}

In this subsection we complete the proof of Proposition \ref{propcasodivergente}. We recall that Proposition \ref{props4ss2} establishes a sufficient shrinking target criterion in terms of the map $\cF_{\eta}$, where $\eta$ is a reference path for some pair $(\be,\al)$. Here we prove the shrinking target property for the refined targets constructed in \S \ref{s4ss3}.

\subsubsection{Borel-Cantelli formulation}\label{s4ss4sss1}

Let $\{X_{k}\}_{k\in\NN}$ be any countable family of sub-sets of some set $X$. We put
$$
\limsup_{k}X_{k}:=\bigcap_{k\geq0}\bigcup_{i\geq k}X_{i}.
$$

From now on we treat jointly all pairs $(\be,\al)$ with $\pi_{0}^{t}(\al)>1$ and $\pi_{0}^{b}(\be)>1$, thus we introduce a simplified and unified notation. Recall that any such pair satisfies either property A or property B, according to Theorem \ref{propcombinatoria}. Consider a reference path $\eta$ for $(\be,\al)$, starting at $\pi_{0}$ and ending in $\pi_{1}$. Let $\cF_{\eta}:\De^{(1)}_{\pi_{1}}\to \De^{(1)}_{\pi_{1}}$ be the map introduced in equation (\ref{eqdefmappaeffe}) and $\psi_{k}$ be the sequence defined in (\ref{eqpsifunzionediphi}). Fix any $k\in\NN$. If the pair $(\be,\al)$ has property A we set $\cE_{k}:=\cE^{A}(\pi_{1},\psi_{k})$ and $\cN_{k}:=\cN^{A}(\pi_{1},\psi_{k})$. Otherwise if $(\be,\al)$ has property B we set $\cE_{k}:=\cE^{B}(\pi_{1},\psi_{k})$ and $\cN_{k}:=\cN^{B}(\pi_{1},\psi_{k})$. We also introduce the simplified notation:
$$
\De^{(1)}(\cE_{k}):=
\bigsqcup_{\ga \in \cE_{k} }
\De^{(1)}_{\ga}
\texttt{ and }
\De^{(1)}(\cN_{k}):=
\bigsqcup_{\ga \in \cN_{k}}
\De^{(1)}_{\ga}.
$$
Proposition \ref{propAs4ss3} and Proposition \ref{propBs4ss3sss4} imply that any path in the families $\cE_{k}$ and $\cN_{k}$ is $\eta$-measurable, that is $\De^{(1)}(\cE_{k})$ and $\De^{(1)}(\cN_{k})$ are $\eta$-measurable sets.

\begin{prop}\label{props4ss4}
Let $(\be,\al)$ be a pair of letter as in Theorem \ref{teoremaa}. Let $\eta$ be a reference path for $(\be,\al)$ and $\cF_{\eta}$ be the associated map. Then we have
$$
\leb_{d-1}(\limsup_{k}\cF^{-k}_{\eta}(\cE_{k}))=1.
$$
\end{prop}

Proposition \ref{props4ss4} implies that almost any $T\in\De^{(1)}_{\pi_{1}}$ belongs to $\limsup_{k}\cF^{-k}_{\eta}(\cE_{k})$, that is there exist infinitely many $k\in\NN$ such that $T\in\cF_{\eta}^{-k}\De^{(1)}(\cE_{k})$. If $(\be,\al)$ satisfies property A, according to Proposition \ref{propAs4ss3} for any such $k$ we have
$$
\cF^{k}_{\eta}(T)\in
\{(\pi_{1},\la)\in\De^{(1)}_{\pi_{1}};\la_{\al}<\psi_{k}\}.
$$
Otherwise, if $(\be,\al)$ satisfies property B, Proposition \ref{propBs4ss3sss4} implies that for any such $k$ we have
$$
\cF^{k}_{\eta}(T)\in
\{(\pi_{1},\la)\in\De^{(1)}_{\pi_{1}};\la_{V}<\min\{\la_{L},\psi_{k}\}\}.
$$
In both cases Proposition \ref{props4ss2} implies Proposition \ref{propcasodivergente} for the pair $(\be,\al)$. The divergent part of Theorem \ref{teoremaa} therefore follows.

\subsubsection{Bounded distortion for positive paths}\label{s4ss4sss2}

Proposition \ref{props4ss4} corresponds to the divergent part of Borel-Cantelli lemma, its proof requires some extra properties of $\cF_{\eta}$, or equivalently of the reference path $\eta$.

\begin{lem}\label{lem1s4ss4sss2}
Let $\eta$ be a positive path ending in $\pi$ and set $M:=\|B_{\eta}\|$. Then for any path $\ga$ ending with $\eta$, if $P_{\ga}$ is the probability measure on $\De^{(1)}_{\pi}$ defined in Equation (\ref{eqprobabilitacondizionatasimplessi}), then we have
$$
\frac{1}{M^{d}}\leq
\|\frac
{dP_{\ga}}
{d\leb_{d-1}}\|
\leq M^{d}.
$$
\end{lem}

\begin{demo}
Equation (\ref{eqprobabilitacondizionatasimplessi}) in paragraph \ref{thedistortionestimate} implies that for any Rauzy path $\nu$ starting at $\pi$ we have
$$
\frac
{P_{\ga}(\De^{(1)}_{\nu})}
{\leb_{d-1}(\De^{(1)}_{\nu})}=
\frac
{\prod_{\xi\in\cA}q^{\ga}_{\xi}q^{\nu}_{\xi}}
{\prod_{\xi\in\cA}q^{\ga\nu}_{\xi}}.
$$
According to Lemma \ref{lem2s4ss2sss2} we have $q^{\ga}_{\xi}<Mq^{\ga}_{\xi'}$ for any $\xi,\xi'\in\cA$. Moreover the concatenation of $\nu$ with $\eta$ gives the relation $q^{\ga\nu}=B_{\nu}q^{\ga}$. It follows that for any $\xi\in\cA$ we have $M^{-1}q^{\ga}_{\xi}q^{\nu}_{\xi}\leq q^{\ga\nu}_{\xi}\leq M q^{\ga}_{\xi}q^{\nu}_{\xi}$ and therefore
$$
\frac{1}{M^{d}}\leq
\frac{P_{\ga}(\De^{(1)}_{\nu})}{\leb_{d-1}(\De^{(1)}_{\nu})}
\leq M^{d}.
$$
When $\nu$ varies among all Rauzy paths starting at $\pi$, the sub-simplices $\De^{(1)}_{\nu}$ form a basis of the Borel sigma-algebra of $\De^{(1)}_{\pi}$, the lemma therefore is proved.
\end{demo}

For any $k\in\NN$ consider the family $\Ga^{(k),\eta}$ of those paths $\ga_{k}$ such that $\De^{(1)}_{\ga_{k}}$ is a connected component of the domain of $\cF_{\eta}^{k}$.

\begin{lem}\label{lem2s4ss4sss2}
There exists a constant $C>0$, depending only on the number of intervals $d$, such that for any $k\in\NN$ and for any $\ga_{k}\in\Ga^{(k),\eta}$ we have
$$
P_{\ga_{k}}(\De^{(1)}(\cN_{k}))
\leq
(1-C\psi_{k}).
$$
\end{lem}
\begin{demo}
Since
$
\De^{(1)}_{\pi_{1}}=
\De^{(1)}(\cE_{k})\sqcup\De^{(1)}(\cN_{k})
\mod 0
$
the statement is equivalent to
$$
P_{\ga_{k}}(\De^{(1)}(\cE_{k}))
\geq
C\psi_{k}.
$$
Set $M:=\|B_{\eta}\|$. By definition of $\cF_{\eta}$ any path $\ga_{k}\in\Ga^{(k),\eta}$ ends with $\eta$, therefore Lemma \ref{lem1s4ss4sss2} applies and we get $\| \frac{dP_{\ga_{k}}}{d\leb_{d-1}}\|\leq M^{d}$ and $\| \frac{d\leb_{d-1}}{dP_{\ga_{k}}}\|\leq M^{d}$.

Recall that $\cE_{k}$ denotes either $\cE^{A}(\pi_{1},\psi_{k})$ (defined in \S \ref{s4ss3sss3}) or $\cE^{B}(\pi_{1},\psi_{k})$ (defined in \S \ref{s4ss3sss4}). We apply directly Theorem \ref{propstimalocalegenerale} in the first case, or Lemma \ref{lem3s4ss3sss4} in the second case. In both cases we have a positive constant $C'$ such that $\leb_{d-1}(\De^{(1)}(\cE_{k}))\geq C'\psi_{k}$. Combining this last estimate with Lemma \ref{lem1s4ss4sss2} we get
$$
P_{\ga_{k}}(\De^{(1)}(\cE_{k}))
\geq
C'M^{-d}\psi_{k}.
$$
Then $C=C'M^{-d}$ is the required constant. The lemma is proved.
\end{demo}

\subsubsection{Weak independence}\label{s4ss4sss3}

We have $\De^{(1)}_{\pi_{1}}=\De^{(1)}(\cE_{k})\sqcup\De^{(1)}(\cN_{k})$ modulo a set of measure zero. For any $k\in\NN$ this is equivalent to
$$
\De^{(1)}_{\pi_{1}}=
\cF_{\eta}^{-k}\De^{(1)}(\cE_{k})
\sqcup
\cF_{\eta}^{-k}\De^{(1)}(\cN_{k})
$$
modulo a set of measure zero. In the next lemma we prove that the family of sets $\cF^{-k}(\De^{(1)}(\cN_{k}))$ with $k\in\NN$ satisfies a weak form of independence.

\begin{lem}\label{lem1s4ss4sss3}
Let $C$ be the constant appearing in Lemma \ref{lem2s4ss4sss2}. For any pair of integers $m,n$ with $m\geq n$ we have:
$$
\leb_{d-1}
\bigg(
\bigcap_{k=n}^{m}
\cF_{\eta}^{-k}\De^{(1)}(\cN_{k})
\bigg)
\leq
\prod_{k=n}^{m}(1-C\psi_{k}).
$$
\end{lem}

\begin{demo}
Fix $k\in\NN^{\ast}$. Condition $T\in\cF_{\eta}^{-k}\De^{(1)}(\cN_{k})$ is equivalent to say that there exists $\ga_{k}\in\Ga^{(k),\eta}$ and $\nu_{k}\in\cN_{k}$ such that $T\in\De^{(1)}_{\ga_{k}\nu_{k}}$. Denote $\cC(k)$ the set of concatenated paths $\ga^{(k)}=\ga_{k}\nu_{k}$ where $\ga_{k}\in\Ga^{(k),\eta}$ and $\nu_{k}\in\cN_{k}$. We have
$$
\cF_{\eta}^{-k}\De^{(1)}(\cN_{k})=
\bigsqcup_{\ga^{(k)}\in\cC(k)}
\De^{(1)}_{\ga^{(k)}},
$$
where the union is disjoint because $\Ga^{(k),\eta}$ and $\cN_{k}$ are disjoint families. With the notation introduced in \S \ref{thedistortionestimate} we can write
$
\leb_{d-1}
\big(
\cF_{\eta}^{-k}\De^{(1)}(\cN_{k})
\big)=
\PP(\cC(k))
$.
Observe that for any $\ga^{(k)}=\ga_{k}\nu_{k}$ we have
$\PP(\ga^{(k)})=\PP(\ga_{k})\PP_{\ga_{k}}(\nu^{k})$, moreover Lemma \ref{lem2s4ss4sss2} implies $\PP_{\ga_{k}}(\cN_{k})\leq (1-C\psi_{k})$ for any $\ga_{k}\in\Ga^{(k),\eta}$. Finally $\PP(\Ga^{(k),\eta})=1$, because $\cF_{\eta}$ is defined almost everywhere on $\De^{(1)}_{\pi_{1}}$. Therefore we have
$$
\PP(\cC(k))=
\sum_{\ga^{(k)}\in\cC(k)}
\PP(\ga^{(k)})
=
\sum_{\ga_{k}\in\Ga^{(k),\eta}}
\sum_{\nu_{k}\in\cN_{k}}
\PP(\ga_{k})\PP_{\ga_{k}}(\nu_{k})=
$$
$$
\sum_{\ga_{k}\in\Ga^{(k),\eta}}
\PP(\ga_{k})\PP_{\ga_{k}}(\cN_{k})
\leq
(1-C\psi_{k})
\PP(\Ga^{(k),\eta})=
(1-C\psi_{n}).
$$

Fix any integer $m>n$. Observe that for any $k\in\{n,\dots,m\}$ any path $\ga_{k}\in\Ga^{(k),\eta}$ admits a decomposition
$
\ga_{k}=\ga_{n}\ga(n+1)\dots\ga(k)
$
with $\ga_{n}\in\Ga^{(n),\eta}$ and $\ga(i)\in\Ga^{\eta}$ for any $i\in\{n+1,\dots,k\}$. Observe also that any $\nu_{k}\in\cN_{k}$ is $\eta$-measurable, that is
$
\De^{(1)}_{\nu_{k}}=
\bigsqcup_{\ga\in\Ga^{\eta}_{\nu_{n}}}
\De^{(1)}_{\ga}.
$

Condition
$
T\in
\bigcap_{k=n}^{m}
\cF_{\eta}^{-k}\De^{(1)}(\cN_{k})
$
is equivalent to say that for any $k\in\{n,\dots,m\}$ there exists $\ga_{k}\in\Ga^{(k),\eta}$ and $\nu_{k}\in\cN_{k}$ such that $T\in\De^{(1)}_{\ga_{k}\nu_{k}}$, moreover for any $k\in\{n,\dots,m-1\}$ there exists a path $\ga(\nu_{k})$ in $\Ga^{\eta}_{\nu_{k}}$ such that $\ga_{k+1}=\ga_{k}\ga(\nu_{k})$. In particular we can write $\ga_{m}=\ga_{n}\ga(\nu_{n})\dots\ga(\nu_{m-1})$. Let $\cC(n,m)$ be the set of concatenated paths
$$
\ga^{(n,m)}=\ga_{n}\ga(\nu_{n})\dots\ga(\nu_{m-1})\nu_{m},
$$
with $\ga_{n}\in\Ga^{(n),\eta}$ and where $\nu_{k}\in\cN_{k}$ for any $k\in\{n,\dots,m\}$ and $\ga(\nu_{k})\in\Ga^{\eta}_{\nu_{k}}$ for any $k\in\{n,\dots,m-1\}$. We have
$$
\bigcap_{k=n}^{m}
\cF_{\eta}^{-k}\De^{(1)}(\cN_{k})
=
\bigsqcup_{\ga^{(n,m)}\in\cC(n,m)}
\De^{(1)}_{\ga^{(n,m)}}
$$
and the statement in the lemma is equivalent to
$\PP(\cC(n,m))\leq\prod_{k=n}^{m}(1-C\psi_{k})$. Introduce the set $\cB(n,m)$ of those path $\ga_{m}$ in $\Ga^{(m),\eta}$ admitting a decomposition $\ga_{m}=\ga_{n}\ga(\nu_{n})\dots\ga(\nu_{m-1})$ with $\ga_{n}\in\Ga^{(n),\eta}$ and where $\nu_{k}\in\cN_{k}$ and $\ga(\nu_{k})\in\Ga^{\eta}_{\nu_{k}}$ for any $k\in\{n,\dots,m-1\}$. In particular fort any $\ga(n,m)$ in $\cC(n,m)$ there exist $\ga_{m}$ in $\cB(n,m)$ and $\nu_{m}$ in $\cN_{m}$ such that $\ga^{(n,m)}=\ga_{m}\nu_{m}$. For these paths we have
$
\PP(\ga^{(n,m)})=
\PP(\ga_{m})
\PP_{\ga_{m}}(\nu_{m})
$,
moreover $\ga_{m}$ belongs to $\Ga^{(m),\eta}$, thus Lemma \ref{lem2s4ss4sss2} implies $\PP_{\ga_{m}}(\cN_{m})\leq (1-C\psi_{m})$. Therefore we have
$$
\PP(\cC(n,m))=
\sum_{\ga_{m}\in\cB(n,m)}
\sum_{\nu_{m}\in\cN_{m}}
\PP(\ga_{m})
\PP_{\ga_{m}}(\nu_{m})=
$$
$$
\sum_{\ga_{m}\in\cB(n,m)}
\PP(\ga_{m})
\PP_{\ga_{m}}(\cN_{m})
\leq
(1-C\psi_{m})
\PP(\cB(n,m))
$$
On the other hand any $\ga_{m}\in\cB(n,m)$ can be decomposed as $\ga_{m}=\ga_{m-1}\ga(\nu_{m-1})$ with $\ga_{m-1}\in\cB(n,m-1)$ and where $\nu_{m-1}\in\cN_{m-1}$ and $\ga(\nu_{m-1})\in\Ga^{\eta}_{\nu_{m-1}}$. Observe that
$
\sum_{\ga\in\Ga^{\eta}_{\nu_{m-1}}}
\PP(\ga_{m-1}\ga)=
\PP(\ga_{m-1}\nu_{m-1})$,
therefore we have
$$
\PP(\cB(n,m))
=
\sum_{\ga_{m-1}\in\cB(n,m-1)}
\sum_{\nu_{m-1}\in\cN_{m-1}}
\PP(\ga_{m-1}\nu_{m-1})
=
\PP(\cC(n,m-1)).
$$
Summing up, we proved that for any $n\in\NN^{\ast}$ and any $m\geq n$ we have
$$
\PP(\cC(n,m))
\leq
(1-C\psi_{m})
\PP(\cC(n,m-1)),
$$
and iterating this last estimate $m-n$ times we get the required estimate for $\bigcap_{k=n}^{m}\cF_{\eta}^{-k}\De^{(1)}(\cN_{k})$. The lemma is proved.
\end{demo}

Recall that the sequence introduced in Equation (\ref{eqpsifunzionediphi}) satisfies $\sum_{k=1}^{+\infty}\psi_{k}=+\infty$ and this condition is equivalent to $\prod_{k=1}^{+\infty}(1-C\psi_{k})=0$ for any positive constant $C$. Fix $n\in\NN$ and take the limit for $m\to+\infty$ in the estimate of Lemma \ref{lem1s4ss4sss3}. We get
$$
\leb_{d-1}
\bigg(
\bigcap_{k=n}^{+\infty}
\cF_{\eta}^{-k}\De^{(1)}(\cN_{k})
\bigg)
=0.
$$
Then we have
$$
\leb_{d-1}
\bigg(
\bigcup_{n=1}^{\infty}
\bigg(
\bigcap_{k=n}^{+\infty}
\cF_{\eta}^{-k}\De^{(1)}(\cN_{k})
\bigg)
\bigg)
=0,
$$
therefore Proposition \ref{props4ss4} follows. The proof of the divergent case of Theorem \ref{teoremaa} is complete.

\section{Main technical results.}

\subsection{Main combinatorial property}\label{Main combinatorial property}

In this paragraph we state and prove a combinatorial property of Rauzy classes which implies that any pair $(\be,\al)$ as in Definition \ref{defproprietaAeB} satisfies either property A or property B (or both). In order to simplify notation, for a Rauzy class $\cR$ over an alphabet $\cA$, we call $X=X(\cR)$ and $Y=Y(\cR)$ the two letters in $\cA$ such that respectively $\pi^{t}(X)=1$ and $\pi^{b}(Y)=1$ for all $\pi\in \cR$.

\begin{thm}\label{propcombinatoria}
Let $\cR$ be any Rauzy class with alphabet $\cA$ and $(\be,\al)$ be any ordered pair of letters with $\be\not =Y$ and $\al\not=X$. Then at least one of the following two statements holds.
\begin{description}
\item[A] There exists an element $\pi$ in $\cR$ such that
\begin{equation}\label{eqcombinatoriasuffa}
\pi^{t}(\al)=\pi^{b}(\be)=d
\end{equation}
\item[B]  There exist two (different) elements $\pi$ and $\pi'$ in $\cR$ and two letters $V$ and $V'$ in $\cA$ such that
\begin{equation}\label{eqcombinatoriasuffbuno}
\begin{array}{ll}
\{\xi\in\cA\,;\,\pi^{t}(\xi)<\pi^{t}(\al)\}\cup \{V\}=
\{\xi\in\cA\,;\,\pi^{b}(\xi)<\pi^{b}(\be) \}\\
\pi^{t}(V)=\pi^{b}(\al)=d
\end{array}
\end{equation}
and
\begin{equation}\label{eqcombinatoriasuffbdue}
\begin{array}{ll}
\{\xi\in\cA\,;\,\pi'^{t}(\xi)<\pi'^{t}(\al)\}=
\{\xi\in\cA\,;\,\pi'^{b}(\xi)<\pi'^{b}(\be)\}\cup \{V'\}\\
\pi'^{b}(V')=\pi'^{t}(\be)=d
\end{array}
\end{equation}
\end{description}
\end{thm}

\emph{Note:} Observe that case \textbf{A} is compatible just with a pair of different letters. In case \textbf{B}, when $\be=\al$ equation (\ref{eqcombinatoriasuffbuno}) implies $\pi^{t}(\al)=d-1$ and $\pi^{b}(\al)=d$ and on the other hand equation (\ref{eqcombinatoriasuffbdue}) implies $\pi'^{t}(\al)=d$ and $\pi'^{b}(\al)=d-1$.\\

\begin{demo}
In \cite{rauzy} it is proven that any Rauzy class $\cR$ contains a \emph{standard} $\tilde{\pi}$, that is a combinatorial datum such that $\tilde{\pi}^{t}(X)=\tilde{\pi}^{b}(Y)=1$ and $\tilde{\pi}^{t}(Y)=\tilde{\pi}^{b}(X)=d$. Denote $A$ and $B$ the second letter respectively in top row and bottom row of $\tilde{\pi}$, that is
$$
\tilde{\pi}=
\left(
\begin{array}{cccc}
X & A & \dots & Y \\
Y & B & \dots & X \\
\end{array}
\right).
$$

\begin{lem}\label{lemdimteoremacombinatoriociclistandard}
Condition (\ref{eqcombinatoriasuffa}) holds for all pairs $(X,\al)$ with $\al\not=X$ and all pairs $(\be,Y)$ with $\be\not=Y$.
\end{lem}
\begin{demo}
The standard element $\tilde{\pi}$ is the base point of two loops in the Rauzy diagram of $\cR$ with length $d-1$. One of these two loops is concatenation of $d-1$ bottom arrows with winner $X$. Any letter $\al\not=X$ loses against $X$ in some arrow in this loop, therefore the lemma follows for the corresponding pair $(X,\al)$. The other loop is the concatenation of $d-1$ top arrows with winner $Y$ and with the symmetric argument we get the statement for the pairs $(\be,Y)$ with $\be\not=Y$. The lemma is proved.
\end{demo}

For other pairs of letters the proof goes on by induction on the number of letters $d$. There are four Rauzy classes with $d\leq4$ letters, which are easily computable (see the pictures of these classes in \S 1.2 of \cite{mmy}). For these classes it is easy to check directly the statement in Theorem \ref{propcombinatoria}. Therefore we consider a Rauzy class $\cR$ on an alphabet $\cA$ with $d\geq5$ letters and we suppose that Theorem \ref{propcombinatoria} holds for any Rauzy class $\cR'$ on an alphabet $\cA'$ with $d'<d$ letters.

\subsubsection*{Pairs $(\be,\al)$ with $\al\not=A,X$ and $\be\not=Y,X$}
Consider the alphabet $\cA_{X}:=\cA\setminus\{X\}$ and the essential $\cA_{X}$-decorated Rauzy class $\cR_{X}$ which contains $\tilde{\pi}$. The only non-essential element of $\cR_{X}$ is $\tilde{\pi}$. Let $\cR^{ess}_{X}\subset\cR_{X}$ be the subset of essential elements and let $\cR^{red}_{X}$ be the associated reduced Rauzy class. The letter $Y$ is first in the bottom line and last in the top line of $\tilde{\pi}$, hence deleting $X$ from any element $\pi$ of $\cR_{X}$ we get an irreducible permutation, thus the alphabet of $\cR^{red}_{X}$ is $\cA_{X}$.

\begin{lem}\label{lemdimteoremacombinatorioriduzioneX}
If Theorem \ref{propcombinatoria} holds for $\cR^{red}_{X}$, then it holds for $\cR$ for all the pairs $(\be,\al)$ with $\al\not=A,X$ and $\be\not=Y,X$.
\end{lem}

\begin{demo}
Consider a pair $(\be,\al)$ with $\al\not=A,X$ and $\be\not=Y,X$. Observe that any element $\widehat{\pi}$ in $\cR^{red}_{X}$ satisfies $\widehat{\pi}^{t}(A)=\widehat{\pi}^{b}(Y)=1$, therefore $(\be,\al)$ satisfies the assumption in Theorem \ref{propcombinatoria} with respect to $\cR^{red}_{X}$. Applying the inductive hypothesis we get either $\widehat{\pi}$ in $\cR^{red}_{X}$ satisfying condition (\ref{eqcombinatoriasuffa}), or a pair of combinatorial data $\widehat{\pi}$ and $\widehat{\pi'}$ in $\cR^{red}_{X}$ satisfying respectively condition (\ref{eqcombinatoriasuffbuno}) and condition (\ref{eqcombinatoriasuffbdue}). In the first case, that is if there exists $\widehat{\pi}$ in $\cR^{red}_{X}$ satisfying (\ref{eqcombinatoriasuffa}), the unique essential pre-image $\pi=red^{-1}(\widehat{\pi})$ of $\widehat{\pi}$ is an element of $\cR$ satisfying condition (\ref{eqcombinatoriasuffa}). Theorem \ref{propcombinatoria} therefore holds for $(\be,\al)$.  Otherwise we have both $\widehat{\pi}$ in $\cR^{red}_{X}$ and $V$ in $\cA_{X}$ satisfying (\ref{eqcombinatoriasuffbuno}) and $\widehat{\pi'}$ in $\cR^{red}_{X}$ and $V'$ in $\cA_{X}$ satisfying (\ref{eqcombinatoriasuffbdue}). This case is more complicated to discuss.

Consider $\widehat{\pi}$ in $\cR^{red}_{X}$ and $V$ in $\cA_{X}$ satisfying (\ref{eqcombinatoriasuffbuno}). Let $\pi$ in $\cR^{ess}_{X}$ be the (unique) essential pre-image of $\widehat{\pi}$. We have $\pi^{t}(X)=1$ and we consider two cases, according to the position of $X$ in the bottom line of $\pi$.
\begin{enumerate}
\item
If $\pi^{b}(X)<\pi^{b}(\be)$ then $\pi$ satisfies (\ref{eqcombinatoriasuffbuno})
\item
If $\pi^{b}(X)>\pi^{b}(\be)$ all pairs $(\be,\al)$ with
$\al=\be$ are automatically excluded. In this case we have
$$
\pi=
\left(
\begin{array}{cccccccccc}
X & A &       & \dots &       & \al &       & \dots  &       &  V  \\
Y &   & \dots &   V   & \dots & \be & \dots &   X    & \dots & \al \\
\end{array}
\right)
\textrm{ with }
\pi^{t}(\al)=\pi^{b}(\be),
$$
therefore $\pi$ does not satisfy (\ref{eqcombinatoriasuffbuno}). Anyway applying the following combinatorial operations
$$
\pi\mapsto
\left(
\begin{array}{cccccccccc}
X & A &       & \dots &       & \al &       & \dots  &       &  V  \\
Y &   & \dots &   V   & \dots & \al & \dots &   \be  & \dots &  X  \\
\end{array}
\right)
$$
$$
\mapsto
\left(
\begin{array}{ccccccccc}
X & \dots &       &   V   &  A  &       & \dots  &       & \al \\
Y & \dots &   V   & \dots & \al & \dots &  \be   & \dots &   X \\
\end{array}
\right)
$$
$$
\mapsto
\left(
\begin{array}{ccccccccc}
X & \dots &       &   V   &  A  &       & \dots  &       &  \al \\
Y & \dots &   V   & \dots & \al & \dots &   X    & \dots & \be\\
\end{array}
\right)
$$
we get a combinatorial datum satisfying (\ref{eqcombinatoriasuffa}). (Note that the argument is compatible with cases $V\not=Y$ and $V=Y$.)
\end{enumerate}

Now consider $\widehat{\pi'}$ in $\cR^{red}_{X}$ and $V'$ in $\cA_{X}$ satisfying (\ref{eqcombinatoriasuffbdue}). Note that we have $Y\neq A,V'$, therefore the general form of $\widehat{\pi'}$ is
$$
\widehat{\pi'}=
\left(
\begin{array}{ccccccc}
A  & \dots &   V'  & \dots &  \al  & \dots & \be \\
Y  &       & \dots &  \be  & \dots &       &  V'  \\
\end{array}
\right)
\textrm{ with }
\widehat{\pi'}^{t}(\al)=\widehat{\pi}'^{b}(\be)+1.
$$
Let $\pi'$ in $\cR^{ess}_{X}$ be the unique essential pre-image of $\widehat{\pi'}$. We have $\pi'^{t}(X)=1$ and we consider separately three cases, according to the position of $X$ in the bottom line of $\pi'$.

\begin{enumerate}
\item
If $\pi'^{b}(X)<\pi'^{b}(\be)$ then $\pi'$ satisfies (\ref{eqcombinatoriasuffbdue}).

\item
If $\pi'^{b}(\be)<\pi'^{b}(X)<\pi'^{b}(\al)$ all the pairs $(\be,\al)$ with $\be=\al$ are automatically excluded. We also observe that we cannot have $\pi'^{t}(V')=\pi'^{t}(\al)-1$, since in this case $\widehat{\pi'}$ would not be admissible. Call $W\neq V'$ the letter which appears just before $\al$ in the top line. The general form of $\pi'$ therefore is
$$
\pi'=
\left(
\begin{array}{ccccccccccccccc}
X &   A   & \dots &   V'  & \dots &   W   & \al &       & \dots
&      & \be \\
Y & \dots &    W  & \dots &  \be  & \dots &  X  & \dots &  \al  & \dots &  V' \\
\end{array}
\right).
$$
with $\pi'^{t}(\al)=\pi'^{b}(\be)+2$, which does not satisfy (\ref{eqcombinatoriasuffbdue}). We apply the following Zorich steps
$$
\pi'\mapsto
\left(
\begin{array}{ccccccccccccccc}
X &   A   & \dots &   V'  & \al &       & \dots &       & \be & \dots & W \\
Y & \dots &   W   & \dots & \be & \dots &   X   & \dots & \al & \dots & V' \\
\end{array}
\right)
$$
$$
\mapsto
\left(
\begin{array}{ccccccccccccccc}
X &   A   & \dots &   V'   & \al &       & \dots &       & \be & \dots & W \\
Y & \dots &   W   & \dots  & \al & \dots &  V'   & \dots & \be & \dots & X \\
\end{array}
\right)
$$
$$
\mapsto
\left(
\begin{array}{ccccccccccccccc}
X &       & \dots &        & \be & \dots & W  &   A   & \dots &   V'  & \al \\
Y & \dots &   W   & \dots  & \al & \dots & V' & \dots &  \be  & \dots & X \\
\end{array}
\right)
$$
$$
\mapsto
\left(
\begin{array}{ccccccccccccccc}
X &       & \dots &       & \be & \dots & W &   A   & \dots &   V'  &   \al \\
Y & \dots &   W   & \dots & \al & \dots & X & \dots &  V'   & \dots &   \be \\
\end{array}
\right)
$$
and we get a combinatorial datum satisfying (\ref{eqcombinatoriasuffa}). Note that this sequence of steps is compatible with cases $A=V'$ and $A\not=V'$.

\item
If $\pi'^{b}(\al)<\pi'^{b}(X)$ both cases $\al=\be$ and $\al\neq\be$ are possible. If $\be=\al$ then
$$
\pi'=
\left(
\begin{array}{ccccc}
 X   &    A   &         &         &   \al    \\
 Y   & \dots  &   \al   &    X    &    V'    \\
\end{array}
\right)
$$
and letting $\al$ win once we get a combinatorial datum satisfying (\ref{eqcombinatoriasuffbdue}). The remaining case is $\be\not =\al$, which we separate into two sub-cases: $A=V'$ and $A\neq V'$. In the sub-case $A=V'$ the general form of $\pi'$ is
$$
\pi'=
\left(
\begin{array}{cccccccccccc}
X   &  A  &  ***  & \al &       & \dots &       & \be \\
*** & \be & \dots & \al & \dots &   X   & \dots &  A  \\
\end{array}
\right)
$$
where $***$ denotes the set of those letters $\xi$ with $\pi'^{t}(A)<\pi'^{t}(\xi)<\pi'^{t}(\al)$, which coincide with the set of $\xi$ with $\pi'^{b}(\xi)<\pi'^{b}(\be)$, and where $\pi'^{t}(\al)=\pi'^{b}(\be)+2$. We get a combinatorial datum satisfying (\ref{eqcombinatoriasuffbdue}) applying the following Zorich steps
$$
\pi'\mapsto
\left(
\begin{array}{cccccccccccc}
X   &  A  & \dots &     &       & \be &  ***  &\al \\
*** & \be & \dots & \al & \dots &  X  & \dots & A  \\
\end{array}
\right)
$$
$$
\mapsto
\left(
\begin{array}{cccccccccccc}
X   &  A  &       & \dots &       & \be &  ***  & \al \\
*** & \be & \dots &  \al  & \dots &  A  & \dots &  X \\
\end{array}
\right)
$$
$$
\mapsto
\left(
\begin{array}{cccccccccccc}
X   & *** &  \al  &  A  &       & \dots &       & \be \\
*** & \be & \dots & \al & \dots &   A   & \dots &  X   \\
\end{array}
\right).
$$
In the sub-case $A\neq V'$ the general form of $\pi'$ is
$$
\pi'=
\left(
\begin{array}{cccccccccccccc}
X &   A   & \dots &   V'  & \dots &  \al  &     & \dots &   &    &\be\\
Y & \dots &   A   & \dots &  \be  & \dots & \al & \dots & X & \dots & V' \\
\end{array}
\right),
$$
with $\pi'^{t}(\al)=\pi'^{b}(\be)+2$. Such $\pi'$ does not satisfies (\ref{eqcombinatoriasuffbdue}). We get a combinatorial datum satisfying (\ref{eqcombinatoriasuffa}) applying the following Zorich steps
$$
\pi'\mapsto
\left(
\begin{array}{cccccccccccccc}
X   &    A    &  \dots  &    V'   &         &  \dots  &         &         &    \be   &  \dots  &   \al   \\
Y   &  \dots  &    A    &  \dots  &   \be   &  \dots  &   \al   &  \dots  &    X     &  \dots  &    V'    \\
\end{array}
\right)
$$
$$
\mapsto
\left(
\begin{array}{cccccccccccccc}
X   &    A    &  \dots  &    V'   &         &  \dots  &         &         &   \be   &  \dots  &   \al   \\
Y   &  \dots  &    A    &  \dots  &   \be   &  \dots  &   \al   &  \dots  &    V'   &  \dots  &    X    \\
\end{array}
\right)
$$
$$
\mapsto
\left(
\begin{array}{cccccccccccccc}
X   &  \dots  &    V'   &         &  \dots  &         &         &   \be   &  \dots  &   \al   &    A    \\
Y   &  \dots  &    A    &  \dots  &   \be   &  \dots  &   \al   &  \dots  &    V'   &  \dots  &    X    \\
\end{array}
\right)
$$
$$
\mapsto
\left(
\begin{array}{cccccccccccccc}
X   &  \dots  &    V'   &         &  \dots  &         &         &   \be   &  \dots  &   \al   &    A    \\
Y   &  \dots  &    A    &  \dots  &   \al   &  \dots  &    V'   &  \dots  &    X    &  \dots  &   \be   \\
\end{array}
\right)
$$
$$
\mapsto
\left(
\begin{array}{cccccccccccccc}
X   &  \dots  &    V'   &         &  \dots  &         &         &   \be   &    A    &  \dots  &   \al   \\
Y   &  \dots  &    A    &  \dots  &   \al   &  \dots  &    V'   &  \dots  &    X    &  \dots  &   \be   \\
\end{array}
\right).
$$
\end{enumerate}
\end{demo}

\subsubsection*{Pairs $(\be,\al)$ with $\al\not=X,Y$ and $\be\not=B,Y$}

Consider the alphabet $\cA_{Y}:=\cA\setminus\{Y\}$ and the essential $\cA_{Y}$-decorated Rauzy class $\cR_{Y}$ which contains $\tilde{\pi}$. The only non-essential element of $\cR_{Y}$ is $\tilde{\pi}$. Let $\cR^{ess}_{Y}\subset \cR_{Y}$ be the subset of essential elements and let $\cR^{red}_{Y}$ be the associated reduced Rauzy class. Since $X$ is first in the top line and last in the bottom line of $\tilde{\pi}$, when we delete $Y$ from any $\pi$ in $\cR_{Y}$ we get an irreducible permutation, therefore the alphabet of $\cR^{red}_{Y}$ is $\cA_{Y}$.

\begin{lem}\label{lemdimteoremacombinatorioriduzioneY}
If Theorem \ref{propcombinatoria} holds for $\cR^{red}_{Y}$, then it holds for $\cR$ for all pairs $(\be,\al)$ with $\al\not=X,Y$ and $\be\not=B,Y$.
\end{lem}

\begin{demo}
Any element $\hat{\pi}$ in $\cR^{red}_{Y}$ satisfies $\widehat{\pi}^{b}(B)=\widehat{\pi}^{t}(X)=1$. The letters $X$ and $A$ play a symmetric role with respect to the letters $Y$ and $B$, therefore changing respectively $A$ with $B$, $X$ with $Y$ and the role of bottom line with the role of the top line, the same argument as in Lemma \ref{lemdimteoremacombinatorioriduzioneX} proves the lemma.
\end{demo}

\subsubsection{The pair $(B,A)$}\label{dimteoremacombinatoriocoppiaBA}

Lemma \ref{lemdimteoremacombinatoriociclistandard} proves the statement in Theorem \ref{propcombinatoria} for all pairs $(X,\al)$ with $\al\not=X$ and all pairs $(\be,Y)$ with $\be\not=Y$. Then Lemma \ref{lemdimteoremacombinatorioriduzioneX} and Lemma \ref{lemdimteoremacombinatorioriduzioneY} prove inductively the statement for all remaining pairs except for the pair $(\be,\al)=(B,A)$. To complete the induction we provide a solution for this pair.

Consider the standard element $\tilde{\pi}$ of $\cR$ as in the beginning of the proof of the theorem. Lemma 20 in \cite{konzo} (or Lemma 3.8 in \cite{av}) shows that it is possible to find a standard $\tilde{\pi}$ which is \emph{good} or \emph{degenerate}, where a standard permutation is said good if the permutation that we get deleting the letters $X$ and $Y$ from $\tilde{\pi}$ is still admissible and is said degenerate if there exists a letter $C\in\cA$ different from $X$ and $Y$ which is second or second to last in both the top and bottom lines. We consider separately the two cases.

If $\tilde{\pi}$ is good then $A\neq B$. Let $\pi$ be the element obtained from $\tilde{\pi}$ letting $Y$ win once, that is
$$
\pi=
\left(
\begin{array}{cccc}
X   &   A   & \dots &       Y   \\
Y   &   X   &   B   &    \dots   \\
\end{array}
\right)
$$
We consider the alphabet $\cA_{Y}=\cA\setminus \{Y\}$ and the $\cA_{Y}$-decorated Rauzy class $\cR_{Y}$ which contains $\pi$. We note that $\cR_{Y}$ is an essential decorated Rauzy class and we call $\cR^{red}_{Y}$ its reduction. Since $\tilde{\pi}$ is good then the alphabet of $\cR^{red}_{Y}$ is $\cA'':=\cA\setminus \{X,Y\}$. Let $\widehat{\pi}_{st}$ be a standard element in $\cR^{red}_{Y}$, that is
$$
\widehat{\pi}_{st}=
\left(
\begin{array}{ccc}
A   & \dots &   B   \\
B   & \dots &   A   \\
\end{array}
\right).
$$
An essential pre-image of $\widehat{\pi}_{st}$ in $\cR_{Y}$ is an element of the form
$$
\left(
\begin{array}{ccccc}
  X   &   A   & \dots &   Y   & \dots    \\
  Y   &   X   &   B   & \dots &   A      \\
\end{array}
\right).
$$
Letting $A$ win the proper number of times we get
$$
\left(
\begin{array}{ccccc}
  X   &   A   & \dots &       &  Y      \\
  Y   &   X   &   B   & \dots &   A      \\
\end{array}
\right),
$$
which satisfies (\ref{eqcombinatoriasuffbuno}). Since the argument is symmetric changing the top line with the bottom one, we can get also a combinatorial datum satisfying (\ref{eqcombinatoriasuffbdue}).

If $\tilde{\pi}$ is degenerate both cases $B=A$ and $B\not=A$ are possible. We consider them separately. If $A=B$ we call $W$ the last letter in the bottom line before $X$, that is the letter such that $\tilde{\pi}^{b}(W)=d-1$. We apply the sequence of Zorich steps
$$
\tilde{\pi}=
\left(
\begin{array}{cccccc}
X   &   A   & \dots &   W   & \dots &   Y   \\
Y   &   A   & \dots &       &   W   &   X   \\
\end{array}
\right)
\mapsto
\left(
\begin{array}{cccccc}
X   & \dots &   W   & \dots &   Y   &   A   \\
Y   &   A   & \dots &       &   W   &   X   \\
\end{array}
\right)
$$
$$
\mapsto
\left(
\begin{array}{cccccc}
X   & \dots &   W   & \dots &   Y   &   A   \\
Y   &   A   &   X   & \dots &       &   W   \\
\end{array}
\right)
\mapsto
\left(
\begin{array}{cccccc}
X   & \dots &   W   &   A   & \dots &   Y   \\
Y   &   A   &   X   & \dots &       &   W   \\
\end{array}
\right)
$$
$$
\mapsto
\left(
\begin{array}{cccccc}
X   & \dots &   W   &   A   & \dots &   Y   \\
Y   & \dots &       &   W   &   A   &   X   \\
\end{array}
\right)
\mapsto
\left(
\begin{array}{cccccc}
X   & \dots &   Y   & \dots &   W   &   A   \\
Y   & \dots &       &   W   &   A   &   X   \\
\end{array}
\right)
$$
and we get a combinatorial datum satisfying (\ref{eqcombinatoriasuffbdue}). Since $\tilde{\pi}$ is standard, with the symmetric argument we can get also a combinatorial datum satisfying (\ref{eqcombinatoriasuffbuno}). If $A\not=B$, since $d\geq5$, there exists a letter $C\neq X,Y,A,B$ which is second to last both in top and bottom lines and the general form of $\tilde{\pi}$ is
$$
\tilde{\pi}=
\left(
\begin{array}{ccccccc}
X   &   A   & \dots &   B   & \dots &   C   &   Y   \\
Y   &   B   & \dots &   A   & \dots &   C   &   X   \\
\end{array}
\right).
$$
We get a combinatorial datum satisfying (\ref{eqcombinatoriasuffa}) for the pair $(B,A)$ applying the following sequence of Zorich steps:
$$
\tilde{\pi}\mapsto
\left(
\begin{array}{ccccccc}
X   &   A   & \dots &   B   & \dots &   C   &   Y   \\
Y   & \dots &   C   &   X   &   B   & \dots &   A   \\
\end{array}
\right)
\mapsto
\left(
\begin{array}{ccccccc}
X   &   A   &   Y   & \dots &   B   & \dots &   C   \\
Y   & \dots &   C   &   X   &   B   & \dots &   A   \\
\end{array}
\right)
$$
$$
\mapsto
\left(
\begin{array}{ccccccc}
X   &   A   &   Y   & \dots &   B   & \dots &   C   \\
Y   & \dots &   C   &   B   & \dots &   A   &   X   \\
\end{array}
\right)
\mapsto
\left(
\begin{array}{ccccccc}
X   & \dots &   B   & \dots &   C   &   A   &   Y   \\
Y   & \dots &   C   &   B   & \dots &   A   &   X   \\
\end{array}
\right)
$$
$$
\mapsto
\left(
\begin{array}{ccccccc}
X   & \dots &   B   & \dots &   C   &   A   &   Y   \\
Y   & \dots &   A   &   X   & \dots &   C   &   B   \\
\end{array}
\right)
\mapsto
\left(
\begin{array}{ccccccc}
X   & \dots &   B   &   Y   & \dots &   C   &   A   \\
Y   & \dots &   A   &   X   & \dots &   C   &   B   \\
\end{array}
\right).
$$
\end{demo}

\emph{Note:} even for a pair $(\be,\al)$ with $\be\not=\al$ in general it is not possible to find a combinatorial datum satisfying (\ref{eqcombinatoriasuffa}). This can be seen for the pair $(B,A)$ when $\cR$ is \emph{hyperelliptic} Rauzy class.

\subsection{Main estimate}\label{teoremastimaprincipale}

Fix any $\pi\in\cR$, any letter $W\in\cA$ and any $\epsilon>0$. Consider the sub-alphabet $\cA_{W}:=\cA\setminus\{W\}$. Let $E(\pi,W,\epsilon)$ be the family of those $\cA_{W}$-colored paths $\ga$ which start at $\pi$, satisfy $q^{\ga}_{W}>1/\epsilon$ and are minimal with this property with respect to $\prec$. Consider the set
$
\De^{(1)}(E(\pi,W,\epsilon)):=
\bigsqcup_{\ga\in E(\pi,W,\epsilon)}
\De^{(1)}_{\ga}
$.

\begin{thm}\label{propstimalocalegenerale}
There exists a positive constant $C$, depending only on $\cA$, such that for any $\pi\in\cR$, any $W\in\cA$ and any $\epsilon>0$ we have
$$
\leb_{d-1}
\big(
\De^{(1)}(E(\pi,W,\epsilon))
\big)
\geq C\epsilon.
$$
\end{thm}

\subsubsection{Preliminary facts}\label{s5ss2sss1}

Denote $\widehat{\De}^{(1)}_{\pi}$ the $(d-2)$-hyperface of $\De^{(1)}_{\pi}$ whose extremal points are the vectors $e_{\xi}$ of the standard basis of $\RR^{\cA}$ with $\xi\not=W$, that is the $(d-2)$-hyperface of $\De^{(1)}_{\pi}$ opposite to the vertex $e_{W}$. Similarly for any $\cA_{W}$-colored path $\ga$ starting at $\pi$ denote $\widehat{\De}^{(1)}_{\ga}$ the $(d-2)$-hyperface of $\De^{(1)}_{\ga}$ spanned by the vectors $(1/q^{\ga}_{\xi})^{t}B_{\ga}e_{\xi}$ with $\xi\not=W$.

\begin{lem}\label{lem1s5ss2sss1}
If $\ga$ is an $\cA_{W}$-colored path then $\widehat{\De}^{(1)}_{\ga}$ is a sub-simplex of $\widehat{\De}^{(1)}_{\pi}$. Moreover we have
$$
\leb_{d-1} (\De^{(1)}_{\ga})=
\frac{1}{q^{\ga}_{W}}\leb_{d-2}
(\widehat{\De}^{(1)}_{\ga}).
$$
\end{lem}

\begin{demo}
Observe that $\spanne\{e_{\xi}\}_{\xi\not=W}$ is invariant under $^{t}B_{\ga}$ if $\ga$ is an $\cA_{W}$-colored path. Therefore $\widehat{\De}^{(1)}_{\ga}$ is a sub-simplex of $\widehat{\De}^{(1)}_{\pi}$ and the first part of the lemma is proved. Moreover $\langle ^{t}B_{\ga}e_{W},e_{W}\rangle=1$, thus the restriction of $^{t}B_{\ga}$ to $\spanne\{e_{\xi}\}_{\xi\not=W}$ preserves the $(d-1)$-volume of the subspace. Applying the same argument that we used to get equation (\ref{eqmisurasimplessi}) in \S \ref{thedistortionestimate} and recalling our normalization of $\leb_{d-2}$ on $\widehat{\De}^{(1)}_{\pi}$ we get
$$
\leb_{d-2}(\widehat{\De}^{(1)}_{\ga})=
\prod_{\xi\not=W}(q^{\ga}_{\xi})^{-1}.
$$
Comparing the last expression with (\ref{eqmisurasimplessi}) we get the second part of the lemma.
\end{demo}

Let $\cR^{col}$ be the $\cA_{W}$-decorated Rauzy class which contains $\pi$ and suppose that it is essential. In this case the reduction map $red$ can be defined on $\cR^{col}$, with image onto its reduced Rauzy class $\cR^{red}$, whose alphabet is a subset $\cA^{red}$ of $\cA_{W}$. Let $red(\pi)$ in $\cR^{red}$ be the image of $\pi$ under $red$. We recall that the reduction map $red:\cR^{col}\to \cR^{red}$ extends to a map $red:\Pi^{col}(\cR^{col})\to \Pi(\cR^{red})$. For any $\ga$ in $\Pi^{col}(\cR^{col})$ let $red(\ga)$ in $\Pi(\cR^{red})$ be its image under $red$. The formalism and notations of reduction of Rauzy classes are exposed in \S \ref{Reduction of rauzy classes} of the background. The inclusions $\cA^{red}\subset\cA_{W}\subset\cA$ induce naturally a decomposition
$$
\RR^{\cA}=
\RR^{\cA^{red}}
\oplus
\RR^{\cA_{W}\setminus\cA^{red}}
\oplus
\RR^{\cA\setminus\cA_{W}}.
$$
Consider the canonical projections $\RR^{\cA}\to\RR^{\cA^{red}}$ and $\RR^{\cA}\to\RR^{\cA_{W}\setminus\cA^{red}}$ respectively on the first and on the second factor in the splitting above.  For any $\ga$ in $\Pi^{col}(\cR^{col})$ the matrix $B_{red(\ga)}$ acts on $\RR^{\cA^{red}}$. Denote $q^{red(\ga)}:=B_{red(\ga)}\vec{1}$, where $\vec{1}$ is the column vector in $\RR^{\cA_{red}}$ with all entries equal to $1$. We remark that $q^{\ga}$ belongs to $\RR^{\cA}$ whereas $q^{red(\ga)}$ belongs to $\RR^{\cA_{red}}$.

\begin{lem}\label{lem2s5ss2sss1}
Consider any $\ga$ in $\Pi^{col}(\cR^{col})$. We have $q^{red(\ga)}_{\xi}=q^{\ga}_{\xi}$ for any $\xi\in\cA^{red}$ and $q^{\ga}_{\xi}=1$ for any $\xi\in\cA_{W}\setminus\cA^{red}$.
\end{lem}

\begin{demo}
Observe that any path $\ga$ in $\Pi^{col}(\cR^{col})$ is $(\cA_{W}\setminus\cA^{red})$-separated, thus it induces two commutative diagrams
$$
\begin{array}{ccccccc}
\RR_{+}^{\cA} & \substack{B_{\ga}\\ \longrightarrow} & \RR^{\cA}_{+}   & &
\RR_{+}^{\cA} & \substack{B_{\ga}\\ \longrightarrow} & \RR^{\cA}_{+}   \\
\downarrow    &                                  &     \downarrow  & &        \downarrow    &                                  &     \downarrow   \\
\RR_{+}^{\cA^{red}} & \substack{B_{red(\ga)}\\ \longrightarrow} & \RR_{+}^{\cA^{red}} &  &
\RR_{+}^{\cA_{W}\setminus\cA^{red}} &  \substack{Id\\ \longrightarrow}  & \RR_{+}^{\cA_{W}\setminus\cA^{red}}.\\
\end{array}
$$
The lemma therefore follows.
\end{demo}

Let $a$ be the number of letters of $\cA^{red}$. We can identify $\De^{(1)}_{red(\pi)}$ with the standard $(a-1)$-simplex in $\RR^{\cA_{red}}$. Then $\De^{(1)}_{red(\ga)}$ is an $(a-1)$-subsimplex of $\De^{(1)}_{red(\pi)}$ for any $\ga$ in $\Pi^{col}(\cR^{col})$ starting at $\pi$. Equation (\ref{eqmisurasimplessi}) in the background and Lemma \ref{lem2s5ss2sss1} implies

\begin{equation}\label{eq1s5ss2sss1}
\leb_{d-2}(\widehat{\De}^{(1)}_{\ga})=
\leb_{a-1}(\De^{(1)}_{red(\ga)}).
\end{equation}

For any $\ga$ in $\Pi^{col}(\cR^{col})$ we introduce the notation
$
\widehat{\PP}
\big(red(\ga)\big):=
\leb_{a-1}(\De^{(1)}_{red(\ga)})
$.
For a family $\Ga$ of paths $\ga$ in $\Pi^{col}(\cR^{col})$ we write
$
\widehat{\PP}
\big(red(\Ga)\big):=
\leb_{a-1}
\big(
\bigcup_{\ga\in\Ga}
\De^{(1)}_{red(\ga)}
\big)
$.
If $red(\Ga)$ is disjoint we have
$
\widehat{\PP}
\big(
red(\Ga)
\big)
:=
\sum_{red(\ga)\in red(\Ga)}
\widehat{\PP}
\big(red(\ga)\big)
\big)
$.

\subsubsection{Proof of the estimate}\label{s5ss2sss2}

Let $\cR^{col}$ be the decorated Rauzy class which contains $\pi$. We first show that if $\cR^{col}$ is not essential then the statement holds trivially.

\begin{lem}\label{lem1s5ss2sss2}
If $\cR^{col}$ is not essential then $\PP(E(\pi,W,\epsilon))>\epsilon/2$.
\end{lem}

\begin{demo}
Since $\cR^{col}$ does not contain essential elements, then $W$ is in last position either in the top row or in the bottom row of $\pi$. Suppose without loss of generality that $\pi^{b}(W)=d$ and let $\ga_{\ast}$ be the top arrow starting at $\pi$ with loser $W$. Since $W$ has to keep in last position in the bottom row of the ending point of $\ga_{\ast}$, we must have
$
\pi=
\left(
\begin{array}{ccc}
      & \dots & \ast \\
\dots &  \ast & W
\end{array}
\right)
$
and then $\ga_{\ast}$ is a length-one loop at $\pi$. Therefore $\cR^{col}=\{\pi\}$ and any $\ga$ in $\Pi^{col}(\cR^{col})$ is of the form $\ga=\ga_{\ast}^{k}$, that is it is the concatenation of $k$ copies of $\ga_{\ast}$. For such $\ga$ we have $q^{\ga}_{W}=k+1$ and $q^{\ga}_{\xi}=1$ for $\xi\not=W$. Therefore $E(\pi,W,\epsilon)=\{\ga\}$, where $\ga=\ga_{\ast}^{k}$ satisfies $k<1/\epsilon\leq k+1$. Thus $\PP(E(\pi,W,\epsilon))=\PP(\ga)=(k+1)^{-1}$ and the statement follows.
\end{demo}

According to Lemma \ref{lem1s5ss2sss2} it only remains to prove Theorem \ref{propstimalocalegenerale} when $\cR^{col}$ is essential. In this case the map $red$ is defined on $\cR^{col}$. Let $\cR^{red}:=red(\cR^{col})$ be the reduced Rauzy class and $\cA^{red}\subset\cA_{W}$ be the alphabet of $\cR^{red}$.

Any $\ga$ in $E(\pi,W,\epsilon)$ is $\cA_{W}$-colored, thus according to Lemma \ref{lem1s5ss2sss1} the $(d-2)$-hyperface $\widehat{\De}^{(1)}_{\ga}$ of $\De^{(1)}_{\ga}$ is contained in $\widehat{\De}^{(1)}_{\pi}$, which is the $(d-2)$-hyperface of $\De^{(1)}_{\pi}$ opposite to $e_{W}$. We have
$$
\PP\big(E(\pi,W,\epsilon)\big)=
\sum_{\ga\in E(\pi,W,\epsilon)}
\frac{1}{q^{\ga}_{W}}
\leb_{d-2}(\widehat{\De}^{(1)}_{\ga}),
$$
and since
$
\widehat{\De}^{(1)}_{\pi}=
\bigsqcup_{\ga\in E(\pi,W,\epsilon)}
\widehat{\De}^{(1)}_{\ga}
$
modulo a set of measure zero, then $\PP(E(\pi,W,\epsilon))$ equals the integral over $\widehat{\De}^{(1)}_{\pi}$ of the piecewise constant function whose constant value over any $\widehat{\De}^{(1)}_{\ga}$ is $1/q^{\ga}_{W}$. The idea of the proof is to show that $\leb_{d-2}$ gives small measure to the subset of $\widehat{\De}^{(1)}_{\pi}$ of those $\widehat{\De}^{(1)}_{\ga}$ such that $q^{\ga}_{W}\gg 1/\epsilon$. The main tool is the map $red$, which converts such estimate into an estimate for $\widehat{\PP}$ over $\Pi(\cR^{red})$. More precisely, equation (\ref{eq1s5ss2sss1}) implies $\leb_{d-2}(\widehat{\De}^{(1)}_{\ga})=\widehat{\PP}(red(\ga))$ for any $\cA_{W}$-colored path $\ga$, thus we have
\begin{equation}\label{eq1s5ss2sss2}
\PP(E(\pi,W,\epsilon))=
\sum_{\ga\in E(\pi,W,\epsilon)}
\frac
{\widehat{\PP}\big(red(\ga)\big)}
{q^{\ga}_{W}}.
\end{equation}

For any $\ga$ in $E(\pi,W,\epsilon)$ consider a sub-path $\nu\prec\ga$ with $q^{\nu}_{W}<1/\epsilon$. Such $\nu$ is of course $\cA_{W}$-colored, being a sub-path of $\ga$. Fix any such $\nu$ and denote $E(\pi,W,\epsilon|\nu)$ the set of those $\ga$ in $E(\pi,W,\epsilon)$ with $\nu\prec\ga$.

\begin{defi}\label{def1s5ss2sss2}
An \emph{intermediate path} is an $\cA_{W}$-colored path $\nu$ starting at $\pi$ which satisfies $q^{\nu}_{W}<1/\epsilon$ and the following extra property: if $\pi'\in\cR^{col}$ is the ending point of $\nu$ then for any $\cA_{W}$-colored path $\eta$ starting at $\pi'$ and containing at least one arrow where $W$ loses, we have $q^{\nu\eta}_{W}\geq 1/\epsilon$.
\end{defi}

Denote $I(\pi,W,\epsilon)$ the set of the intermediate paths starting at $\pi$ which are minimal with respect to the ordering $\prec$.

\begin{lem}\label{lem2s5ss2sss2}
For any $\ga\in E(\pi,W,\epsilon)$ there exists an unique path $\nu$ in $I(\pi,W,\epsilon)$ such that $\ga\in E(\pi,W,\epsilon|\nu)$. On the other hand for any $\nu \in I(\pi,W,\epsilon)$ the set $E(\pi,W,\epsilon|\nu)$ is not empty.
\end{lem}

\begin{demo}
Fix $\ga\in E(\pi,W,\epsilon)$ and decompose it as $\ga=\ga'\ga_{last}$, where $\ga_{last}$ is the last arrow of $\ga$. The loser of $\ga_{last}$ is $W$ by minimality of paths in $E(\pi,W,\epsilon)$. The path $\ga'$ is of course $\cA_{W}$-colored and satisfies $q^{\ga'}_{W}<1/\epsilon$. Call $\pi'$ the ending point of $\ga'$. Any $\cA_{W}$-colored path $\eta$ starting at $\pi'$ has $\ga_{last}$ as first arrow, since $W$ is the winner of the other arrow starting at $\pi'$. It follows that we can decompose any such $\eta$ as $\eta=\ga_{last}\eta'$. Since $B_{\eta}q^{\ga'}=B_{\eta'}q^{\ga}$, then $q^{\ga'\eta}_{W}\geq q^{\ga}_{W}>1/\epsilon$, thus $\ga'$ is intermediate. We proved that the set of intermediate paths $\ga'$ with $\ga'\prec\ga$ is not empty, thus the minimal path $\nu=\nu(\ga)$ is well defined. Uniqueness of $\nu$ follows by minimality. The second statement is evident and the lemma is proved.
\end{demo}


For any $k\in\NN$ denote $I(\pi,W,\epsilon|k)$ the set of paths $\nu$ in $I(\pi,W,\epsilon)$ such that $M(q^{\nu})\geq2^{k}/\epsilon$. We have $I(\pi,W,\epsilon)=\bigcup_{k=1}^{\infty}I(\pi,W,\epsilon|k)$ (the union is not disjoint).

\begin{lem}\label{lem3s5ss2sss2}
There exist two positive constants $C$ and $\te$, depending only on the number of intervals $d$, such that for any $k\in\NN^{\ast}$ we have
$$
\widehat{\PP}
\big(red(I(\pi,W,\epsilon|k))\big)
\leq
Ck^{\theta}2^{-(k-1)}.
$$
\end{lem}
\begin{demo}
Decompose any $\nu$ in $I(\pi,W,\epsilon)$ as $\nu=\nu'\nu_{last}$, where $\nu_{last}$ is the last arrow in $\nu$, then denote $I'(\pi,W,\epsilon)$ the family of paths $\nu'$ obtained from $\nu$ in $I(\pi,W,\epsilon)$. Since any $\nu$ in $I(\pi,W,\epsilon)$ is minimal intermediate, then paths in $I'(\pi,W,\epsilon)$ are not intermediate. Observe that for any two different paths $\nu_{1}$ and $\nu_{2}$ in $I(\pi,W,\epsilon)$ we have $\nu_{1}'\not=\nu_{2}'$, otherwise the path
$\nu':=\nu_{1}'=\nu_{2}'$ would be an intermediate element of $I'(\pi,W,\epsilon)$, which is absurd. Therefore the map $\nu\mapsto\nu'$ is injective, thus a bijection between $I(\pi,W,\epsilon)$ and $I'(\pi,W,\epsilon)$.

Fix $k\in\NN$ and denote $I'(\pi,W,\epsilon|k)$ the set of paths $\nu'$ obtained from a path $\nu$ in $I(\pi,W,\epsilon|k)$. Consider $\nu'$ in $I'(\pi,W,\epsilon|k)$ and let $\pi'$ in $\cR^{col}$ be its ending point. Since $\nu'$ is not intermediate then there exists a $\cA_{W}$-colored path $\eta'$ starting at $\pi'$ which contains one arrow with loser $W$ and such that the concatenation $\nu'\eta'$ satisfies $q^{\nu'\eta'}_{W}<1/\epsilon$. If $X\in\cA^{red}$ is the letter which wins against $W$ in $\eta'$, we obviously have $q^{\nu'}_{X}<1/\epsilon$. Since $X$ belongs to $\cA^{red}$, then Lemma \ref{lem2s5ss2sss1} implies $q^{\nu'}_{X}=q^{red(\nu')}_{X}$, where $red(\nu')$ is the reduced path of $\nu$ and $q^{red(\nu')}=B_{red(\nu')}\vec{1}$. In particular we have
$$
m(q^{red(\nu')})<1/\epsilon.
$$

On the other hand for any $\nu$ in $I(\pi,W,\epsilon|k)$ we have $M(q^{\nu})<2M(q^{\nu'})$, thus $M(q^{\nu'})>2^{k-1}/\epsilon$. Since $M(q^{\nu'})=M_{\cA^{red}}(q^{\nu'})$, applying again Lemma \ref{lem2s5ss2sss1} we have
$$
M(q^{red(\nu')})=
M(q^{\nu'})>
2^{k-1}/\epsilon.
$$
Denote $red(I'(\pi,W,\epsilon|k))$ the image of $I'(\pi,W,\epsilon|k)$ under the map $red$. We proved that for any $\nu'$ in $I'(\pi,W,\epsilon|k)$ we have
$$
M(q^{red(\nu)})\geq
2^{k-1}m(q^{red(\nu)}).
$$
According to equation (\ref{eq1teo1thedistorsionestimate}) in the background there exist two positive constants $C$ and $\te$, depending only the cardinality of $\cA^{red}$, such that for any $k\in\NN$ we have
$$
\widehat{\PP}
\big(
red(I'(\pi,W,\epsilon|k))
\big)
\leq Ck^{\theta}2^{k-1}.
$$
For any $\nu$ in $I(\pi,W,\epsilon|k)$ we trivially have $red(\nu')\prec red(\nu)$, thus
$$
\widehat{\PP}\big(red(I(\pi,W,\epsilon|k))\big)\leq
\widehat{\PP}\big(red(I'(\pi,W,\epsilon|k))\big)
$$
and the lemma is proved.
\end{demo}

Fix any $\nu$ in $I(\pi,W,\epsilon)$ and let $\pi'$ in $\cR^{col}$ be the point where $\nu$ ends. Recall that we call $E(\pi,W,\epsilon|\nu)$ the set of those $\ga$ in $E(\pi,W,\epsilon)$ with $\nu\prec\ga$. Then we denote $S(\pi,W,\epsilon|\nu)$ the set of paths $\eta$ in $\Pi(\cR^{col})$ starting at $\pi'$ and such that the concatenation $\ga=\nu\eta$ belongs to $E(\pi,W,\epsilon|\nu)$. For such $\ga$ we have $q^{\ga}=B_{\eta}q^{\nu}$. Now fix any integer $m\geq 1$ and consider the set $S(\pi,W,\epsilon|\nu,m)$ of those paths $\eta$ in $S(\pi,W,\epsilon|\nu)$ such that the concatenation $\nu\eta$ satisfies $q^{\nu\eta}_{W}\geq2^{m}M(q^{\nu})$. Let $red\big(S(\pi,W,\epsilon|\nu,m)\big)$ be the image of $S(\pi,W,\epsilon|\nu,m)$ under the map $red$. It is a family of paths in $\Pi(\cR^{red})$ starting at the element $red(\pi')$ of $\cR^{red}$ where $red(\nu)$ ends.

\begin{lem}\label{lem4s5ss2sss2}
There exist two positive constant $C$ and $\te$, depending only on the cardinality $\cA$, such that for any $\nu$ in $I(\pi,W,\epsilon)$ and any integer $m\geq 1$ we have
$$
\widehat{\PP}_{red(\nu)}
\big(red(S(\pi,W,\epsilon|\nu,m))\big)
\leq Cm^{\theta}2^{-(m-1)}.
$$
\end{lem}

\begin{demo}
Fix $\nu\in I(\pi,W,\epsilon)$ and a positive integer $m$. Consider any $\eta\in S(\pi,W,\epsilon|\nu,m)$ and decompose it as $\eta=\eta'\eta_{last}$, where $\eta_{last}$ is its last arrow. Since the concatenation $\ga=\nu\eta$ belongs to $E(\pi,W,\epsilon)$ then the loser of $\eta_{last}$ is $W$. Moreover $\eta'$ is $\{W\}$-separated, since $\nu$ is minimal intermediate. Let $Y\in\cA^{red}$ be the letter which wins against $W$ in the arrow $\eta_{last}$. Since $\eta'$ is $\{W\}$-separated and obviously $q^{\nu}_{W}\leq M(q^{\nu})$, then $(B_{\eta'}q^{\nu})_{Y}\geq (2^{m}-1)M(q^{\nu})\geq 2^{m-1}M(q^{\nu})$, that is
$$
M(B_{\eta'}q^{\nu})\geq 2^{m-1}M(q^{\nu}).
$$
Call $\pi'$ the ending point of $\nu$ and for any $M\geq m-1$ let $\Ga_{M}$ be the set of $\{W\}$-separated paths $\eta'$ in $\Pi(\cR^{col})$ starting at $\pi'$ such that
$$
2^{M}M(q^{\nu})
\leq M(B_{\eta'}q^{\nu})
<2^{M+1}M(q^{\nu}).
$$
Denote $\widehat{\Ga}_{M}:=red(\Ga_{M})$. Since any $\eta'\in \Ga_{M}$ is $\{W\}$-separated, remark \ref{rem1standarddecompositionofseparatedpaths} in \S \ref{standarddecompositionofseparatedpaths} implies that there exists a positive integer $s$ with $s\leq2(d-1)$ and $s$ paths $\widehat{\eta}_{1},\dots,\widehat{\eta}_{s}$ in $\Pi(\cR^{red})$, not complete with respect to the alphabet $\cA^{red}$ and such that
\begin{equation}\label{eq2lem3propstimalocalegenerale}
red(\eta')=\widehat{\eta}_{1}\dots\widehat{\eta}_{s}.
\end{equation}
We put $\widehat{q}^{(0)}:=q^{red(\nu)}$ and $\widehat{\eta}^{0}:=red(\nu)$ and for any $i=1,\dots,s$ we define inductively $\widehat{\eta}^{i}:=\widehat{\eta}^{i-1}\widehat{\eta}_{i}$ and $\widehat{q}^{(i)}:=B_{\widehat{\eta}_{i}}\widehat{q}^{(i-1)}$. We can find $s$ non-negative integers $m_{1},\dots,m_{s}$ such that for any $i\in\{1,\dots,s\}$ we have:
\begin{equation}\label{eq3lem3propstimalocalegenerale}
2^{m_{i}}M(\widehat{q}^{(i-1)})
\leq M(\widehat{q}^{(i)})
<2^{m_{i}+1}M(\widehat{q}^{(i-1)}).
\end{equation}
Moreover $m_{1},\dots,m_{s}$ satisfy the relation:
\begin{equation}\label{eq4lem3propstimalocalegenerale}
M-s-1\leq m_{1}+\cdots +m_{s}\leq M.
\end{equation}
Fix a positive integer $s$ with $s\leq2(d-1)$ and $s$ non-negative integers $m_{1},\dots,m_{s}$ satisfying (\ref{eq4lem3propstimalocalegenerale}).
For any $i\in\{1,\dots,s\}$ define the set $\widehat{\Ga}(m_{1},\dots,m_{i})$ of those $\widehat{\eta}^{i}$ which satisfy the first $i$ conditions in equation (\ref{eq3lem3propstimalocalegenerale}) for the first $i$ integers $m_{1},\dots,m_{i}$. Fix $i\in\{1,\dots,s-1\}$ and $\widehat{\eta}^{i}\in \widehat{\Ga}(m_{1},\dots,m_{i})$. Equation (\ref{eq2teo1thedistorsionestimate}) in \S \ref{thedistortionestimate} implies that there exist two positive constants $C$ and $\te$, depending only on the cardinality of $\cA^{red}$, such that
$$
\widehat{\PP}_{\widehat{\eta}^{i}}
\{\widehat{\eta}_{i+1}\textrm{ not complete };
M(\widehat{q}^{(i+1)})\geq2^{m_{i}}M(\widehat{q}^{(i)})\}
\leq C(m_{i+1}+1)^{\theta}2^{-m_{i+1}}.
$$
Applying this last equation $s$ times we get
$$
\widehat{\PP}_{red(\nu)}
(\widehat{\Ga}(m_{1},\dots,m_{s}))\leq \prod_{i=1}^{s}C(m_{i}+1)^{\theta}2^{-m_{i}}\leq C^{s}M^{s\theta}2^{-M+s+1}.
$$
For any $s$ with $s\leq2(d-1)$ the number of vectors $(m_{1},\dots,m_{s})\in \NN^{s}$ satisfying condition (\ref{eq4lem3propstimalocalegenerale}) is proportional to $M^{s-1}$, therefore summing over all admissible $(m_{1},\dots,m_{s})$ and all $s=1,\dots,2(d-1)$, modulo changing the constants $C$ and $\theta$, we get
$$
\widehat{\PP}_{red(\nu)}(\widehat{\Ga}_{M})\leq C(M+1)^{\theta}2^{-M}.
$$
Since $\{red(\eta');\eta\in S(\pi,W,\epsilon|\nu,m)\}$ is contained in $\bigcup_{M\geq m-1}\widehat{\Ga}_{M}$, summing over $M\geq m-1$ we get
$$
\widehat{\PP}_{red(\nu)} \{\widehat{\eta}';\eta \in S(\pi,W,\epsilon|\nu,m)\} \leq Cm^{\theta}2^{-(m-1)}
$$
and it follows trivially $P_{\widehat{\nu}} (\widehat{S}(\pi,W,\epsilon|\nu,m))\leq Cm^{\theta}2^{-(m-1)}$. The lemma is proved.
\end{demo}

Here we finish the proof of Theorem \ref{propstimalocalegenerale}. Let $C$ and $\te$ be the constant appearing in Lemma \ref{lem3s5ss2sss2} and Lemma \ref{lem4s5ss2sss2}. Take $N\in\NN$ such that $CN^{\theta}2^{N-1}\ll1$ and we set $c:=1-CN^{\theta}2^{N-1}$. For convenience of notation we define $I_{N}:=I(\pi,W,\epsilon)\setminus I(\pi,W,\epsilon|N)$. Similarly, for any $\nu$ in $I(\pi,W,\epsilon)$ we define $S_{N}(\nu):=S(\pi,W,\epsilon|\nu)\setminus S(\pi,W,\epsilon|\nu,N)$. Lemma \ref{lem3s5ss2sss2} implies $\widehat{\PP}\big(red(I_{N})\big)\geq c$. For any $\nu\in I(\pi,W,\epsilon)$ Lemma \ref{lem4s5ss2sss2} implies
$\widehat{\PP}_{red(\nu)}\big(red(S_{N}(\nu))\big)\geq c$.

For any $\nu$ in $I(\pi,W,\epsilon)$ and any $\eta$ in $S(\pi,W,\epsilon|\nu)$ the concatenation $\ga=\nu\eta$ satisfies $red(\ga)=red(\nu)red(\eta)$, thus
$\widehat{\PP}\big(red(\ga)\big)=
\widehat{\PP}\big(red(\nu)\big)
\widehat{\PP}_{red(\nu)}\big(red(\eta)\big)
$. Moreover for $\eta$ in $S_{N}(\nu)$ the concatenation $\ga=\nu\eta$ satisfies $q^{\ga}_{W}\leq2^{N}M(q^{\nu})$. Finally the inclusion $S_{N}(\nu)\subset S(\pi,W,\epsilon|\nu)$ implies trivially
$
\widehat{\PP}
\big(red(S_{N}(\nu))\big)
\leq
\widehat{\PP}
\big(red(S(\pi,W,\epsilon|\nu))\big)
$,
thus we have
$$
\sum_{\ga\in E(\pi,W,\epsilon|\nu)}
\frac
{\widehat{\PP}\big(red(\ga)\big)}
{q^{\ga}_{W}}
\geq
\frac
{\widehat{\PP}\big(red(\nu)\big)}
{2^{N}M(q^{\nu})}
\widehat{\PP}_{red(\nu)}\big(red(S_{N}(\nu))\big)
\geq
\frac
{c}
{2^{N}}
\frac
{\widehat{\PP}\big(red(\nu)\big)}
{M(q^{\nu})}.
$$
On the other hand Equation (\ref{eq1s5ss2sss1}) and Lemma \ref{lem2s5ss2sss1} imply
$$
\PP\big(E(\pi,W,\epsilon)\big)
=
\sum_{\nu\in I(\pi,W,\epsilon)}
\sum_{\ga\in E(\pi,W,\epsilon)}
\frac
{\widehat{\PP}\big(red(\ga)\big)}
{q^{\ga}_{W}},
$$
thus we have
$$
\PP\big(E(\pi,W,\epsilon)\big)
\geq
\frac
{c}
{2^{N}}
\sum_{\nu\in I(\pi,W,\epsilon)}
\frac
{\widehat{\PP}\big(red(\nu)\big)}
{M(q^{\nu})}
\geq
\frac
{c}
{2^{N}}
\frac
{\widehat{\PP}\big(red(I_{N})\big)}
{2^{N}}
\epsilon
\geq
\big(\frac{c}{2^{N}}\big)^{2}
\epsilon,
$$
where the second inequality holds because $M(q^{\nu})\leq 2^{N}/\epsilon$ for $\nu$ in $I_{N}$ and the third inequality follows because the inclusion $I_{N}\subset I(\pi,W,\epsilon)$ implies trivially
$
\widehat{\PP}\big(red(I_{N})\big)
\leq
\widehat{\PP}\big(red(I(\pi,W,\epsilon))\big).
$
Theorem \ref{propstimalocalegenerale} is proved.

\end{document}